\documentclass[10pt,a4paper]{article} 
\usepackage{amsmath,amssymb}
\usepackage[dvips]{epsfig}
\usepackage[francais]{babel}

\newtheorem{theo}{Th\'eor\`eme}
\newtheorem{prop}{Proposition}
\newtheorem{lemme}{Lemme}
\newtheorem{coro}{Corollaire}
\newcommand {\csq}{\noindent {\bf Cons\'equence :} \rm }
\newcommand {\rem}{\noindent {\bf Remarque :} \rm }

\newcommand {\dem}{\noindent {\sc Preuve : } \rm
}\newcommand{\findem}{\hfill$\square$\par\medskip}

\def\mddefault{}
\def\updefault{}
\def\lims{\,\overline{\lim}\;}
\def\limi{\,\underline{\lim}\;}

\def\M#1{{\mathbb #1}}

\def\R{\M{R}}
\def\T{\M{T}}
\def\t{\M{T}}
\def\Z{\M{Z}}
\def\rond{{\scriptstyle\circ}}
\def\e#1{{\hbox{e}^{#1}}}
\def\unde#1{{\hbox{ \bf 1}_{#1}}}
\def\DSM{\displaystyle}
\def\tend#1{\mathop{\longrightarrow}\limits_{#1}}
\def\inv#1{{\frac{1}{#1}}}

\def\C#1{{\cal #1}}
\def\<{{\hspace{-.5pt} <\hspace{-.5pt}}}
\def\>{{\hspace{-.5pt} >\hspace{-.5pt}}} 
\def\I{{\cal I}}
\def\J{{\cal J}}
\def\O{{\cal O}}
\def\Oun{\C{O}_1}
\def\scal(#1,#2){\langle#1,#2\rangle}
\def\vect#1{\overrightarrow{#1}}
\def\mat(#1,#2){ \begin{smallmatrix}#1\\#2\end{smallmatrix}}
\def\pmat(#1,#2){
\left(\begin{smallmatrix}#1\\#2\end{smallmatrix}\right)}
\def\Pmat(#1,#2){ \begin{pmatrix}#1\\#2\end{pmatrix}}
\def\Det#1{\left|{#1} \right|}

\title{\Large Valeurs propres de transformations li\'ees aux
rotations 
irrationnelles et aux fonctions en escalier.\\
Eigenvalues of transformations arising from irrational rotations and step functions.}
\author{M\'elanie Guenais et Fran\c cois Parreau}

\begin{document}

\maketitle
\begin{center}
{\it \small Laboratoire de Math\'ematiques, \'equipe de Topologie et
Dynamique,
 UMR 8628 du CNRS,\\
 Universit\'e PARIS-SUD, 91400 Orsay.
{\small e-mail:} Melanie.Guenais@math.u-psud.fr\\
Laboratoire d'Analyse,
 G\'eom\'etrie et
Applications, UMR 7539 du CNRS,\\
 Universit\'e PARIS 13, 93430 Villetaneuse.
{\small e-mail:} parreau@math.univ-paris13.fr}
\end{center}

\vspace{1cm}

{\bf Abstract. }\\
Given an irrational rotation $T$ on $\M T$ we settle necessary and
sufficient conditions on a step function $\phi$ and $t\in \M T$ for
the existence of measurable solutions to the cohomogical equation

$$\exp{(2i\pi\phi)}=\e{2i\pi t}f/f\rond T.$$

This yields a characterization of eigenvalues and eigenfunctions for
several transformations arising from irrational rotations and step
functions~: cylinder flows, special flows, induced maps...

From there we give constructions of special flows and three-interval
exchange transformations with unusual spectral properties.  In both
cases we exhibit examples with Kronecker factors of infinite rank.  We
also construct three-interval exchange transformations which are
non-trivially conjugate to irrational rotations or to odometers.
Similarly there exist special flows over irrational rotations which
are non-trivially conjugate to translations flows on $\M T^2$ or on
solenoids.

Finally, we prove a regularization property which allows us to give
similar examples of special flows with smooth ceiling functions, under
natural Diophantine conditions for the rotation.

\smallskip

{\bf R\'esum\'e. }\\
Soit $T$ une rotation irrationnelle de $\M T$.  Nous \'etablissons
des  conditions n\'ecessaires et suffisantes sur une
fonction en escalier $\phi$ et $t\in \M T$  pour que l'\'equation

$$\exp{(2i\pi\phi)}=\e{2i\pi t}f/ f\rond T$$

admette une solution mesurable.

Ceci donne une caract\'erisation des valeurs propres et fonctions
propres de syst\`emes dynamiques tels que les \'echanges de trois
intervalles et certains flots sp\'eciaux.  \`A partir de l\`a nous en
donnons diff\'erentes constructions avec des propri\'et\'es tr\`es
particuli\`eres.  Nous montrons d'abord que ces syst\`emes peuvent
admettre un facteur de Kronecker de rang infini.  Nous contruisons
ensuite des \'echanges de trois intervalles conjugu\'es non
trivialement \`a une rotation irrationnelle ou \`a un odom\`etre, et
des flots sp\'eciaux conjugu\'es non trivialement \`a un flot de
translations sur $\M T^2$ ou sur un sol\'eno\"\i de.

Enfin nous montrons une propri\'et\'e de r\'egularisation qui permet
d'obtenir les m\^emes exemples de flots sp\'eciaux avec des fonctions
plafond r\'eguli\`eres, sous des conditions naturelles d'approximation
diophantienne.

\tableofcontents

\part{Pr\'esentation g\'en\'erale}

\section{Introduction.}
\subsection{L'\'equation de cohomologie.}

On note $\M T=\M R/\M Z$ et $\lambda$ la mesure de Lebesgue sur $\M
T$.  Le point de d\'epart de ce travail est une caract\'erisation,
\'etant donn\'e un irrationnel $\alpha$ de $\M T$, des triplets
$(\beta, s,t) \in [0,1[{}\times \M R\times \M T$ pour lesquels il
existe une fonction mesurable $f$ non nulle satisfaisant pour
presque tout $x \in \M T$, l'\'egalit\'e
\begin{equation}\label{eq}
\exp{(2i\pi s \unde{[0,\beta[}(x))}=\e{2i\pi t}\frac{f(x)} 
{f(x+\alpha)}
\end{equation}
Ce probl\`eme s'inscrit dans l'\'etude plus g\'en\'erale des
\'equations de cohomologie
\begin{equation}\label{eq'}
\e{2i\pi\phi}=\e{2i\pi t}\frac{f}{f\rond T}.
\end{equation}
Ici $T$ d\'esignera toujours la translation de $\alpha$ sur $\M T$ et
$\phi$ une fonction r\'eelle mesurable sur $\M T$. On appelle
quasi-cobords les fonctions $\e{2i\pi \phi}$ cohomologues \`a une
constante, c'est-\`a-dire pour lesquelles l'\'equation (\ref{eq'})
admet des solutions $(f, t)$ avec $f$ mesurable non nulle et $t$
constante, et lorsque $t=0$ on parle de cobords. Comme $T$ est
ergodique, toute fonction $f$ satisfaisant (\ref{eq'}) est de module
constant, et par cons\'equent il suffit d'\'etudier l'existence de
solutions avec des fonctions $f$ de module 1.

Cette \'equation est l'\'equation aux valeurs propres de l'op\'erateur
de $L^2(\M T)$ d\'efini par $V_\phi(f)= \e{2i\pi\phi} f\rond T$.  Elle
intervient naturellement dans l'\'etude de l'ergodicit\'e et des
valeurs propres des extensions de la rotation par un groupe ab\'elien
m\'etrisable compact (cf \cite{anzai}, \cite{martine}).
\medskip

L'\'equation (\ref{eq}) appara\^ \i t pour la premi\`ere fois en
th\'eorie ergodique en 1967 dans \cite{katste0} (ou \cite{katste}),
avec $s=\frac{1}{2}$~: il s'agit de montrer que l'extension \`a 2
points associ\'ee \`a $\unde{[0,\beta[}$ au dessus d'une rotation
irrationnelle d'angle $\alpha$ (d\'efinie pour $(x,y) \in \M T \times
\M Z/2\M Z$ par $S_\beta(x,y)=(Tx,\; y+\unde{[0,\beta[}(x)\mod 2)$)
est faiblement m\'elangeante (et en particulier ergodique) pour
presque tout $\beta $.

Une \'etude plus approfondie est faite par W. Veech dans \cite{vee} en
1969~: il \'enonce des conditions n\'ecessaires et des conditions
suffisantes pour l'existence d'une solution de l'\'equation (\ref{eq})
lorsque $s=1/2$.  Il prouve l'ergodicit\'e ou la non-ergodicit\'e de
$S_\beta$ sous des conditions arithm\'etiques sur $\alpha$ et $\beta$,
sans toutefois donner de condition n\'ecessaire et suffisante.  Ces
r\'esultats ont \'et\'e \'etendus \`a tout $s \in ]0,1[$ par M.
Stewart en 1981 (\cite{ste}), puis retrouv\'es, avec des arguments
plus simples, par K. Merrill en 1985 (\cite{mer})~; elle \'etablit
\'egalement une condition suffisante pour l'\'equation (\ref{eq'})
lorsque $\phi$ est une fonction en escalier avec 3 sauts.
\medskip

Nous donnons ici une solution compl\`ete \`a ce probl\`eme.  Fixons
d'abord quelques notations~:

\begin{itemize}
\item Pour un irrationnel $\alpha$ de $]0,1[$, on note $(a_n)_{n\geq
1}$ les quotients partiels et $(p_n/q_n)_{n\geq 0}$ la suite des
r\'eduites de sa fraction continue ($p_0=0$, $q_0=1$).

\item Pour tout r\'eel $x$, on d\'esigne par $\|x\|$ sa distance \`a
l'entier le plus proche $[x]$, et on pose $\<x\>=x-[x]$.  On notera
aussi $\lfloor x\rfloor$ la partie enti\`ere et $\{x\}$ la partie
fractionnaire usuelles de $x$.

\item Pour $\beta \in [0,1[$,  on note $\phi_\beta$ la fonction de $\M
T$ dans $\M R$ d\'efinie par 
$$\phi_\beta= \unde{[0,\beta[}-\beta.$$ 
\end{itemize}
On a aussi $\phi_\beta(x)=\{x-\beta\}-\{x\}$ pour tout $x\in \M T$, et
il est naturel d'\'etendre la d\'efinition de $\phi_\beta$ \`a tout
$\beta\in \M R$ (ou tout $\beta\in \M T$), en posant
$\phi_\beta=\phi_{\{\beta\}}$.  Il sera plus commode de travailler
avec $\phi_{\beta}$ qu'avec $\unde{[0,\beta[}$, ce qui revient \`a
changer $t$ en $t+s\beta$ dans l'\'equation.

\begin{theo}\label{marche}
Soient $T$ la translation par un irrationnel $\alpha$ sur $\M T$, $s$
un r\'eel non entier, $\beta$ et $t\in \M T$.  Il existe $f$ une
fonction mesurable de module 1 sur $\M T$ v\'erifiant
$\lambda$-presque partout
\begin{equation}\label{eqbis}
\exp{(2i\pi s \phi_\beta)}=\e{2i\pi t}f/ f\circ T,
\end{equation}
si et seulement si $\beta$, $s$ et $t$ satisfont les conditions
suivantes~:
il existe une suite d'entiers $(b_n)$ et un entier $k'$ tels que
$$
\beta = \sum_{n\geq 0}b_nq_n \alpha \mod 1
\quad\hbox{et}\quad
t = k'\alpha - 
\sum_{n\geq 0}[b_ns] q_n  \alpha \mod 1,$$
avec
$$\sum_{n\geq 0} \frac{|b_n|}{a_{n+1}}<\infty \qquad\hbox{et}\qquad
\sum_{n\geq 0}\|b_ns\|^2<\infty.
$$
\end{theo}

\noindent{\bf Remarque~:} 
Il est bien connu (voir \cite{sos} par exemple) que tout $\beta \in \M
T$ s'\'ecrit de fa\c con unique sous la forme $\sum_0^\infty b_n
q_n\alpha$ avec $b_n \in \{0,..,a_{n+1}\}$ avec $b_{n+1}=0$ lorsque
$b_n=a_{n+1}$~: on parle alors de d\'ecomposition d'Ostrowski.  Il
est plus naturel dans notre cas de choisir une d\'ecomposition
sym\'etrique (cf \cite{vee}) et la condition $\sum
|b_n|/a_{n+1}<\infty$ garantit l'unicit\'e \`a un nombre fini de
termes pr\`es.  L'ensemble des $\beta$ v\'erifiant la premi\`ere
condition du th\'eor\`eme~\ref{marche} d\'efinit un sous-groupe qu'on
notera $H_1(\alpha)$.  Ce groupe peut aussi se d\'efinir au moyen
d'une propri\'et\'e d'approximation par les rationnels de la forme
$k_n/q_n$.  Notons plus g\'en\'eralement pour tout $\gamma >0$,
$$
H_\gamma(\alpha)=\bigl\{
\sum_{n\geq 0}b_n q_n\alpha \mod 1 ,\quad (b_n)_n \in \M Z^{\M N}
\hbox{ et } \sum_{n\geq 0}(|b_n|/a_{n+1})^\gamma <+\infty
\bigr\},
$$
et 
$ H_\infty(\alpha)=
\{\sum_{n\geq 0} b_nq_n\alpha \mod 1,\; |b_n|/a_{n+1}\rightarrow 0\}$.
Nous montrerons  dans l'appendice la
propri\'et\'e suivante~:
\begin{prop}\label{hsat}
Soit $\alpha$ un irrationnel. On a pour tout $\gamma>0$
\begin{equation*}
H_\gamma(\alpha)=\{x \in \M T, \; \sum_n \|q_nx\|^\gamma <+\infty\}
\; \hbox{}
\end{equation*}
et
\begin{equation*}
H_\infty(\alpha)=\{x \in \M T, \; \|q_nx\|\tend{n \rightarrow
\infty}0\}.
\end{equation*}
\end{prop}
Par abus de langage, nous noterons de la m\^eme fa\c con les relev\'es
de ces groupes dans $\M R$.  Pour l'\'equation (\ref{eq}), o\`u
$\beta$ est un r\'eel de $[0,1[$ (et non un r\'eel modulo 1), $s$
\'etant un r\'eel non entier donn\'e, il existe une solution si et
seulement si $\beta$ et $t+s\beta$ satisfont les conditions du
th\'eor\`eme~\ref{marche}.  Lorsque $s=1/2$, ou plus g\'en\'eralement
lorsque $s$ est rationnel, la derni\`ere condition signifie que $sb_n$
doit \^etre entier pour tout $n$ assez grand.  On peut alors
compl\'eter les r\'esultats de \cite{ste} et \cite{vee} sur les
extensions finies associ\'ees  \`a l'indicatrice d'un intervalle.

Rappelons que si $q$ est un entier${}>1$, une fonction mesurable
$\phi$ \`a valeurs enti\`eres d\'efinit une extension \`a $q$ points
de $T$ sur $\M T \times \M Z/q\M Z$ par
$T_{\phi,q}(x,y)=(Tx,\;y+\phi(x)\mod q)$, dont l'\'etude spectrale se
ram\`ene \`a celle des op\'erateurs $V_{k\phi/q}$ pour $0\leq k <q$.
Elle a un spectre purement discret si et seulement si
$\e{2i\pi\phi/q}$ est un quasi-cobord et elle est ergodique si et
seulement si $\e{2i\pi k\phi/q}$ n'est un cobord pour aucun entier $k$
non multiple de $q$.  De l\`a, on obtient le corollaire suivant du
th\'eor\`eme~\ref{marche}~ :

\begin{coro} Soient $T$ la translation par un irrationnel $\alpha$ sur
$\M T$, $\beta\in[0,1[$ et $q$ un entier${}>1$.  L'extension \`a $q$
points de $T$ d\'efinie par $\unde{[0,\beta[}$ admet un type spectral
discret si et seulement si $\beta\in q H_1(\alpha)+\M Z+ \M
Z\alpha$.  Cette extension est ergodique si et seulement s'il
n'existe pas de diviseur $d>1$ de $q$ tel que $\beta\in
dH_1(\alpha)$.
\end{coro}

De fa\c{c}on plus g\'en\'erale, \'etant donn\'ee une fonction
mesurable r\'eelle $\phi$, le groupe $$\Sigma_{\phi}=\{s \in \M R, \;
\e{2i\pi s\phi} \hbox{ est un quasi-cobord}\}$$
joue un r\^ole dans plusieurs questions de th\'eorie ergodique (voir
par exemple \cite{mosch} \cite{lemles}, \cite{lepaII}). Il est li\'e
au groupe des valeurs propres $e(T_\phi)$ du flot cylindrique d\'efini
sur $\M T\times \M R$ par $T_\phi (x,y)=(Tx, y+\phi(x))$~: si
$s\in\Sigma_\phi$, alors les constantes $t$ correspondantes
(d\'efinies modulo $\M Z\alpha$) sont des valeurs propres de $T_\phi$
et dans le cas o\`u $T_\phi$ est ergodique on obtient ainsi toutes les
valeurs propres de $T_\phi$ (voir \cite{lemles}).

Dans le cas o\`u $\phi=\phi_\beta$, le th\'eor\`eme~\ref{marche} donne
une caract\'erisation du groupe $\Sigma_{\phi}$~: si $\beta\notin
H_1(\alpha)$ il est r\'eduit aux entiers et, si $\beta=\sum_{n\geq
0}b_nq_n\alpha\in H_1(\alpha)$,
$$\Sigma_{\phi_\beta}=
\bigl\{s \in \M R, \;\sum_{n\geq 0}\|b_ns\|^2<\infty\bigr\}.$$
C'est un groupe de type $H_2$ au sens
de \cite{hmp}.  D'autre part, I. Oren a montr\'e dans \cite{oren}
l'ergodicit\'e du flot cylindrique d\`es que $\beta \notin \M Z\alpha
\mod 1$ et on obtient alors aussi une description du groupe
$e(T_{\phi_\beta})$.

Ceci fournit en particulier des exemples pour lequels $e(T_\phi)$ est non
d\'e\-nom\-bra\-ble, permettant ainsi une construction de \cite{lepaII}~:
 
\begin{coro}
Pour tout irrationnnel $\alpha$ \`a quotients partiels non born\'es,
il existe des fonctions mesurables r\'eelles $\phi$ pour
lesquelles le flot cylindrique $T_\phi$ associ\'e est ergodique 
et admet un groupe de valeurs propres non d\'enombrable.
\end{coro}
\dem
D'apr\`es ce qui pr\'ec\`ede, il suffit de v\'erifier qu'il existe
$\beta\in H_1(\alpha)\setminus \M Z\alpha$ (dans $\M T$) tel que
$\Sigma_{\phi_\beta}$ ne soit pas d\'enombrable.  Un r\'esultat de
\cite{par} (voir aussi \cite{hmp}) montre que le groupe $\{s\in \M R,
\sum_n \|b_ns\|^2<\infty\}$ d\'efini par une suite d'entiers non nuls
$(b_n)$ est non d\'enombrable lorsque cette suite satisfait la
condition de lacunarit\'e $\sum (b_n/b_{n+1})^2<\infty$. Cela reste 
vrai si les entiers $b_n$ sont nuls sauf pour une sous-suite
v\'erifiant cette condition. Quitte \`a choisir une sous-suite
$(n_j)$ tels que les rapports $a_{n_{j+1}+1}/a_{n_j+1}$ des quotients 
partiels soient suffisamment grands, il est clair qu'on peut choisir
des coefficients $b_{n_j}$ non nuls avec $\sum
(b_{n_j}/b_{n_{j+1}})^2<\infty$ et $\sum |b_{n_j}|/a_{n_j+1}<\infty$.
En posant $\beta= \sum b_{n_j}q_{n_j}\alpha$, on obtient l'exemple
cherch\'e (le fait qu'une infinit\'e de termes soient non nuls assure 
que $\beta\notin \M Z\alpha$).
\findem

\bigskip

Nous \'etendons ensuite le th\'eor\`eme~\ref{marche} \`a toutes les
fonctions en escalier.  On peut bien s\^ur se ramener au cas d'une
fonction $\phi$ d'int\'egrale nulle.  Par coh\'erence avec le cas de
$\phi_\beta$, on notera $-s$ le saut de $\phi$ en un point de
discontinuit\'e $\beta$.

\begin{theo}\label{cn}
Soient $T$ la translation par un irrationnel $\alpha$ sur $\M T$ et
$\phi$ une fonction r\'eelle en
escalier sur $\M T$ d'int\'egrale nulle, de sauts $-s_j$ aux
points distincts $\beta_j$ ($m\geq 1$, $0\leq j\leq m$), et soit $t\in\M T$.
\begin{itemize}
\item[$\bullet$]  
On suppose qu'il existe une partition $\cal P$ de $\{0,..,m\}$ telle
qu'on ait pour tout $J\in \cal P$ et  $\beta_J\in \{\beta_j,j \in J\}$~:
\begin{enumerate}
\item  la somme $\sum_{j\in J}s_j$ est enti\`ere, 
\item  pour tout $j \in J$ il existe une suite d'entiers $(b_n^j)_n$ telle que 
$$\beta_j= \beta_{J}+ \sum_{n\geq 0} b_n^jq_n\alpha \mod 1$$
$$\hbox{avec}\qquad\sum_{n\geq 0} \frac{|b_n^j|}{a_{n+1}}<+\infty\qquad
\hbox{et}\qquad
\sum_{n\geq 0}\Bigl\|\sum_{j \in J}b_n^js_j\Bigr\|^2 <+\infty\;, 
$$
\item  il existe un entier $k'$ tel que $t=k'\alpha -\sum_{J\in \C
P}t_J$ o\`u 
$$t_J=\beta_J\sum_{j\in J}s_j
+\sum_{n\geq 0}\Bigl[\sum_{j\in J}b_n^js_j\Bigr]q_n\alpha \mod 1.
$$
\end{enumerate}
Alors il existe une fonction mesurable $f$ de module 1 satisfaisant
l'\'equation
$$\e{2i\pi \phi}= \e{2i\pi t}\frac{f}{f\rond T}\;\cdotp$$
\item[$\bullet$]
 R\'eciproquement lorsque $\sum_{j \in J} s_j
\notin \M Z$ pour toute partie stricte non vide $J$ de $\{0,..,m\}$,
ces conditions sont n\'ecessaires pour que l'\'equation admette une
solution.
\end{itemize}
\end{theo}
\rem
Ce th\'eor\`eme  montre qu'on peut trouver des fonctions en escalier
qui sont des quasi-cobords, mais qui ne se d\'ecomposent pas en somme
de quasi-cobords plus simples.  Il s'applique en particulier \`a
toutes les fonctions en escalier avec 3 sauts, ce qui am\'eliore
consid\'erablement les conditions propos\'ees par K. Merrill dans
\cite{mer}.

\subsection{Application aux flots sp\'eciaux.}
Si $\phi$ est une fonction mesurable strictement positive sur $\M T$,
on note $\tau_{\alpha,\phi}$ le flot sp\'ecial de fonction plafond
$\phi$ au dessus de la translation irrationnelle $T$ par $\alpha$
(\cite{cofosi},\cite{nad1}).  On rappelle qu'on peut le d\'efinir sur
le quotient de $\M T\times \M R$ par la relation
d'\'equivalence $(x, y+\phi(x))\sim (Tx,y)$ par~: pour tout $(x,y) \in
D_\phi$ et pour tout $t\in \M R$,
$$\tau_{\alpha,\phi}^t(x,y)=(x,y+t).$$
Si on identifie ce quotient au domaine fondamental $D_\phi=\{(x,y)\in
\M T\times \M R,\; 0\leq y<\phi(x)\}$, la mesure de Lebesgue
restreinte \`a $D_\phi$ est invariante.  On supposera ici $\phi$
int\'egrable et donc cette mesure finie.

Tout flot sp\'ecial $\tau_{\alpha,\phi}$ est topologiquement
conjugu\'e \`a un flot lin\'eaire de $\M T^2$ de direction
$(\alpha,1)$, not\'e $R_{\alpha,1}$, reparam\'etr\'e par une fonction
$\varphi>0$~: le champ de vitesse de ce flot est \'egal en chaque
point $(x,y) \in \M T^2$ \`a $(\alpha/\varphi(x,y), 1/\varphi(x,y))$
et la relation entre $\phi$ et $\varphi$ est donn\'ee pour tout $x \in
\M T$ par $$ \phi(x)=\int_0^1 \varphi(x+\alpha t,t)dt.$$
Quitte \`a faire une homoth\'etie de temps, on supposera
toujours $\phi$ normalis\'e de sorte que $\int_{\M T}\phi d\lambda
=1$, ce qui permet de comparer les valeurs propres de
$\tau_{\alpha,\phi}$ \`a celles du flot de translations non
renormalis\'e.

Comme les orbites du flot reparam\'etr\'e sont celles du flot
lin\'eaire, celui-ci est toujours minimal et uniquement ergodique. Les
propri\'et\'es spectrales des flots sp\'eciaux, ou des flots
lin\'eaires reparam\'etr\'es, d\'ependent de deux param\`etres, la
r\'egularit\'e de $\phi$ (ou, ce qui revient au m\^eme, de $\varphi$,
cf \cite{fakawi}), et les propri\'et\'es d'approximation
diophantiennes de $\alpha$. De fa\c con g\'en\'erale, l'absence de
m\'elange pour des fonctions $\varphi$ \`a variation born\'ee est
classique (voir \cite{koc}), et le probl\`eme du faible m\'elange est
donc naturel.

Le premier r\'esultat sur les valeurs propres est d\^u \`a  J. Von
Neumann (\cite{vn}). Il montre  le m\'elange faible de ces flots
sp\'eciaux lorsque $\phi$ est absolument continue par morceaux et
la somme de ses sauts est non nulle.
Au contraire, lorsque  $\alpha$ est diophantien et $\varphi \in
C^\infty$, le flot $R_{\alpha,1}$ reparam\'etr\'e par $\varphi$   est
conjugu\'e  au flot lin\'eaire $R_{\alpha,1}$ de $\M T^2$ et admet
donc dans ce cas un spectre purement discret (\cite{kol}).\\
Les \'etudes r\'ecentes de ces syst\`emes, dans le cas o\`u  $\alpha$
est Liouville,  montrent  leur diversit\'e spectrale~:  A. Katok et
J. Robinson ont donn\'e dans \cite{karo} des crit\`eres de m\'elange
faible, la g\'en\'ericit\'e du m\'elange faible a \'et\'e \'etablie
par B. Fayad dans  \cite{fa}, l'existence de flots \`a spectre mixte
avec  des fonctions plafond analytiques est un r\'esultat de B.
Fayad, A. Katok  et A. Windsor (\cite{fakawi}). En g\'en\'eral,  le
probl\`eme de savoir quelles sont les valeurs propres possibles,
quelle est la structure des facteurs Kronecker, ou si un tel  flot
peut-\^etre isomorphe \`a un flot de rotation autre que
$R_{\alpha,1}$  sont des questions anciennes  (\cite{vn},
\cite{kol}), qui restent d'actualit\'e (\cite{karo}, \cite{fakawi}).

Nous obtenons ici un certain nombre de r\'eponses \`a ces questions,
qui sont pr\'ecis\'ees dans le th\'eor\`eme~\ref{flot}. Avant d'en
donner l'\'enonc\'e, on pr\'ecise d'abord quelques notations et les
liens avec le th\'eor\`eme~\ref{marche}.

\medskip
Notons $e(\tau_{\alpha,\phi})$ le sous-groupe de $\M R$ des valeurs
propres du flot $\tau_{\alpha,\phi}$.  Rappelons qu'un r\'eel $t$ est
une valeur propre du flot s'il existe une fonction mesurable $F$ non
presque partout nulle v\'erifiant,
$F(\tau^{t'}_{\alpha,\phi}x)=\e{2i\pi t t'}F(x)$ presque partout pour
tout $t'$ r\'eel, et que pour le flot sp\'ecial cette condition
revient \`a une \'equation de cohomologie (\cite{vn}, voir
\cite{karo})~: on a $t \in e(\tau_{\alpha,\phi})$ si et seulement s'il
existe une fonction mesurable $f$ non presque partout nulle
satisfaisant
$$  \e{2i\pi t \phi} f=f\circ T . $$ 

On s'int\'eressera d'abord au flot sp\'ecial de fonction plafond
$\phi=1+\gamma\phi_\beta$, avec $\gamma >0$, $\beta \in \M R$ et
$\gamma\{\beta\}<1$, qu'on notera $\tau_{\alpha,\beta,\gamma}$. Les
valeurs propres de ce flot sont donc exactement les r\'eels $t$ pour
lesquels l'\'equation (\ref{eqbis}) avec le triplet $(\beta,s,t)$ ou 
$s=-t\gamma$ admet
une solution. La construction de tels triplets et des fonctions $f$
correspondantes fournira plusieurs exemples de flots sp\'eciaux
$\tau_{\alpha,\beta,\gamma}$ avec des propri\'et\'es particuli\`eres
de leurs valeurs propres et de leur facteur de Kronecker.

Nous \'etablirons ensuite un r\'esultat de r\'egularisation
(proposition~\ref{regul}, chapitre 8) pour montrer que ces
constructions conduisent \`a des exemples de flots sp\'eciaux avec des
fonctions plafonds lisses. Plus pr\'ecis\'ement nous montrons,
sous des conditions presque optimales d'approximation diophantienne
pour $\alpha$ et selon le d\'eveloppement d'Ostrowski de $\beta$, que 
$\phi_\beta$ est cohomologue \`a une fonction de classe $C^k$,
ou analytique, entrainant alors la conjugaison entre les flots sp\'eciaux. 

Avec ce dernier r\'esultat, le th\'eor\`eme~\ref{marche} permet
d\'ej\`a de retrouver le th\'eor\`eme 2 de \cite{fakawi} avec une
hypoth\`ese plus faible, et de r\'epondre au probl\`eme de l'existence
de spectres mixtes pour des nombres $\alpha$ non Liouville~:

\begin{coro}\label{mixte}
Pour tout $\alpha$ irrationnel il existe des flots sp\'eciaux  au
dessus de la rotation d'angle $\alpha$ admettant un groupe cyclique
de valeurs propres et donc un spectre mixte, avec de plus
\begin{itemize}
\item
Si $\inf_{q> 0}  q/(-\ln\|q\alpha\|)=0$, une
fonction plafond analytique sur $\M T$.

\item Si  $\inf_{q> 0} q^{k+1}\|q\alpha\|=0$, une
fonction plafond dans $C^k(\M T)$.
\end{itemize}
\end{coro}
 
\dem On consid\`ere le flot $\tau_{\alpha, \beta, \beta^{-1}}$, avec
$\beta>1$. Ses valeurs propres sont les r\'eels $t$ pour lesquels
l'\'equation $\e{2i\pi s\phi_\beta}=\e{2i\pi t} f/f\rond T$ avec
$s=-t\beta^{-1}$ admet des solutions.\\
Si $t \in \M Z\beta$ et donc $s\in \M Z$, l'\'equation est
trivialement satisfaite par les fonctions constantes (car
$s\phi_\beta=-s\beta\mod 1$)~: ceci montre que $\M Z\beta$ est contenu
dans le groupe des valeurs propres. Pour que le spectre soit mixte, il
suffit d'obtenir qu'il n'y ait pas d'autres valeurs propres (car le
flot sp\'ecial ne peut \^etre conjugu\'e \`a un flot de translations
sur $\M T$, qui admet une orbite de mesure pleine, cf
\cite{fakawi}).\\
D'apr\`es le th\'eor\`eme~\ref{marche}, il n'y a pas de solutions avec
$s$ non entier lorsque $\beta \notin H_1(\alpha)$, ou bien lorsque
$\beta \in H_1(\alpha)$ mais que $\{s, \sum_n\|b_ns\|^2 <\infty\}$ est
trivial. On r\'ealise cette derni\`ere condition, en choisissant
$b_n=0$ sauf pour une infinit\'e de $n$ pour lesquels $b_n=1$ et,
selon l'hypoth\`ese, $q_n/(-\ln\|q_n\alpha\|)$ ou
$q_n^{k+1}\|q_n\alpha\|$ suffisamment petit. Dans chacun des cas, cela
permet d'appliquer la proposition~\ref{regul}, qui donne le r\'esultat
souhait\'e.
\findem

\medskip
Le th\'eor\`eme suivant rassemble les autres r\'esultats d'existence
que nous obtenons pour les flots sp\'eciaux~:
\begin{theo}\label{flot} Soient les propri\'et\'es, pour un flot
sp\'ecial $\tau_{\alpha,\phi}$ au dessus de la translation par un
irrationnel $\alpha$~:
\begin{enumerate}
\item 
$\tau_{\alpha,\phi}$ admet une infinit\'e de valeurs propres
ind\'ependantes.
\item 
 $\tau_{\alpha,\phi}$  est conjugu\'e  \`a un flot de translations de
$\M T^2$, diff\'erent de $R_{\alpha,1}$. 
\item 
$\tau_{\alpha,\phi}$ est conjugu\'e  \`a un flot de translations sur
un sol\'eno\" \i de. 
\end{enumerate}
\medskip
\begin{itemize}
\item
  (i), (ii) et (iii) sont r\'ealisables avec $\phi \in C^\omega(\M
T)$ pour tout $\alpha$ v\'erifiant $\inf_{q>1} \left|q\ln
(q)/\ln(\|q\alpha\|)\right|=0$.
\item  
 (i), (ii) et (iii) sont r\'ealisables avec $\phi \in C^k(\M T)$,
pour tout $\alpha$ v\'erifiant  $\inf_{q\neq 0}
q^{k+2}\|q\alpha\|=0$.
\item
(i) est r\'ealisable avec $\tau_{\alpha,\beta,\gamma}$ pour un
ensemble non d\'enombrable dense de $\beta$, pour tout $\alpha$ tel
que $\inf_{q\neq 0}q\|q\alpha\|=0$ et tout $\gamma\in]0,1[$.
\item
(ii) et (iii) sont r\'ealisables avec $\tau_{\alpha,\beta,\beta^{-1}}$
pour un ensemble non d\'e\-nom\-bra\-ble dense de $\beta>1$ pour tout
$\alpha$ tel que $\inf_{q\neq 0} q^2\|q\alpha\|=0$.
\end{itemize}
\end{theo}
Des constructions r\'eciproques permettent de donner des conditions
sur les flots de translations entra\^{\i}nant qu'ils soient
conjugu\'es \`a des flots sp\'eciaux~:
\begin{theo}\label{reci}
On a les propri\'et\'es suivantes~:
 \begin{enumerate}
\item 
Tout flot de translations $R_{s,1}$ o\`u $s$ v\'erifie $\inf_{q>1}
q^{k+2}\|q s\|=0$ est mesurablement conjugu\'e, \`a une homoth\'ethie
de temps pr\`es, \`a un flot sp\'ecial $\tau_{\alpha,\phi}$
 avec $\alpha\neq s$ et  $\phi \in C^k$. Si  $\inf_{q\neq 0}
\left|q\ln (q)/\ln \|qs\|\right|=0$, alors $R_{s,1}$ est conjugu\'e
\`a un flot sp\'ecial $\tau_{\alpha,\phi}$ avec $s \neq\alpha$ et
$\phi \in C^\omega$.
\item 
Tout flot de translations sur un sol\'eno\"\i de est mesurablement
conjugu\'e,  \`a une homoth\'ethie de temps pr\`es, \`a un flot
sp\'ecial $\tau_{\alpha,\phi}$ de fonction plafond analytique.
\end{enumerate}
\end{theo}

\noindent{\bf Remarque~:}
Toutes les constructions des th\'eor\`emes~\ref{flot} et \ref{reci} fournissent \'egalement des exemples pour les flots $R_{\alpha,1}$ reparam\'etr\'es par une fonction $\varphi$ de m\^eme r\'egularit\'e que $\phi$ (voir \cite{fakawi}).

\subsection{Applications aux \'echanges de 3 intervalles.}

Un \'echange d'intervalles est une transformation bijective de
l'intervalle $[0,1]$ affine par morceaux, et qui pr\'eserve les
longueurs. Ces transformations classiques ne sont jamais
m\'elangeantes (\cite{ka80}). L'unique ergodicit\'e et  le m\'elange
faible sont des r\'esultats presque s\^urs   difficiles et connus
depuis longtemps pour le premier gr\^ace aux travaux de H. Masur et
W. Veech (\cite{ve1}, \cite{mas}) mais d\'emontr\'e tr\`es
r\'ecemment pour le second par A. Avila et G. Forni  (\cite{avfo}).

Le cas le plus simple apr\`es les \'echanges de deux intervalles, qui
sont des translations sur le tore, est  celui des \'echanges de trois
intervalles. Pour ceux-ci,  le m\'elange faible presque s\^ur  est un
r\'esultat classique (\cite{katste}) qui provient  d'une remarque
simple~:  tout \'echange de trois intervalles est l'induit d'une
translation irrationnelle sur un intervalle. Dans ce cas, la stricte
ergodicit\'e est \'egalement connue depuis longtemps (lorsque les
discontinuit\'es sont rationnellement ind\'ependantes, voir
\cite{ke75}). En revanche la caract\'erisation des valeurs propres et
la description du facteur Kronecker en fonction des param\`etres
reste un probl\`eme ouvert.  L'existence de valeurs propres
rationnelles dans certains cas particuliers est une cons\'equence de
\cite{vee} et \cite{ste}. Des  travaux  r\'ecents de S. Ferenczi, C.
Holton et L. Zamboni  (\cite{fehozaIII}) fournissent des r\'eponses
partielles \`a ce probl\`eme, en donnant en particulier des exemples
d'\'echanges de 3 intervalles non triviaux conjugu\'es \`a une
rotation irrationnelle du cercle. Ils montrent \'egalement que tout
r\'eel quadratique est une valeur propre d'un \'echange de trois
intervalles non trivial.

Inversement n'importe quel induit d'un rotation irrationnelle sur un
intervalle est aussi un \'echange de trois intervalles. C'est sous cette
forme que nous construirons les \'echanges de trois
intervalles. Pour $\beta \in ]0,1[$, nous noterons  $T_{\alpha,\beta}$
la transformation induite par $T$ sur $[0,\beta[$, d\'efinie par 
$$T_{\alpha,\beta}(x)=T^{n(x)}(x) \quad \hbox{o\`u } n(x)=\min(k>0,
T^kx \in [0,\beta[\}$$
pour tout $x \in [0,\beta[$, 
la probabilit\'e invariante \'etant la mesure de Lebesgue renormalis\'e
qu'on note $\lambda_\beta$.\\
Les valeurs propres de $T_{\alpha,\beta}$ sont encore d\'etermin\'ees
par une \'equation de cohomologie (cf. \cite{katste}, \cite{karo})~:
notons $e(T_{\alpha,\beta})$ le sous-groupe de $\M T$ des valeurs
propres de $T_{\alpha,\beta}$. Il est constitu\'e des $s $ de $\M T$
pour lesquels il existe une fonction mesurable $f$ non presque partout
nulle v\'erifiant
$$\e{2i\pi s \unde{[0,\beta[}}f=f\rond T.$$
Donc le relev\'e dans $\M R$ de $e(T_{\alpha,\beta})$ est
l'ensemble des r\'eels $s$ pour lesquels l'\'equa\-tion (\ref{eqbis})
avec $t=-s\beta \mod 1$ admet une solution.

Nous obtenons \`a partir du th\'eor\`eme~\ref{marche} des r\'esultats
pour les \'echanges de trois intervalles analogues \`a ceux pour les
flots sp\'eciaux~:

\begin{theo}\label{3i}
On a les propri\'et\'es suivantes~:
\begin{enumerate}
\item
Pour tout $\alpha$  tel que $\inf_{q\neq 0} q\|q\alpha\|=0$, il
existe un ensemble non d\'e\-nom\-bra\-ble dense de $\beta$ pour lesquels
$T_{\alpha,\beta}$ admet une infinit\'e de valeurs propres
rationnellement ind\'ependantes.
\item
Pour tout $\alpha$  tel que $\inf_{q\neq 0}q^2\|q\alpha\|=0$, il
existe un ensemble non d\'e\-nom\-bra\-ble dense de $\beta$ tels que
$T_{\alpha,\beta}$ soit conjugu\'e \`a une rotation irrationnelle du
cercle. 
\item 
Pour tout $\alpha$ tel que $\inf_{q\neq 0}q^2\|q\alpha\|=0$, il
existe un ensemble non d\'e\-nom\-bra\-ble dense de $\beta$ tels que
$T_{\alpha,\beta}$ soit conjugu\'e \`a un odom\`etre. 
\item 
Toute rotation irrationnelle d'angle $s$ du cercle tel que
$\inf_{q\neq 0}q^2\|qs\|=0$ est conjugu\'ee \`a un \'echange de 3
intervalles.
\item
Tout odom\`etre est conjugu\'e \`a un \'echange de 3 intervalles.
\end{enumerate}
\end{theo}

\noindent{\bf Remarque~:} Notons que  (ii) est une g\'en\'eralisation
du th\'eor\`eme 1 de \cite{fehozaIII}.

\subsection{Plan de l'article.}

Ce travail est constitu\'e de 3 parties et d'un appendice.  La
premi\`ere partie regroupe l'introduction et un chapitre
pr\'eliminaire (chapitre 2) qui pr\'esente la d\'emarche et les outils
de la d\'emonstration du th\'eor\`eme~\ref{marche}.  Nous y discutons
aussi les repr\'esentations en tours d'intervalles de la translation
irrationnelle li\'ees \`a la fraction continue utiles pour les
chapitres 3, 4 et 7.

La seconde partie regroupe les chapitres 3 et 4.  Le chapitre 3
compl\`ete la preuve du th\'eor\`eme~\ref{marche} en en traitant les
points les plus techniques.  Le chapitre 4, qui s'appuie sur le
pr\'ec\'edent, est consacr\'e au cas des fonctions en escalier
g\'en\'erales et contient la preuve du th\'eor\`eme \ref{cn}.

Les chapitres 5, 6, 7 et 8 forment la troisi\`eme partie~: ils
contiennent les applications aux valeurs propres des flots sp\'eciaux
et des \'echanges de trois intervalles. Le chapitre 5 est une approche
g\'en\'erale du probl\`eme qui met en relation ces deux types de
transformations gr\^ace aux tours de Kakutani.
Le chapitre 6 est consacr\'e aux constructions avec une infinit\'e de
valeurs propres ind\'ependantes (points (i) des th\'eor\`emes
\ref{flot} et \ref{3i}) et le chapitre 7 aux isomorphismes avec des
translations~: on y propose d'abord une construction  de tours de Kakutani
o\`u l'on peut identifier le facteur Kronecker.  Nous
montrons ensuite les points (ii) et (iii) des th\'eor\`emes
\ref{flot} et \ref{3i}, ainsi que leurs r\'eciproques (th\'eor\`eme
\ref{reci}, points (iv) et (v) du th\'eor\`eme~\ref{3i}).
Le chapitre 8 est constitu\'e de la preuve de la proposition
\ref{regul},  qui permet d'achever les preuves des th\'eor\`emes~\ref{flot} et \ref{reci} par r\'egularisation des fonctions en 
escalier.\\
Enfin on trouvera dans l'appendice la preuve de la proposition
\ref{hsat}.
 
\medskip
Les parties 2 et 3 peuvent \^etre lues ind\'ependamment l'une de
l'autre~: les objets indispensables au chapitre 7 (tours de la
translation, approximations des fonctions de transfert) sont introduits
dans le chapitre 2. Dans la partie 3, les chapitres 6 et 7 sont
ind\'ependants.
  
\section{Pr\'eliminaires.}

% Nous pr\'esentons dans ce chapitre les outils n\'ecessaires \`a la
% preuve du th\'eor\`eme~\ref{marche}. Comme certains objets d\'efinis
% dans la preuve sont 
% indispensables pour le chapitre 7, et par souci de clart\'e pour le
% lecteur nous commen\c cons cette partie par une description
% heuristique de la d\'emonstration du th\'eor\`eme~\ref{marche}.
% 
\subsection{Des notations et d\'efinitions}

Dans tout l'article, l'irrationnel $\alpha$ est suppos\'e donn\'e sauf
pour les constructions r\'eciproques aux chapitres 6 et 7, et $T$
d\'esigne la translation par $\alpha$ sur $\M T$.  Pour la preuve du
th\'eor\`eme~\ref{marche}, dans ce chapitre et le suivant, le r\'eel
non entier $s$ est aussi suppos\'e fix\'e.

\subsubsection{Cocycles. Les quasi-cobords triviaux}\label{cobord}

Deux fonctions r\'eelles mesurables $\phi$ et $\psi$ sur $\M T$ sont
cohomologues (additivement) s'il existe une fonction r\'eelle
mesurable $f$ telle que $\phi-\psi=f-f\rond T$. On dit que $\phi$ est
un cobord si elle est cohomologue \`a $0$ et qu'elle est un
quasi-cobord si elle est cohomologue \`a une constante. La fonction
$f$ correspondante, unique \`a une constante pr\`es, est appel\'ee
fonction de transfert.\\
On note, pour tout entier $n> 0$,
$$\phi^{(n)}= \sum_{0 \leq j <n}\phi\rond T^j, \quad 
\phi^{(-n)}= -\phi^{(n)} \rond T^{-n}, \quad 
\hbox{et }\phi^{(0)}=0.$$
On utilisera souvent les relations $\phi-\phi\rond
T^n=\phi^{(n)}-\phi^{(n)}\rond T$ quel que soit $\phi$, et
$\phi^{(n)}=f-f\rond T^n$ lorsque $\phi=f-f\rond T$.
\\
Pour une fonction $\varphi= \exp{2i\pi \phi}$ de module 1 sur $\M T$,
on adopte les m\^emes d\'efinitions et notations, en rempla\c{c}ant
l'addition par la multiplication, ou on parlera de cohomologie modulo 
1 pour $\phi$. On a alors $\varphi^{(n)}= \exp{2i\pi \phi^{(n)}}$ pour
tout entier $n$.

On pose $\varphi_\beta=\exp{2i\pi s \phi_{\beta}}$ pour tout $\beta \in
\M T$ (ou $\M R$).  On rappelle que $\phi_\beta$ est d\'efinie par
$\phi_\beta=\unde{[0,\{\beta\}[}-\{\beta\}$, avec aussi
$\phi_{\beta}(x)=\{x-\beta\}-\{x\}$.  En particulier $\phi_\alpha$ est
un cobord.\\
Pour tout entier $n$, on a la propri\'et\'e ``de cocycle"
\begin{equation}\label{coc}
\phi_{n\alpha+\beta}=\phi_{n\alpha}+\phi_\beta\circ T^{-n}
\end{equation}
et, en notant $\omega$ la fonction
$x\mapsto \{x\}$ de $\M T$ dans $\M R$, on a aussi 
$$\phi_{n\alpha}=\omega\rond T^{-n}-\omega=\omega_n-\omega_n \circ T$$
avec
$$\omega_n=-\omega^ {(-n)}.$$
Donc $\phi_{n\alpha}$ est un cobord avec une fonction de transfert
$\omega_n$ affine par morceaux.  Si $n\ge 0$, alors
$\omega_n=\sum_{1\le j<n}\omega\rond T^{-j}$ admet $n$
discontinuit\'es aux points $T^j0$ pour $1 \leq j \leq n$, les sauts
en chacun de ces points valant $-1$~; si $n<0$ les points de
discontinuit\'e sont les $T^{j}0$ pour $n<j\leq 0$ avec des sauts
\'egaux \`a $+1$, et dans les deux cas $\omega_n$ est de pente $n$
dans chaque intervalle o\`u elle est continue.

D'apr\`es cela, $\varphi_{n\alpha}$ est un cobord multiplicatif pour
tout entier $n$, avec la fonction de transfert $\e{2i\pi s\omega_n}$.
En multipliant cette fonction par les exponentielles $\e{2i\pi mx}$,
on voit qu'elle est aussi cohomologue \`a toutes les valeurs propres
$\e{2i\pi m\alpha}$ de $T$.  En notation additive, quels que soient
les entiers $m$ et $n$,
$$
s\phi_{n\alpha}=  m\alpha+
\omega_{n,m}-\omega_{n,m}\circ T\mod 1,
$$
avec
$$\omega_{n,m}= s\omega_n + m\omega.$$
Pour la suite, on observe que $\omega_{n,m}$ est encore une
fonction affine par morceaux de m\^eme points de discontinuit\'e
que $\omega_n$, avec des sauts de $\mp s$ selon le signe de $n$ et 
la pente $ns+m$, entre deux sauts.

\subsubsection{Fraction continue et d\'eveloppements
d'Ostrowski}\label{ostrow}

Pour le d\'eveloppement en fraction continue, on identifie $\alpha$
\`a un r\'eel de $]0,1[$.  On notera $\alpha_0=\alpha$ et, pour $n\geq
1$,
$$\alpha_n=\|q_n\alpha\|.$$
On rappelle les propri\'et\'es classiques utiles pour la suite (voir
\cite{kin} par exemple).  La suite $(\alpha_n)$ est strictement
d\'ecroissante.  Pour les premiers termes, $q_0=1$, $p_0=0$,
$q_1=a_1$, $p_1=1$ et, pour $n\ge 1$,
\begin{eqnarray*}
    q_{n+1} =a_{n+1}q_n+q_{n-1},&\quad& p_{n+1}=a_{n+1}p_n+p_{n-1},\\
    p_nq_{n-1}-p_{n-1}q_n=(-1)^n,&\quad&
 \<q_n\alpha \> = q_n\alpha-p_n=(-1)^n\alpha_n,\\
 \alpha_{n-1}=a_{n+1}\alpha_n+\alpha_{n+1},&\quad&
 q_n\alpha_{n-1}+q_{n-1}\alpha_n = 1
\end{eqnarray*}
et
\begin{equation}\label{qalpha}
\|q\alpha\|\geq   \alpha_{n-1}  \quad \hbox{lorsque }
0<|q|<q_{n}.
\end{equation}
Dans le cas o\`u la suite des quotients partiels est non born\'ee, on
utilisera fr\'e\-quem\-ment les \'equivalences, le long d'une suite
d'indices $n$ tels que $a_{n+1}\rightarrow +\infty$~:
$$ q_n\alpha_n \sim \inv{a_{n+1}}\qquad \hbox{ et }\qquad
\alpha_{n-1}\sim \inv{q_n},$$
avec dans tous les cas $\inv{2}<q_n\alpha_{n-1}<1$ et
$\inv{2}\alpha_{n-1}<a_{n+1}\alpha_n<\alpha_{n-1}$.

\medskip
En particulier pour un d\'eveloppement $\beta=\sum_{n\geq
0}b_nq_n\alpha\mod 1$, la condition $\sum_{n\geq
0}|b_n|/a_{n+1}<\infty$ (soit $\beta\in H_1(\alpha)$) \'equivaut \`a
$\sum_{n\geq0}|b_n|q_n\alpha_n <\infty$.\\
Pour un r\'eel $\beta$ non modulo 1, un tel d\'e\-ve\-lop\-pe\-ment
d'Ostrowski sym\'etrique s'\'ecrit $$\beta=k_1\alpha-l_1+\sum_{n\geq
1}b_n\<q_n\alpha\>,$$
avec $k_1$, $l_1$ et les $b_n$ entiers, et on note pour tout $n\geq 1$
$$\beta_n=k_n\alpha-l_n=
k_1\alpha-l_1+\sum_{1\leq j<n}b_j\<q_j\alpha\>,$$
o\`u $k_n=k_1+\sum_{1\leq j<n}b_jq_j$ et de m\^eme pour $l_n$.  On
utilisera ces d\'eveloppements uniquement dans le cas o\`u $\beta\in
H_{\infty}(\alpha)$, c'est-\`a-dire $b_n/a_{n+1}\to 0$.  On a alors
clairement $k_n=o(q_n)$ et $\beta-\beta_n=o(\alpha_{n-1})$, et plus
pr\'ecis\'ement~:
\begin{lemme}\label{unicite}
Avec les notations ci-dessus, soient $n\geq 1$ et $\varepsilon>0$. Si
$|b_j|/a_{j+1}\leq\varepsilon$ pour tout $j\geq n$, alors
$$|\beta-\beta_n|< 2\varepsilon\alpha_{n-1}.$$
En supposant de plus $|k_n|<q_n$ et $\varepsilon<\inv{2}$, si
$\beta\in\M Z$ alors $k_n$ et tous les $b_j$ pour $j\geq n$ sont nuls.
\end{lemme}
\dem L'in\'egalit\'e r\'esulte de
$$\sum_{j\geq n}|b_j|\alpha_j\leq\varepsilon\sum_{j\geq n}a_{j+1}\alpha_j
=\varepsilon\sum_{j\geq n}(\alpha_{j-1}-\alpha_{j+1})
=\varepsilon (\alpha_{n-1}+\alpha_n).$$
Maintenant, si $|k_n|<q_n$, $\varepsilon<\inv{2}$ et $\beta\in\M Z$,
on a $\|k_n\alpha\|=\|\beta_n-\beta\|<\alpha_{n-1}$ donc
n\'ecessairement $k_n=0$ d'apr\`es (\ref{qalpha}).  On a alors aussi
$|k_{n+1}|=|b_n|q_n<q_{n+1}$ donc $k_{n+1}=b_n=0$ et, par une
r\'ecurrence imm\'ediate, $b_j=0$ pour tout $j\geq n$.
\findem
Ce lemme montre l'unicit\'e du d\'eveloppement
$\beta=k_n\alpha-l_n+\sum_{j\geq n}b_n\<q_n\alpha\>$ d\`es que
$|k_n|<q_n/2$ et $|b_j|\leq a_{j+1}/4$ pour $j\geq n$. En particulier 
$\beta\in \M Z\alpha+\M Z$ si et seulement les $b_j$ sont nuls sauf
pour un nombre fini.

Ainsi tout r\'eel $\beta$ de $H_{\infty}(\alpha)$ admet de bonnes
approximations $\beta_n=k_n\alpha\mod 1$ avec $k_n=o(q_n)$.
R\'eciproquement, on obtient $\beta\in H_{\infty}(\alpha)$ sous une
condition plus faible~:
\begin{lemme}\label{betainHinfini}
Soit $\beta$ un r\'eel.  Si pour tout $n\geq 1$ il existe un entier
$k_n=o(q_{n+1})$ tel que $\|\beta-k_n\alpha\|=o(\alpha_{n-1})$, alors
$\beta\in H_{\infty}(\alpha)$.
\end{lemme}
\dem On utilise la proposition~\ref{hsat}.  Sous les conditions
donn\'ees $\|q_n\beta\|\leq
\|q_nk_n\alpha\|+q_n\|\beta-k_n\alpha\|\leq
|k_n|\|q_n\alpha\|+o(q_n\alpha_{n-1})=o(q_{n+1}\alpha_n)
+o(q_n\alpha_{n-1})$, et cette quantit\'e tend vers $0$.
\findem

\subsubsection{Tours de la translation}
Soit $B$ un bor\'elien de $\M T$~; on appelle tour de base $B$ une
suite $(T^jB)_{0 \leq j <h}$ d'images successives de $B$ deux-\`a-deux
disjointes.  La hauteur de la tour est $h$ et les $T^jB$ sont les
\'etages de la tour.  Par abus de langage on appelera aussi tour la
r\'eunion des \'etages.

Il est bien connu que le d\'eveloppement de $\alpha$ en fraction
continue d\'efinit par induction, \`a chaque ordre $n\geq 1$, une
partition de $\M T$ en deux tours de la translation, form\'ees
d'intervalles.  En fait, si $B_n$ est un intervalle quelconque de
longueur $\alpha_{n-1}$, les $T^jB_n$ pour $0\leq j<q_n$ sont
deux-\`a-deux disjoints d'apr\`es (\ref{qalpha}) et forment donc une
premi\`ere tour de hauteur $q_n$ et de mesure
$q_n\alpha_{n-1}>\inv{2}$, qu'on appellera tour majeure d'ordre $n$.
Pour la m\^eme raison, pour tout $x\in\M T$, la suite $(T^jx)_{0\leq
j<q_n}$ ne contient aucun point entre $T^jx$ et
$T^{j+q_n-q_{n-1}}x=T^jx+(-1)^n(\alpha_n+\alpha_{n-1})$ lorsque $0\leq
j<q_{n-1}$~; si par exemple $n$ est pair et
$B_n=[x_n,x_n+\alpha_{n-1}[$, on voit ainsi que le compl\'ementaire
de la tour majeure est la r\'eunion des intervalles $T^j[x_n+\alpha_{n-1},
x_n+\alpha_{n-1}+\alpha_n[$ pour $0\leq j<q_{n-1}$, qui forment une
tour mineure de hauteur $q_{n-1}$ et de mesure
$q_{n-1}\alpha_n<\inv{2}$.

\medskip Enfin, nous introduisons deux notations qui serviront surtout
\`a majorer les variations de fonctions sur les \'etages des tours~:

-- Si $f$ est une fonction mesurable sur $\M T$
et $B$ un bor\'elien de $\M T$ de mesure non nulle, on appelle
variation moyenne de $f$
sur $B$ la quantit\'e
$$V_B(f)=\inv{\lambda(B)^2} \iint_{B^2}|f(x)-f(y)|dxdy\;.$$
\\
-- Nous utilisons les notations habituelles $o$ et ${\cal O}$.
On notera
$$g(x)=\Oun(f(x))$$
lorsque $|g(x)|\leq |f(x)|$ pour tout $x$ et  de m\^eme pour des
suites.

\subsection{Outils et sch\'ema de la preuve du th\'eor\`eme
\ref{marche}.}\label{idee}

Nous pr\'esentons ici un plan d\'etaill\'e de la d\'emonstration du
th\'eor\`eme~ \ref{marche}, avec les principaux arguments.
L'irrationnel $\alpha$ et le r\'eel non entier $s$ \'etant fix\'es,
nous discutons les conditions sur un couple $(\beta,t)$ de r\'eels
modulo 1 pour que l'\'equation~(\ref{eqbis}), $\varphi_\beta=\e{2i\pi
t}f/f\rond T$, admette une solution.

\subsubsection{Principe -- approximations de $\beta$ et $t$.}

La d\'emonstration sera bas\'ee sur l'\'etude de la convergence des
solutions approch\'ees fournies par les fonctions de transfert
``triviales" obtenues pour par les approximations de $\beta$ et $t$
dans $\M Z\alpha$, selon le paragraphe~\ref{cobord}.

Pour montrer que les conditions du th\'eor\`eme suffisent, c'est la
m\'ethode de K.~Merrill, qui montre que $\varphi_\beta$ est un
quasi-cobord sous une hypoth\`ese tr\`es voisine mais plus forte~:
dans \cite{mer}, on trouve la condition $\sum\|b_ns\|<\infty$
au lieu de $\sum\|b_ns\|^2<\infty$. Nous aurons surtout \`a pr\'eciser
exactement les conditions de la convergence.

Dans l'autre sens, on montrera que les solutions approch\'ees doivent
converger vers $f$, \`a un choix de constantes pr\`es, et on est
ramen\'e au probl\`eme de caract\'eriser cette convergence. Mais il
sera plus simple de savoir d'abord que $\beta$ et $t$ sont bien
approch\'es par $\M Z\alpha$, c'est-\`a-dire $\beta$, $t\in
H_\infty(\alpha)$.

Le fait que $\beta$ doit appartenir \`a $H_\infty(\alpha)$ est connu 
depuis W.~ÊVeech \cite{vee} et M.~Stewart \cite{ste} qui
montrent aussi que ses
coefficients $b_n$ doivent satisfaire $\|b_ns\|\to 0$ -- et ce sont
les seules conditions n\'ecessaires d\'ej\`a connues.  Ils utilisent
une m\'ethode diff\'erente, bas\'ee sur les temps de rigidit\'e
de $T$~: lorsque $\varphi=\e{2i\pi t}f/f\rond T$, on a
$\varphi^{(n)}=\e{2i\pi nt}f/f\rond T^n$ pour tout $n$ et $\e{-2i\pi
nt}\varphi^{(n)}$ doit tendre vers 1 en mesure lorsque $n$ tend vers
l'infini avec $\|n\alpha\|\to 0$, donc on obtient des conditions sur
la r\'epartition des $\varphi^{(n)}$ (c'est aussi l'utilisation de la
r\'eciproque qui conduit aux crit\`eres suffisants \'etablis dans
\cite{vee} et \cite{ste}).

Nous avons besoin de cette propri\'et\'e uniquement lorsque $n$ est un
d\'e\-no\-mi\-na\-teur de la fraction continue de $\alpha$, \`a
savoir~: si $\varphi_{\beta}$ est cohomologue \`a $\e{2i\pi t}$, on doit avoir
\begin{equation}\label{rigide}
 \e{-2i\pi q_nt}\int_\M T \varphi_\beta^{(q_n)}d\lambda \to1,
\end{equation}
ou de fa\c{c}on \'equivalente $s\phi_\beta^{(q_n)}-q_nt\to 0\mod 1$ en
mesure.\\
On a $\phi_\beta^{(q_n)}(x)=\sum_{0\leq
j<q_n}\unde{[0,\{\beta\}[)}(x)-q_n\beta$~: c'est l'\'ecart \`a la
moyenne du nombre de passages de $(x+j\alpha)_{0\leq j<q_n}$ dans
$[0,\{\beta\}[$ et il est connu que, pour un d\'enominateur $q_n$, on
obtient au plus 3 valeurs. Ces valeurs sont de la forme $j-q_n\beta$
avec des entiers $j$ cons\'ecutifs car les sauts de
$\phi_\beta^{(q_n)}$ valent $\pm1$. Lorsque $s\in \inv{2}\M Z$, la
condition (\ref{rigide}) implique qu'il existe $j_n$ tel que
$\phi_\beta^{(q_n)}-(j_n-q_n\beta)$ tende vers 0 en mesure~; comme
$\int_{\M T}\phi_\beta^{(q_n)}d\lambda=0$, on doit avoir
$j_n-q_n\beta\to 0$ d'o\`u $\|q_n\beta\|\to 0$ et de plus
$s\phi_\beta^{(q_n)}\to 0$ en mesure, donc aussi $\|q_nt\|\to 0$.

Nous discuterons plus pr\'ecis\'ement la r\'epartition de
$\phi_\beta^{(q_n)}$ et les cas restants au \S~\ref{betaqn}.  Pour la
suite de ce paragraphe, on suppose que $\beta$ et $t$ appartiennent
\`a $H_\infty(\alpha)$.  On note les
d\'eveloppements d'Ostrowski sym\'etriques de $\beta$ et $t$ dans $\M T$ et
leurs approximations dans $\M Z \alpha$  comme
au \S~\ref{ostrow}~:
\begin{equation}\label{ostro}
\beta=\sum_{j\geq 0}b_nq_n\alpha,\qquad
t =\sum_0^\infty b'_n q_n\alpha,
\end{equation}
et
$$\beta_n=k_n\alpha=\sum_{0\leq j<n}b_jq_j\alpha,\qquad
t_n= k'_n\alpha=\sum_0^{n-1}b'_jq_j\alpha$$
pour $n\geq 1$, les $b_n$ et les $b'_n$ \'etant des entiers v\'erifiant
$b_n/a_{n+1}\to0$, $b'_n/a_{n+1}\to0$. On rappelle les estimations  
$$
\left\lbrace{ 
\begin{array}{lcl}
k_n=o(q_n), & &\beta-\beta_n=o(\alpha_{n-1}),\\
k'_n=o(q_n), &\hbox{et} &t-t_n=o(\alpha_{n-1}).
\end{array}
}\right.$$

\subsubsection{Description de la suite des fonctions de
transfert.}\label{transfert}
Pour tout $n\geq 1$, $s\phi_{\beta_n}=\phi_{k_n\alpha}$ est
cohomologue \`a $t_n=k'_n\alpha$ modulo 1. On
d\'ecrit ici la suite des fonctions de transfert, d\'efinies \`a
une constante pr\`es, en vue d'\'etudier leur convergence. On pose
$$\tilde f_n=\omega_{k_n,k'_n}+c_n, \quad
f_n=\exp{2i\pi \tilde f_n}$$
o\`u $\omega_{k_n,k'_n}$ est la fonction affine par morceaux d\'efinie
au \S~\ref{cobord} et $c_n$ est une constante qu'on choisira ensuite
selon les besoins~; alors
\begin{equation}\label{cob}
s\phi_{\beta_n}= t_n+\tilde f_n-\tilde f_n\circ T \mod 1.
\end{equation}
On pose aussi 
$$\tilde \theta_n=\tilde f_{n+1}-\tilde f_n \quad \hbox{et }
\theta_n=f_{n+1}/f_n=\exp{2i\pi \tilde \theta_n}.$$
Comme $t_{n+1}-t_n=b'_nq_n\alpha$ et
$\phi_{\beta_{n+1}}-\phi_{\beta_n}=\phi_{b_nq_n\alpha}\rond T^{-k_n}$
d'apr\`es (\ref{coc}),
\begin{equation}\label{thetatransf}
s\phi_{b_nq_n\alpha}\circ T^{-k_n}= b'_nq_n\alpha+\tilde \theta_n-\tilde
\theta_n\circ T \mod 1
\end{equation}
et $ \tilde \theta_n$ est \'egal \`a $\omega_{b_nq_n,b'_nq_n}\circ
T^{-k_n}$ \`a une constante pr\`es.

On suppose $n$ assez grand pour que $b_n$ et $b'_n$ soient petits
devant $a_{n+1}$~; alors $\|b_nq_n\alpha\|=|b_n|\alpha_n$ est petit
devant $\alpha_{n-1}$. Pour d\'ecrire le comportement de
$\tilde f_n$ et $\tilde \theta_n$, on doit distinguer plusieurs cas,
selon les signes de $k_n$ et $b_n$, et selon la parit\'e de $n$. On
suppose d'abord ici $k_n\geq 0$, $b_n\geq 0$ et $n$ pair.

\begin{figure}
\begin{center}\caption{Tour majeure d'ordre $n$, partie
principale.}\label{2tours} 
\input{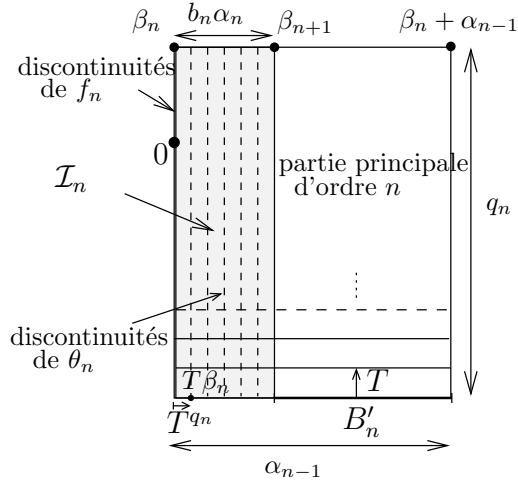}
\end{center}\end{figure}

Alors les discontinuit\'es de $\tilde f_n$ et $\tilde \theta_n$ sont
situ\'ees aux points $(T^j0)_{1\leq j\leq k_n}$ et $(T^j
\beta_n)_{1\leq j\leq b_nq_n}$ respectivement~; ces points se
repr\'esentent bien dans la tour majeure
$${\C T}_n=\bigl(T^{-j}[\beta_n,\beta_n+\alpha_{n-1}[\bigr)_{0\leq j
<q_n},$$
de base $B_n=T^{-q_n+1}[\beta_n,\beta_n+\alpha_{n-1}[$, qu'on
appellera simplement dans la suite la tour majeure d'ordre $n$.
Plus pr\'ecis\'ement, les discontinuit\'es de $\tilde f_n$ se trouvent
\`a une extr\'emit\'e des derniers \'etages et celles de $\tilde
\theta_n$, qu'on peut aussi \'ecrire $T^{-j}(\beta_n+k\alpha_n)$ avec
$1\leq k<b_n$ et $0\leq j<q_n$, se trouvent dans la sous-tour
\begin{equation}\label{in}
\C I_n=\bigl(T^{-j}[\beta_n, \beta_{n+1}]\bigr)_{0\leq j<q_n}
\end{equation}
de base $I_n=T^{-q_n+1}[\beta_n, \beta_{n+1}]$ et de mesure
$b_nq_n\alpha_n$ convergeant vers $0$.\\
La sous-tour compl\'ementaire $\C T'_n=(T^jB'_n)_{0\leq j<q_n}$ de
base $B'_n=B_n\setminus I_n$ sera appel\'ee la {\em tour principale}
d'ordre $n$. On a repr\'esent\'e figure \ref{2tours} la tour majeure
avec les deux sous-tours et les points de discontinuit\'e de $\tilde
f_n$ et $\tilde \theta_n$.

On obtient les propri\'et\'es suivantes~:
\begin{enumerate}
\item Sur chaque \'etage $T^jB_n$ de $\C T_n$ ($0\leq j<q_n$), la
fonction $\tilde f_n$ est affine de pente $k_ns+k'_n$, donc sa
variation y est petite~:
\begin{equation}\label{varfn}
V_{T^j B_n}(\tilde f_n)\leq |k_ns+k'_n|\alpha_{n-1}<|k_ns+k'_n|/q_n\to
0.
\end{equation}
\item  Sur chaque \'etage $T^jB'_n$ de $\C T'_n$ ($0\leq j<q_n$), la
fonction $\tilde \theta_n$ est affine de pente $(b_ns+b'_n)q_n$, et sa 
variation est de l'ordre de $|b_ns+b'_n|$ (car
$\inv{2}<q_n\alpha_{n-1}<1$).
\item  Hors de $[\beta_n,\beta_{n+1}[$, on a
$\phi_{b_nq_n\alpha}\circ T^{-k_n}=-\{b_nq_n\alpha\}$ donc 
$\tilde \theta_n\circ T-\tilde \theta_n= (b_ns+b'_n)q_n\alpha\mod 1$
d'apr\`es (\ref{thetatransf}). En it\'erant, pour tout point $x\in B_n$, 
\begin{equation}\label{invariant} 
 \|\tilde \theta_n\circ T^j(x)-\tilde \theta_n(x)\| \leq
|b_ns+b'_n|q_n\alpha_n\quad\mbox{pour }0\leq j<q_n.
\end{equation}
\end{enumerate}
Comme $q_n\alpha_n$ est de l'ordre de $1/a_{n+1}$, cette derni\`ere
in\'egalit\'e signifie que $\tilde \theta_n$ est presque invariante
par $T$ sur $\C T_n$.

Pour les autres cas ($k_n<0$, $b_n<0$ ou $n$ impair), on pr\'ecisera
au \S~\ref{df} les modifications \`a apporter pour obtenir les m\^emes
propri\'et\'es.

\subsubsection{Pourquoi les conditions du th\'eor\`eme~\ref{marche}
suffisent-elles?}\label{csidee}

On suppose que $\beta$ et $t$ satisfont  les conditions du
th\'eor\`eme~\ref{marche}, c'est-\`a-dire $\sum_n |b_n|q_n\alpha_n
<\infty$, les coefficients du d\'eveloppement de $t$
v\'erifient $b'_n=-[b_ns]$ pour $n$ assez grand et
$$\sum_{n\geq 0}(b_ns+b'_n)^2 <\infty.$$
Pour montrer que l'\'equation (\ref{eqbis}) admet une solution, il
suffit de montrer que la suite des fonctions $f_n= \exp{2i\pi\tilde
f_n}$ admet une valeur d'adh\'erence faible $f$ non nulle dans $L^2(\M
T)$ (par exemple). En effet, comme $\phi_{\beta_{n}}$ converge
ponctuellement vers $\phi_{\beta}$, sauf peut-\^etre en $\beta$, et que
$t_n$ converge vers $t$, on v\'erifie imm\'ediatement d'apr\`es
(\ref{cob}) qu'alors $f$ satisfait l'\'egalit\'e
$\varphi_{\beta}\, f\rond T=\e{2i\pi t}f$, et l'ergodicit\'e de $T$
assure que $f$ est de module constant, non nul puisque $f$ est
suppos\'ee non nulle.

\medskip
Avec l'hypoth\`ese plus forte de K.~Merrill (\cite{mer}), avec toujours
$\sum_n |b_n|q_n\alpha_n <\infty$ mais $\sum_n |b_ns+b'_n| <\infty$,
on obtient la convergence de la s\'erie $\sum\|1-\theta_n\|_1$, donc
la convergence de $(f_n)$ dans $L^1(\M T)$, moyennant un choix des
constantes $c_n$. On peut se restreindre \`a la sous-suite des indices
pour lesquels $\tilde\theta_n$ est non nulle, c'est-\`a-dire $b_n\neq
0$ ou $b'_n\neq 0$. Avec la premi\`ere hypoth\`ese, la somme des
mesures des tours mineures correspondantes est finie et on a aussi
$\sum_n \lambda(\C I_n)<\infty$, donc il suffit d'estimer
l'int\'egrale de $|1-\theta_n|$ sur la tour principale. En choisissant
les constantes $c_n$ de fa\c{c}on que $\tilde \theta_n$ s'annule en un
point de $B'_n$, on a $|\tilde\theta_n|\leq|b_ns+b'_n|$ sur $B'_n$
d'apr\`es la propri\'et\'e (ii) et
$\|\tilde\theta_n\|\leq|b_ns+b'_n|(1+q_n\alpha_n)$d'apr\`es
(\ref{invariant}) sur $\C T'_n$. on obtient gr\^ ace \`a la seconde
condition la convergence annonc\'ee.

\medskip

Dans notre cas, il faut \^etre plus pr\'ecis et on montrera seulement
qu'il existe une fonction $g\in L^2(\M T)$ telle que $\lims| \int_{\M
T} gf_n d\lambda|>0$, pour prouver que $(f_n)$ admet une valeur
d'adh\'erence faible non nulle. On peut toujours n\'egliger
l'int\'egrale sur le compl\'ementaire de la tour principale, qui reste
de masse sommable. Sur la tour principale, on a une propri\'et\'e
d'ind\'ependance asymptotique de $f_n$ et $\theta_n$, du fait que
$f_n$ est pratiquement constante sur les \'etages et que $\theta_n$
est presque invariante par $T$. Plus pr\'ecis\'ement, on montrera que 
si $g$ est une fonction constante sur les \'etages de $\C T_n$
l'int\'egrale de $f_{n+1}$ est proche  du produit des int\'egrales de
$f_n$ et $\theta_n$. 
Il s'ensuit que la condition $\lim |\int f_nd\lambda|>0$ est
v\'erifi\'ee d\`es que le produit $\prod \inv{\lambda(B'_n)}|\int_{B'_n}\theta_nd\lambda|$ est
non nul.  
Comme d'apr\`es (\ref{invariant}), l'\'ecart des valeurs de
$\theta_n$ entre les \'etages de la tour d'ordre $n$ est encore
sommable et il suffit donc de montrer la convergence de la s\'erie
$1-\inv{\lambda(B'_n)}|\int_{B'_n}\theta_nd\lambda|$. Le calcul sur
$B'_n$ est parfaitement explicite, et la valeur trouv\'ee est bien de
l'ordre de $(b_ns+b'_n)^2$ qui est suppos\'ee sommable.

\subsubsection{Pourquoi ces conditions sont-elles n\'ecessaires?}\label{cnidee} 

Il s'agit maintenant de comprendre pourquoi les conditions  sont
\'egalement n\'ecessaires. 
On suppose donc l'existence d'une solution \`a l'\'equation
(\ref{eqbis}), et on note $\varphi_\beta=\exp{2i\pi s\phi_\beta}$. 
On a d\'eja vu au d\'ebut du chapitre que n\'ecessairement 
 $\beta$ et $t$ se d\'ecomposent sous la forme (\ref{ostro}),
et on peut donc reprendre les notations et
hypoth\`eses du paragraphe \ref{transfert}.
\\
Pour tout $n$, l'\'equation de cohomologie pour $f$ et pour $f_n$
donne en dehors de l'intervalle $[\beta_n,\beta[$ l'\'egalit\'e
$$(\tilde f_n-\tilde f)\circ T=(\tilde f_n-\tilde
f)+(t_n-t)-s(\beta-\beta_n)  \mod 1. $$
De fa\c con analogue au  paragraphe \ref{transfert}, 
on peut donc affirmer que~:
\begin{itemize}
\item
En dehors de l'ensemble $\cup_0^{q_n-1}T^{-j}[\beta_n,\beta[$, qui
est de mesure petite (car $(\beta_n-\beta)$ est n\'egligeable devant
$\alpha_{n-1}$), si $f_n$ est proche de $f$ sur la base de la tour
majeure d'ordre $n$, alors elle approche $f$ sur n'importe quel autre
\'etage \`a $2\pi q_n(\|t-t_n\|+|s|\|\beta-\beta_n\|)$ pr\`es qui
reste n\'egligeable.
\item
$(\tilde f-\tilde f_n)$ est pratiquement constante sur les \'etages
de la tour majeure d'ordre $n$ repr\'esent\'ee figure~\ref{2tours}.
\end{itemize}
Ces propri\'et\'es permettent de montrer que $f$ est la limite forte
$L^2$ de la suite des fonctions $(f_n)$  ( voir proposition
\ref{l2}). \\
Il reste enfin \`a estimer le comportement de la suite $(\theta_n)$.
Comme la suite $(f_n)$ converge dans $L^2$, il est n\'ecessaire que
$\|1-\theta_n\|_2$ soit petite. Par suite, comme sa pente est
constante sur les \'etages de la tour principale, et que sa valeur
moyenne est \`a peu pr\`es la m\^eme  sur chaque \'etage (par
(\ref{invariant})), on en d\'eduit que sa variation totale sur un
\'etage doit \^etre petite~: ceci signifie que
$(sb_n+b'_n)q_n\alpha_{n-1}\rightarrow 0$, et on retrouve la seconde
condition obtenue par Veech et Stewart, 
$|sb_n+b'_n|=\|sb_n\|\rightarrow 0$ (\cite{vee}, \cite{ste}).\\
Reprenons maintenant la figure~\ref{2tours}~: les discontinuit\'es de
$\theta_n$ sont r\'eguli\`erement r\'eparties dans $\C I_n$, et les
sauts y sont de taille constante $-s$. Pour que le produit des
$\theta_n$ converge, il est n\'ecessaire que les sauts ne se
propagent pas~: autrement dit il faut que $\lambda(\lims \C I_n)=0$.
Comme les ensembles $\C I_n$ sont \`a peu pr\`es ind\'ependants, on
obtient  comme condition n\'ecessaire la convergence de la s\'erie
$(b_nq_n\alpha_n)$. \\
Enfin, la derni\`ere condition provient comme pr\'ec\'edemment de
l'``ind\'ependance'' entre $f_n$ et $\theta_n$, qui permet d'affirmer
que, comme $\lim|\int_\M T f_nd\lambda|>0$, et que
$f_{n+1}=f_n\theta_n$, on a encore $\prod |\int_\M T
\theta_nd\lambda|>0$. On obtient que la s\'erie $(1-|\int_\M T
\theta_nd\lambda|)$ converge, d'o\`u la derni\`ere condition
n\'ecessaire $\sum |sb_n+b'_n|^2<\infty$.

\subsection{Les tours de la rotation}

Nous pr\'ecisons dans ce paragraphe la construction des tours qui
interviennent tout au long de notre travail. Les choix des tours
repr\'esent\'ees plus loin varient selon les valeurs de $n$ et de la
suite $(k_n)$~: il ne s'agit que d'unifier les discussions lors de la
preuve du th\'eor\`eme \ref{marche}. Par confort, le lecteur pourra se
restreindre \`a l'une des repr\'esentations donn\'ees ci-apr\`es (on
pourra par exemple supposer $n$ pair, $b_n$ et $k_n$ positifs).

\subsubsection{Repr\'esentation sur $\M T\times \M Z$}\label{releve} 

\begin{figure}[h]
\begin{center}\caption{Domaines fondamentaux d'ordre $n$.}\label{dn} 
\input{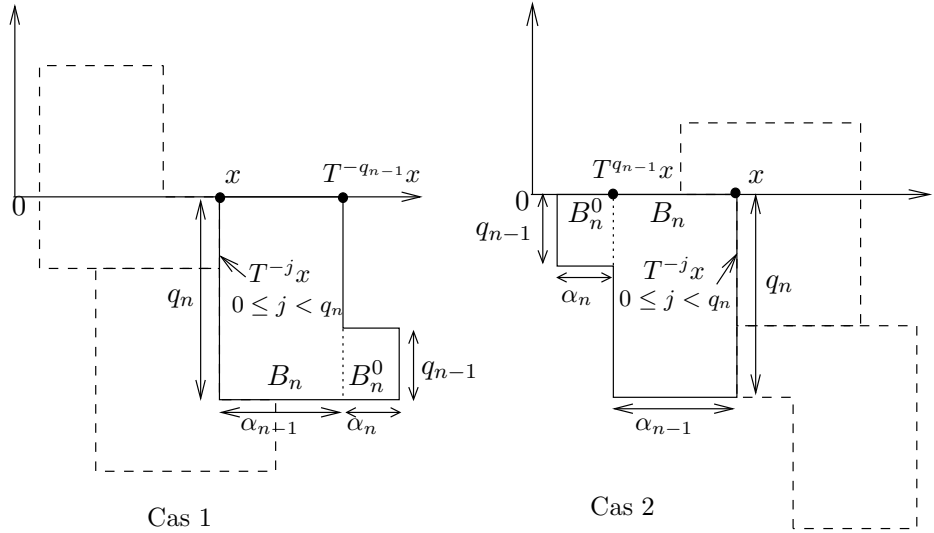}
\end{center}\end{figure}
Nous rappelons bri\`evement les \'el\'ements n\'ecessaires \`a la
suite, 
pour une description plus d\'etaill\'ee, nous r\'ef\'erons \`a
\cite{mel4}.
On d\'efinit sur  $\M T\times \M Z$ la relation d'\'equivalence
$(Tx,n)\sim (x,n+1)$~; alors $\M T$ est isomorphe \`a l'espace
quotient $\M T\times \M Z/\sim$. 
Ceci revient \`a identifier la rotation de $\alpha$ sur $\M T$ 
avec, ou bien 
la translation verticale de $(0,1)$, ou bien la translation
horizontale de $(\alpha, 0)$ dans $\M T\times \M Z$.\\
On choisit maintenant des domaines fondamentaux li\'es \`a la
fraction continue de $\alpha$ de la fa\c con suivante~:
pour tout $n$ et pour tout $x \in \M T$, la partition d\'etermin\'ee
par les
points $(T^{-j}x)_{0 <j\leq q_n}$ se d\'ecompose en  deux tours
d'intervalles~: on lui associera l'un ou l'autre des domaines
fondamentaux 
repr\'esent\'es pour $n$ pair sur la figure~\ref{dn}. 
Ils peuvent s'\'ecrire, aux translations par $(\alpha,-1)$ pr\`es,
sous la forme~: 
\\
$[x,T^{-q_{n-1}}x[\times\{- q_n+1,..,0\} \cup
[T^{-q_{n-1}}x,T^{-q_{n-1}+q_n}x[\times \{-q_n+1,..,-q_n+q_{n-1}\}$
(cas 1),
\\
$]T^{q_{n-1}}x,x]\times\{- q_n+1,..,0\} \cup 
]T^{-q_n+q_{n-1}}x,T^{q_{n-1}}x[\times\{-q_{n-1}+1,..,0\}$ (cas 2).
\\  
Remarquons que le cas 2 s'obtient \`a partir du cas 1 en rempla\c
cant $T$ par $T^{-1}$. Nous adopterons pour la suite les
terminologies suivantes~:
\begin{itemize}
\item
La tour majeure associ\'ee au domaine fondamental d'ordre $n$~: c'est
la tour 
de hauteur $q_n$ et de base $B_n$ \'egale \`a
$T^{-q_n+1}[0,T^{-q_{n-1}}0[$ dans le cas 1. Dans le cas 2, nous
appellerons encore base l'intervalle  $B_n$ d\'efini par
$]T^{q_{n-1}}0,0]$ (qui est bien la base de la tour majeure pour la
transformation $T^{-1}$).  
\item
La tour mineure associ\'ee au domaine fondamental d'ordre $n$ est la
tour de hauteur $q_{n-1}$ et de base $B_n^0$ \'egale \`a
$T^{-q_n+1}[T^{-q_{n-1}}0, T^{-q_{n-1}+q_n}0[$ dans le cas 1 et
$]T^{q_{n-1}-q_n}0,T^{q_{n-1}}0]$  dans le cas 2 (selon les m\^emes
conventions que ci-dessus).
\end{itemize}
On notera que les bases des  tours mineure et majeure sont deux
intervalles admettant une extr\'emit\'e commune.

\subsubsection{Une suite de domaines fondamentaux
privil\'egi\'ee.}\label{df}
\begin{figure}[h]
\begin{center}
\caption{Repr\'esentation de $\C D_n(\beta)$ (cas $n$
pair).}\label{tourpropre}
\input{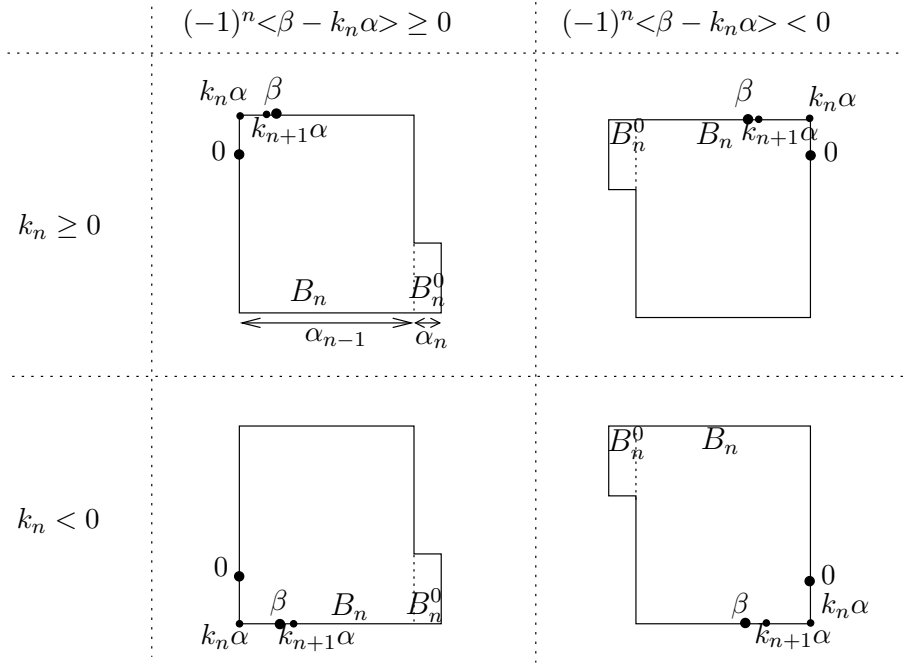}
\end{center}\end{figure}
Soit $\beta \in \M T$. On suppose que $\|\beta q_n\|\rightarrow 0$ et
on reprend les notations (\ref{ostro}) pour $\beta$. 
%$(b_n)$ est une suite d'entiers v\'erifiant $b_n/a_{n+1} \rightarrow
%0$.
On choisit une suite de domaines fondamentaux qu'on note $(\C
D_n(\beta))$ et qui v\'erifie les conditions suivantes~:
\begin{itemize}
\item
$\C D_n(\beta)$ est un des domaines fondamentaux associ\'e \`a la
partition d\'etermin\'ee par les points $(T^j0)_{k_n\leq j<q_n+k_n}$
si $k_n<0$ et $(T^j0)_{k_n-q_n<j\leq k_n}$ sinon.
\item
$\C D_n(\beta)$ est repr\'esent\'e selon le cas 1 de la figure
\ref{dn} si 
$\<\alpha q_n><\beta-k_n\alpha>\geq 0$ et selon le cas 2 sinon.
\end{itemize}
Nous obtenons pour $\C D_n(\beta)$ l'un des cas de la figure
\ref{tourpropre}.
On remarquera ais\'ement que la suite des partitions associ\'ee \`a
$(\C D_n(\beta))$ se raffine, et bien sur qu'elle engendre la tribu
bor\'elienne.
Enfin, pour tout $n $ suffisamment grand (de sorte que
$|b_n|<a_{n+1}/2$ par exemple), lorsque $k_n \neq k_{n+1}$,
$\|k_n\alpha-k_{n+1}\alpha\|$ est inf\'erieur \`a
$\alpha_{n-1}$, et $\<\beta-k_n\alpha\>$ et
$\<k_{n+1}\alpha-k_n\alpha\>$ sont de m\^eme signe. On en d\'eduit
que les points 
$k_n\alpha$, $k_{n+1}\alpha$ et $\beta$  sont toujours sur le m\^eme
\'etage, 
soit le premier soit le dernier, de la tour majeure d'ordre $n$. 
Nous retrouverons ces domaines fondamentaux dans les preuves des
th\'eor\`emes \ref{marche} et \ref{cn}.

\subsubsection{$\|\beta q_n\|$ dans un domaine fondamental d'ordre
$n$.}
\label{betaqn}

\begin{figure}[h]
\caption{R\'epartition de $\phi_\beta^{(q_n)}$ dans $\C
D_n(0)$.}\label{bqn}
\begin{center}
\input{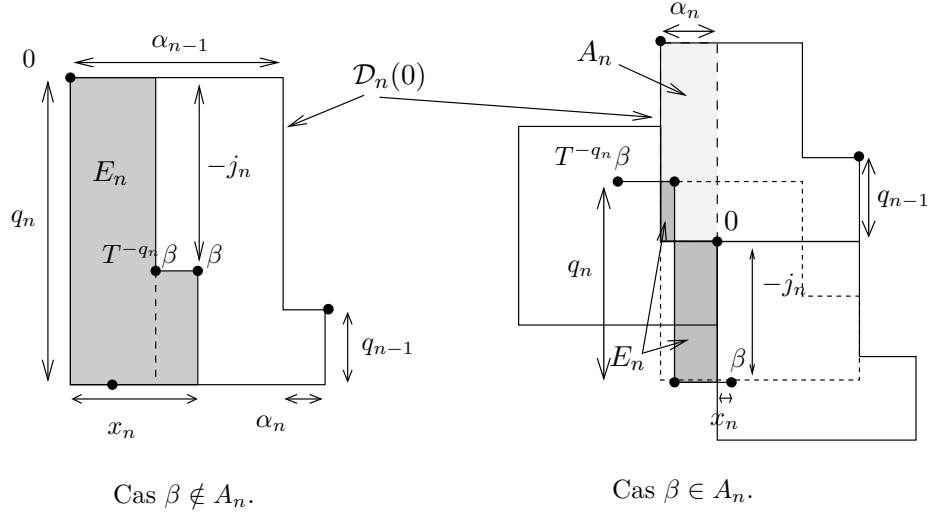} 
\end{center}
\end{figure}

Soit $\beta$ dans $\M T$. Pour tout $n$, ce point admet un unique
repr\'esentant (encore appel\'e $\beta$) dans le domaine fondamental
$\C D_n(0)$. Ses  coordonn\'ees $(x_n,j_n)$ dans $\C D_n(0)$
satisfont 
$\beta= j_n\alpha + (-1)^nx_n $ (cf figure~\ref{bqn}).
Notons \'egalement que 
$$\beta q_n=j_n\alpha q_n +(-1)^nx_nq_n=(-1)^n(j_n
\alpha_n+q_nx_n).$$
On s'int\'eresse maintenant \`a la distribution
de $\phi_\beta^{(q_n)}$ dans $\M T$~:
$\phi_\beta^{(q_n)}$ poss\`ede $2q_n$ discontinuit\'es $(T^{-j}0)_{0
\leq j<q_n}$ et $(T^{-j}\beta)_{0 \leq j <q_n}$, qui  se
repr\'esentent facilement
 dans $\C D_n(0)$. 

Comme les \'etages de $\C D_n(0)$ sont des intervalles, 
$\phi_\beta^{(q_n)}$ est constante sur les segments des \'etages de
$\C D_n(0)$ sans discontinuit\'e, et constante
sur les translat\'es de ces segments d\`es qu'ils ne contiennent pas
$[T^{-q_n}\beta,\beta[$ ni $[T^{-q_n}0, 0[$ 
(car $\phi_\beta^{(q_n)}\rond T = \phi_\beta^{(q_n)} +\phi_\beta\rond
T^{q_n}-\phi_\beta$).  En remarquant enfin que les sauts de
$\phi_\beta$
sont constants et valent $\pm 1$, on constate que
$\phi_\beta^{(q_n)}$ ne prend que 2 ou 3 valeurs qui se r\'epartissent
selon les domaines 
repr\'esent\'es sur la figure~\ref{bqn} (cas $n$ pair).
Soit
 $E_n=\{\phi_\beta^{(q_n)}\neq \phi_\beta^{(q_n)}(\beta)\}$
(l'ensemble hachur\'ee en gris sur la figure~\ref{bqn}). La condition (\ref{rigide})  montre que les seules valeurs d'adh\'erence possibles de $\lambda(E_n)$ sont $0$ ou $1$.\\
Posons  $A_n$ l'ensemble de $\C D_n(0)$ des $\beta$ tels que $x_n\leq 0$, alors $\|\beta q_n\|$ et $\lambda(E_n)$ sont li\'es par les remarques suivantes~:
\\
Lorsque $\beta \notin A_n$, $\phi^{(q_n)}$ n'a que 2 valeurs et  
$\lambda(E_n)=\alpha_nj_n+x_nq_n=(-1)^n\beta q_n \hbox{ mod }1$ ce
qui donne
$$\|\beta q_n\|=\min(\lambda(E_n), 1-\lambda(E_n)).$$
Lorsque $\beta \in A_n$, on a  
$\lambda(E_n)= -x_n(q_n-j_n)+(\alpha_n+x_n)j_n$.
ALors  on voit ais\'ement que 
$|(-1)^n \beta q_n|=|x_n(q_n-j_n)+(\alpha_n+x_n)j_n|\leq \lambda(E_n)$.
D'autre part, \`a l'aide de la repr\'esentation g\'eom\'etrique des 2 quantit\'es, on v\'erifie facilement que 
$\lambda(E_n)-|<\beta q_n>|=2 \min(-x_n(q_n-j_n),j_n(\alpha_n+x_n)\leq 2(\alpha_nq_n-\lambda(E_n))\leq 2(1-\lambda(E_n))$.   
On obtient donc
$$3\lambda(E_n)-2\leq |\{\beta q_n\}|\leq \lambda(E_n).$$
Par cons\'equent, comme
$(\min (\lambda(E_n), 1-\lambda(E_n)))$ converge vers $0$,  on en
d\'eduit aussi que  $(\|\beta q_n\|)$
converge vers $0$.

\subsection{Quelques relations utiles}\label{preli}
Nous aurons besoin, pour la preuve du th\'eor\`eme~\ref{marche},
d'exprimer la fonction de transfert solution de l'\'equation
(\ref{eqbis}) comme limite forte dans $L^2$ d'une suite de fonctions
de transfert bien approch\'ees par leurs projections sur des suites
de tours associ\'ees \`a la rotation.  
Ce paragraphe a pour objectif d'\'etablir des conditions suffisantes
sur 
une suite de fonctions de transfert pour qu'elle admette une limite
forte 
$L^2$. Les conditions obtenues ne sont pas sp\'ecifiques aux
rotations 
irrationnelles,  et on pourra remarquer que  tous les r\'esultats
\'enonc\'es ici sont valables pour n'importe quelle transformation 
$T$ de $\M T$, 
bijective, bimesurable et pr\'eservant $\lambda$.\\
On suppose dans les lemmes qui suivent que $\varphi$ est 
une fonction de module 1. 

S'il existe $\theta$ de module 1 et $\delta$ solutions de l'\'equation
$\DSM{\theta\rond T \over \theta} = \e{i\delta}\varphi$, alors pour
tout
ensemble mesurable $B$ on peut \'ecrire
\begin{eqnarray}\label{etage}
\int_{T^jB}\theta d\lambda = \int_B\theta d \lambda +
{\Oun}\left(\delta j
\lambda(B)\right)+ {\Oun}\left(2\lambda(\{\varphi^{(j)}\neq
1\}\cap B)\right).
\end{eqnarray}

\begin{lemme}\label{l1}
Supposons $\theta$ et $\delta$ solutions de l'\'equation
pr\'ec\'edente.
Soit $(T^jB)_{0 \leq j <h}$ une tour telle que
 $\varphi=1$ sur $\C T= \cup_{0 \leq j <h}T^jB$.
Pour toute fonction $g$ de module inf\'erieur \`a 1,  constante sur
chacun des \'etages de la tour et pour toute fonction $f$  de module
1, on a  
\begin{eqnarray}\label{tour}
\int_{\C T}g f\theta d\lambda =\frac{1 }{ \lambda(B)} \int_{\C T}g
fd\lambda
\int_{B}\theta d\lambda + 
\lambda(\C T){\Oun}\left(\inv{h}\sum_{0 \leq j <h}V_{T^jB}(f) 
+ \frac{\delta}{ \lambda(B)} \right)\; .
\end{eqnarray}
\end{lemme}

\dem On commence par d\'ecomposer le terme de gauche sur les
\'etages~:
\begin{multline*}
\int_{\C T}g f \theta d\lambda
=  \sum_{0 \leq j <h}\int_{T^jB}g f \theta d\lambda
\\
= \sum_{0 \leq j <h}\int_{T^jB}g\theta d \lambda
\inv{\lambda(B)}\int_{T^jB}f d\lambda +
\sum_{0 \leq j <h}\int_{T^jB}g \theta
\left(f- \inv{\lambda(B)}\int_{T^jB}f d\lambda \right)d\lambda.
\end{multline*}
Pour tout $0\leq j <h$, comme $g$ est constante sur $T^jB$ et que
$\varphi^{(j)}=1$ sur $B$, on obtient en utilisant  (\ref{etage}) que
~:
\begin{eqnarray*}
\int_{T^jB}g\theta d \lambda\inv{\lambda(B)}\int_{T^jB}f d\lambda 
&=& \inv{\lambda(B)}\int_{T^jB} gf
d\lambda \int_{T^jB}\theta d\lambda \\
&=&
\int_{T^jB}g f d\lambda
\left( \inv{\lambda(B)}\int_{B}\theta d\lambda + 
{\Oun}(\delta h)\right).
\end{eqnarray*}
D'autre part on a aussi, en remarquant que $|\theta g|\leq 1$,  
$$
\left|\int_{T^jB}g \theta
\left(f- \inv{\lambda(B)}\int_{T^jB}f d\lambda \right)d\lambda\right|
\leq  V_{T^jB}(f)\lambda(B)= \frac{\lambda(\C T)}{h}  V_{T^jB}(f).$$
Comme $|fg|\leq 1$, on obtient  par sommation l'\'egalit\'e
annonc\'ee.
\findem

\begin{lemme}
Soit $(T^jB)_{0 \leq j <h}$ une tour, on a pour tout entier $k$ 
l'in\'egalit\'e~:
\begin{eqnarray}\label{rem}
\lambda(B \cap \{\varphi^{(k)}\neq 1\}) \leq
\left[\frac{|k|}{h}+1\right] 
\lambda(\{\varphi\neq 1\}).
\end{eqnarray}
\end{lemme}

\dem
Prenons $k>0$ par exemple.
Il suffit de remarquer que $\{\varphi^{(k)}\neq 1\}\subset \cup_{0
\leq
j<k}\{\varphi\rond T^j\neq 1\}$, en posant $A= \{\varphi \neq 1\}$ 
il vient alors naturellement~:
\begin{eqnarray*}
\lambda(B\cap \{\varphi^{(k)}\neq 1\}) & \leq & \sum_{0 \leq j
<k}\lambda(B\cap
T^{-j}A)
\leq  
\sum_{0 \leq j <k}\lambda(T^j B\cap A)\\ 
&\leq &
\sum_{0 \leq i <k/h}\sum_{0 \leq j <h}\lambda(T^{ih+j} B\cap A)\\
&\leq &
\sum_{0 \leq i <k/h}\lambda(T^{ih}(\cup_{0 \leq j <h}T^j B)\cap A)
\leq 
\left[\frac{k}{h}+1\right]\lambda(A).
\end{eqnarray*}
\findem

\begin{lemme}\label{convergefort}
Pour tout entier $n$ on suppose qu'il
existe deux fonctions de module 1, $f_n$ et $\gamma_n$, 
et un r\'eel $x_n$ v\'erifiant
$$f_n \rond T= f_n \e{i x_n} \gamma_n\;. $$
Si on peut trouver $(T^jB_n)_{0 \leq j <h_n}$ une suite de  tours et 
$N$ un entier ind\'ependant de $n$ de sorte 
 qu'on ait~:
\begin{enumerate}
\item{ $\lambda(\cup_{0 \leq j <Nh_n}T^jB_n)\tend{n
\rightarrow\infty} 1$,}
\item{$x_n h_n \tend{n \rightarrow \infty} 0$,}
\item{$\lambda(\gamma_n \neq 1)h_n \tend{n \rightarrow \infty}0$,}
\item{$V_{B_n}(f_n) \tend{n\rightarrow \infty}0$.}
\end{enumerate}
Alors, il existe une suite $(c_n)$ de constantes de module 1 telles
que
$\|c_nf_n-1\|_2\tend{n \rightarrow \infty} 0$.
\end{lemme}

\dem
Soit $n \in \M N$ et $c_n$ une constante de module 1, on a
l'in\'egalit\'e
$$
\int_{\M T}|1-c_nf_n|^2 d\lambda 
 \leq
\sum_{0 \leq j <Nh_n}\int_{T^{j}B_n} |1-c_nf_n|^2d\lambda +
4\left(1-\lambda(\cup_{0 \leq j<Nh_n}T^jB_n)\right)\; .$$
En remarquant que 
$|f_nc_n-1|^2= 2-f_nc_n-\overline{f_nc_n}$, et que $\overline {f_n}$
est 
une fonction de transfert associ\'ee \`a $\overline{\gamma_n}$ et
$-x_n$, 
on obtient \`a l'aide de (\ref{etage}), pour tout $0 \leq j<Nh_n$~:
\begin{flalign*}
\int_{T^jB_n}|1-c_nf_n|^2 d\lambda 
&\leq 
 \int_{B_n}|1-c_nf_n|^2 d\lambda +2  |x_n| j\lambda(B_n) + 
4 \lambda\left(B_n \cap \{\gamma^{(j)}_n \neq 1\}\right)\\
&\leq  
\int_{B_n}|1-c_nf_n|^2d\lambda +2N|x_n| +
4 N \lambda(\{\gamma_n \neq 1\})\quad \hbox{par } (\ref{rem}).
\end{flalign*}
D'autre part on a 
$$V_{B_n}(f_n)=\inv{\lambda(B_n)}\int_{B_n}\left(\inv{\lambda(B_n)}
\int_{B_n}|f_n(x)-f_n(y)|dx\right)dy\; . $$ 
On peut trouver $y_n$ tel que si $\overline c_n= f_n(y_n)$, on ait
$\inv{\lambda(B_n)}\int_{B_n}|f_n-\overline{c_n}|d\lambda\leq
V_{B_n}(f_n)$.
Par cons\'equent, pour une constante $c_n$ ainsi choisie on obtient
$$\int_{B_n}|1-c_nf_n|^2d\lambda \leq 2 \lambda(B_n)V_{B_n}(f_n)\;,$$
et en sommant dans les in\'egalit\'es pr\'ec\'edentes il vient
\begin{multline*}
\int_{\M T}|1-c_nf_n|^2d\lambda \leq
2NV_{B_n}(f_n)+2N^2h_n|x_n|+4N^2h_n\lambda(\{\gamma_n\neq 1\})\\
+4(1-\lambda(\cup_{0}{Nh_n}T^jB_n)).
\end{multline*}
Pour une telle suite $(c_n)$, les hypoth\`eses du lemme 
permettent donc clairement d'assurer la convergence forte dans $L^2$
de 
$(c_nf_n)$ vers 1.
\findem

\begin{lemme}\label{proj}
Soit $(T^jB_n)_{0\leq j<h_n}$ une suite de tours v\'erifiant
$\limi h_n\lambda(B_n) >0$. Alors pour toute fonction mesurable $f$ de
module 1, les deux propositions sont
\'equivalentes~:
\begin{itemize}
\item{
La distance dans $L^2$ de la restriction de $f$ \`a
$\cup_{j<h_n}T^jB_n$ sur le
sous-espace engendr\'e par les fonctions indicatrices des \'etages de
$(T^jB_n)_{j<h_n}$ converge vers 0.}
\item{
La moyenne des variations moyennes de $f$ sur $T^jB_n$ pour $0\leq
j<h_n$
tend vers 0.}
\end{itemize} 
\end{lemme}

\dem
En consid\'erant la projection orthogonale dans $L^2$ de $f$ au
sous-espace
$\C H_n$ engendr\'e par les $(\unde{T^jB_n})_{j<h_n}$, la distance de
la 
restriction de $f$ \`a $\cup_{j<h_n}T^jB_n$  au sous-espace $\C H_n$
peut 
s'\'ecrire  
\begin{flalign*}
d^2(f\unde{\cup_{j<h_n}T^jB_n}, \;&\C H_n)\\
 &= \sum_{0 \leq j
<h_n}\left(\int_{T^jB_n}|f|^2d\lambda-\inv{\lambda(B_n)}
\left|\int_{T^jB_n}fd\lambda\right|^2\right) \\
 &= 
\sum_{0\leq j <
h_n}\inv{\lambda(B_n)}\iint_{(T^jB_n)^2}(1-f(x)\overline{f(y)})dxdy\\
 &=  
\frac{h_n \lambda(B_n)} {2}\inv{h_n}\sum_{0 \leq j <h_n}
\inv{\lambda(B_n)^2}\iint_{(T^jB_n)^2}|f(x)-f(y)|^2dxdy\; .
\end{flalign*}
Notons pour tout bor\'elien $B$ de mesure non nulle, 
$W_B(f)=\inv{\lambda(B)^2}\iint_{B^2}|f(x)-f(y)|^2dxdy$. 
Comme $f$ est une fonction de module 1, on a les in\'egalit\'es~:
$$V_B(f)^2\leq W_B(f)\leq 2V_B(f).$$
Par cons\'equent,  $(\inv{h_n}\sum_{j<h_n}W_{T^jB_n}(f))$ converge
vers 
0 si et seulement si la suite $(\inv{h_n}\sum_{j<h_n}V_{T^jB_n}(f))$
converge 
vers 0. 
Par ailleurs on avait 
$$d^2(f\unde{\cup_{j<h_n}T^jB_n}, \;\C H_n)=
\frac{h_n \lambda(B_n)} {2}\inv{h_n}\sum_{0 \leq j
<h_n}W_{T^jB_n}(f).$$ 
Comme $\limi \lambda(B_n)h_n>0$ et $h_{n}\lambda(B_{n})\leq 1$, on
obtient bien l'\'equivalence annonc\'ee.  
\findem

\break
\part{L'\'equation de cohomologie pour les fonctions en escalier.}

\section{Cas des fonctions \`a deux discontinuit\'es.}

Il s'agit dans cette partie de prouver le th\'eor\`eme~\ref{marche}~:
la
preuve se d\'ecompose en plusieurs \'etapes, dont la premi\`ere est
l'\'etude des cas triviaux trait\'ee aux paragraphes \ref{cobord} et
\ref{transfert}.
Le paragraphe qui suit est consacr\'e \`a la recherche de conditions
n\'ecessaires \`a l'existence d'une solution \`a l'\'equation
(\ref{eqbis}). 
On y montre que la suite des fonctions de
transfert d\'efinie au paragraphe \ref{transfert}  converge fortement
vers une 
fonction de transfert associ\'ee \`a $\varphi_\beta$. Comme ces
fonctions sont proches de fonctions constantes sur les \'etages des
tours de la rotation, 
on \'etablit \`a l'aide des
propri\'et\'es d'ind\'ependance asymptotique de ces tours, une
convergence plus forte qui aboutit, apr\`es un calcul explicite, 
aux conditions du th\'eor\`eme~\ref{marche}.
Le second paragraphe montre le caract\`ere suffisant de ces
hypoth\`eses.

\subsection{Conditions n\'ecessaires.}

Soit $\beta \in \M T$ et $s \in \M R \backslash \M Z$, on note
$\varphi_\beta=\exp{(2i\pi s\phi_\beta)}$. 
On suppose qu'il existe $t \in \M T$ et $f \in L^\infty(\M T)$ de
module 1
tels que $$\varphi_\beta= \e{2i\pi t} {f \over f \rond T}\; .$$
On note  $(q_n)$ la suite des d\'enominateurs de la fraction continue 
associ\'ee \`a $\alpha$.
\\
Comme $T^{q_n}\to I$ et que  $\e{-2i\pi q_n
t}\varphi^{(q_n)}_\beta=f/f\circ T^{q_n}$, on doit avoir la
condition (\ref{rigide})~:
$$ \e{-2i\pi t q_n}\int_\M T \varphi_\beta^{(q_n)}d\lambda \tend{n
\rightarrow \infty} 1.$$
On obtient en particulier que  $|\int_{\M T} \varphi^{(q_n)}_\beta
d\lambda |\rightarrow 1$, d'o\`u on d\'eduit que 
$\|\beta q_n\|\rightarrow 0$ (cf paragraphe 
\ref{betaqn} ). 
Estimons maintenant la valeur principale de $\varphi_\beta^{(q_n)}$.
On a la propri\'et\'e suivante
\begin{lemme}\label{phiqn}
Pour tout $\beta \in \M T$ tel que  
$\|\beta q_n\|\tend{n \rightarrow \infty}0$,  la suite des fonctions 
$(\phi_\beta^{(q_n)})_n$ converge vers 0 en mesure. 
\end{lemme}
\csq
Il en r\'esulte imm\'ediatement la convergence de
$(\varphi_\beta^{(q_n)})$ 
vers 1 dans $L^1$.\\
\dem
Comme $\|\beta q_n\|\tend{n \rightarrow \infty}0$, on peut trouver
une suite 
d'entiers $(k_n)$ v\'erifiant $\lim k_n\alpha= \beta \mod 1$, avec
$k_n/q_n\rightarrow 0$ et $\|\beta -k_n\alpha\|q_n \rightarrow 0$.  
On peut donc approcher $\phi_\beta$ par
$\phi_{k_n\alpha}$. On obtient \`a l'aide du  paragraphe \ref{cobord}
$$\phi_\beta^{(q_n)}=\phi_{k_n\alpha}^{(q_n)}+
\phi_{\beta-k_n\alpha}^{(q_n)}\circ T^{-k_n} =  \omega_{k_n}-
\omega_{k_n}\circ T^{q_n}+
\phi_{\beta-k_n\alpha}^{(q_n)}\circ T^{-k_n}.$$
On sait que $\phi_{\beta-k_n\alpha}= -\<\beta-k_n\alpha\>$ en dehors
d'un intervalle de longueur $\|\beta-k_n\alpha\|$.  
Par cons\'equent $\phi_{\beta-k_n\alpha}^{(q_n)}= -
q_n\<\beta-k_n\alpha\>$ 
sauf sur un ensemble de mesure inf\'erieure \`a
$q_n\|\beta-k_n\alpha\|$. 
Comme la suite $(\|\beta-k_n\alpha\|q_n) $ converge 
vers $0$, on trouve que $(\phi_{\beta-k_n\alpha}^{(q_n)})$ converge 
vers 0 en mesure.
\\ 
D'autre part, $\omega_{k_n}$ est de pente constante $k_n$, et
poss\`ede $k_n$ discontinuit\'es, aux points $(T^j0)_{0<j\leq k_n}$
(en supposant par exemple $k_n>0$)~; alors $\omega_{k_n}-\omega_{k_n}\circ
T^{q_n} $ est  constante et \'egale \`a
$k_n(-1)^n\alpha_n$, sauf sur des intervalles de longueur 
$\alpha_n$  de la forme $[T^j0, T^{j-q_n}0]$ pour $0 <j\leq k_n$.
Comme $k_n\alpha_n\tend{n\rightarrow \infty}0$,  on en d\'eduit  que
$\omega_{k_n}-\omega_{k_n}\circ
T^{q_n} $ converge  vers 0 en
mesure.  
\findem

D'apr\`es le lemme~\ref{phiqn}, la condition (\ref{rigide}) entra\^\i
ne donc la convergence de  $(\|tq_n\|)$ vers $0$. On peut alors
\'ecrire 
$$\beta = \sum_0^\infty b_n q_n \alpha  \qquad \hbox{et }\qquad t =
\sum_0^\infty b'_nq_n\alpha,$$
o\`u $(b_n)$ et $(b'_n)$ sont des suites d'entiers v\'erifiant $\lim
b_n/a_{n+1}= \lim b'_n/a_{n+1} =0$. Nous noterons \'egalement comme
dans la pr\'esentation (cf \ref{idee})~:
$$
\left\lbrace{ 
\begin{array}{lcl}
\beta_n=k_n\alpha, & & k_n=\sum_0^{n-1}b_jq_j,\\
t_n= k'_n\alpha &\hbox{et} & k'_n=\sum_0^{n-1}b'_jq_j.
\end{array}
}\right.$$
On rappelle que $k_n$ et $k'_n$ sont infiniment petits devant $q_n$,
et
que $\|\beta-\beta_n\|$ et $\|t-t_n\|$ sont infiniment petits devant 
$\alpha_{n-1}$.
\\
Notons, toujours selon le paragraphe \ref{transfert} pour tout $n$, 
$f_n = \exp{(2i\pi \omega_{k_n,k'_n})}$~: cette fonction 
 est bien approch\'ee par sa projection sur les  \'etages de la tour
majeure d'ordre $n$ de la rotation, dont la mesure tend vers 1. De
plus 
les suites $(k_n\alpha)$ et $(k'_n\alpha)$ convergent assez vite vers
$\beta$ et $t$. Ceci permet de montrer le r\'esultat suivant~: 

\begin{prop}\label{l2}
S'il existe une solution \`a l'\'equation $(1)$, quitte \`a
multiplier $f_n$ par une constante de module 1, la suite $(f_n)$
d\'efinie pr\'ec\'edemment 
converge fortement vers $f$ dans $L^2$. 
\end{prop}

\dem
Les \'egalit\'es (\ref{eqbis}) et (\ref{cob}) permettent d'\'ecrire~:
$$(f \overline{f_n})\rond T=(f \overline{f_n}) 
\e{2i\pi [(t-t_n)+s(\beta-\beta_n)]}\e{-2i\pi s\unde{[\beta_n,
\beta[}}\; . $$
On  va appliquer  le lemme~\ref{convergefort} pour la suite
$(f\overline
f_n)$ en choisissant la  suite des tours majeures des domaines $\C
D_n(\beta)$ (cf figure~\ref{tourpropre}), qu'on note  $((T^jB_n)_{0
\leq j <q_n-1})_n$.
V\'erifions les hypoth\`eses du lemme~\ref{convergefort}~: 
la condition $(i)$ 
est bien r\'ealis\'ee en prenant $N=2$, et les conditions $(ii)$
et $(iii)$ r\'esultent imm\'ediatement du fait que
$(\|\beta-\beta_n\|q_n)$ et
$(\|t-t_n\|q_n)$ convergent vers 0. 
Enfin, il reste \`a v\'erifier que $V_{B_n}(f \overline f_n)
\rightarrow 0$.
Comme
les discontinuit\'es de $f_n$ se trouvent aux extr\'emit\'es des
\'etages
de $(T^jB_n)_{0 \leq j <q_n-1}$, $f_n$ est continue sur chacun des
\'etages de la tour d'ordre $n$. Sur $B_n$, $\tilde f_n$ est donc
affine et de pente $sk_n+k'_n$ d'o\`u 
$V_{B_n}(f_n)\leq 2\pi|sk_n+k'_n|\alpha_{n-1} \rightarrow 0$. 
De plus 
$V_{B_n}(f\overline f_n) \leq V_{B_n}(f_n)+V_{B_n}(f)$, 
et il suffit pour montrer le r\'esultat de v\'erifier que la
variation de $f$ sur $B_n$ converge bien vers $0$.
\\
$(B_n)$ \'etant une
suite d'intervalles dont la mesure tend vers $0$, la restriction de
$f$ \`a
$\cup_{0 \leq j <q_n-1}T^jB_n$  est proche dans $L^2$ 
de sa projection orthogonale sur l'espace engendr\'e par 
$(\unde{T^jB_n})_{0 \leq j <q_n-1}$. De plus $(q_n-1)\alpha_{n-1}\geq
1/2$,
et le lemme~\ref{proj} assure donc  que 
$$\inv{q_n-1}\sum_{0\leq j<q_n-1}V_{T^jB_n}(f) \tend{n \to \infty}0\;
.$$  
Par ailleurs $\varphi_{\beta}$ est constante sur les intervalles ne
contenant 
ni 0 ni $\beta$, elle est donc constante sur chaque $T^jB_n$, pour 
$j \in \{0,..,q_n-2\}$. Alors pour tout $x$ et $y$ dans $B_n$ on a 
$\varphi^{(j)}_\beta(x)=\varphi_\beta^{(j)}(y)$. Par cons\'equent, 
quelque soit $j \in \{0,..,q_n-2\}$, pour tous $x$ et $y$ dans $B_n$, 
\begin{eqnarray*}
|f(T^jx)-f(T^jy)| & = & 
|f(x)\e{2i\pi jt} \varphi_\beta^{(j)}(x) -f(y)\e{2i\pi j t}
\varphi_\beta^{(j)}(y)|\\
& = &
|f(x)-f(y)|.
\end{eqnarray*}
Il en r\'esulte que la variation moyenne de $f$ ne d\'epend pas de
l'\'etage de
la tour, et ceci entra\^ \i ne bien que $V_{B_n}(f) \rightarrow 0$. 
\findem

On s'int\'eresse maintenant \`a $\theta_n=f_{n+1}\overline f_n$~:
d'apr\`es le paragraphe \ref{transfert}, $\tilde \theta_n$ admet des
discontinuit\'es aux points de la forme $(T^j0)_{k_n<j\leq k_{n+1}}$
(si $b_n>0$), dont les sauts en ces points valent $-s$, et qui est
affine et de pente constante \'egale \`a $(b_ns+b'_n)q_n$ entre 2
discontinuit\'es. 
Le but de ce paragraphe est d'expliciter \`a l'aide de $\theta_n$ le
r\'esultat de la proposition pr\'ec\'edente. C'est ce calcul qui
aboutira aux conditions  du th\'eor\`eme~\ref{marche}.\\
Comme $(\|f_n-f_{n+1}\|_2)$ converge vers 0, on obtient~:
$$\|1-\theta_n\|_2 \tend{n \rightarrow \infty}0.$$ 
\\
{\bf Remarque~:}
lorsque $b_n=0$,  $\theta_n$ est simplement
une fonction propre de $T$ pour la valeur propre  $b'_nq_n\alpha$.
Comme $(\theta_n)_n$
converge vers $1$ dans $L^2$, pour $n$ assez grand $\theta_n$ est
alors
constante et \'egale \`a 1.  On en d\'eduit que pour $n$ assez grand
l'hypoth\`ese $b_n=0$ entra\^ \i ne
$b'_n=0$ et $f_{n+1}=f_n$.
Nous pourrons donc par la suite nous restreindre \`a l'ensemble des
indices $\Lambda$ d\'efini par~:
$$\Lambda=\{n, b_n\neq 0\}.$$ 

\begin{figure}
\begin{center}\caption{Les tours  de $\C D_n(\beta)$.}\label{3tours} 
\input{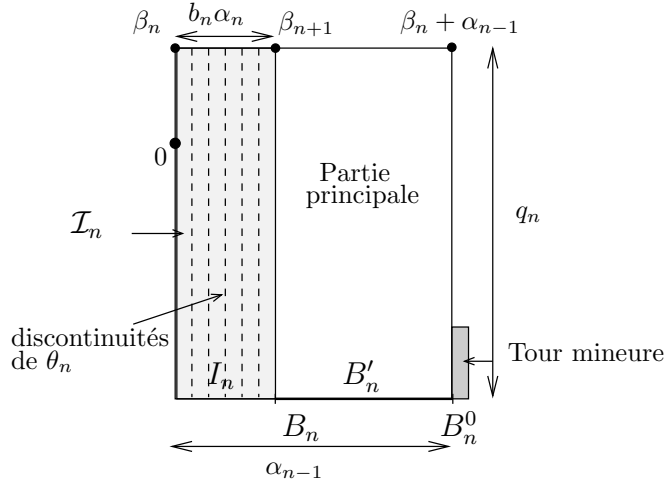}
\end{center}\end{figure}

Notons pour la suite  $\psi_n=\exp{-2i\pi s\unde{[0,b_n\<
q_n\alpha\>[}}\circ T^{-k_n}$. D'apr\`es le paragraphe
\ref{transfert}, et  par analogie avec le paragraphe \ref{preli} on
\'ecrit~:
$$\frac{\theta_n \circ T}{\theta_n}=
\e{2i\pi(sb_n+b'_n)\<q_n\alpha\>} 
\psi_n.$$
Gr\^ ace \`a la repr\'esentation de $\theta_n$ dans le domaine
fondamental $\C D_n(\beta)$ (voir figure~\ref{3tours}) on \'etablit
les propri\'et\'es  suivantes (on rappelle que les ensembles $(\C
I_n)$ sont d\'efinis par (\ref{in}), voir figure~\ref{2tours} et
\ref{3tours}).
\begin{lemme}\label{bs-b'}
Lorsque $(b_nq_n\alpha_n)$ et $(b'_nq_n\alpha_n)$ convergent vers
$0$, toutes les propri\'et\'es suivantes sont \'equivalentes ~:
\begin{enumerate}
\item
$\|1-\theta_n\|_2 \tend{n\rightarrow \infty} 0,$ 
\item
Pour toute suite d'intervalles $(J_n)_n$ de longueur
$\alpha_{n-1}$, on a 
$$\inv{\lambda(J_n)}\int_{J_n}(1-\theta_n)d\lambda \tend{n \rightarrow
\infty}0,$$
\item
$sb_n+b'_n  \tend{n \rightarrow \infty}  0$,  et $\exists x_n \in
B'_n, \; \theta_n(x_n)\tend{n \rightarrow \infty}  1$,
\item
$ \|(1-\theta_n)\unde{\C I_n^c}\|_\infty \tend{n \rightarrow \infty }
0$.
\end{enumerate}
\end{lemme}

\dem
D'apr\`es la remarque pr\'ec\'edente, il suffit de d\'emontrer le
lemme pour $n \in \Lambda$.\\
Supposons (i). Si $J_n$ est un intervalle de longueur $\alpha_{n-1}$,
alors 
$(T^jJ_n)_{0 \leq j <q_n}$ est une tour dont la mesure
$\alpha_{n-1}q_n$ est
sup\'erieure \`a 1/2. Les relations (\ref{etage}) et (\ref{rem})
permettent 
d'\'ecrire pour tout $j \in \{0,..,q_n-1\}$~:
$$\int_{J_n}\theta_n d\lambda=\int_{T^jJ_n}\theta_nd\lambda+
{\Oun} (2\pi |sb_n+b'_n|\alpha_n
q_n\lambda(J_n)+4\lambda(\{\psi_n\neq 1\})).$$
Par sommation, on obtient alors 
\begin{multline*}
q_n\left|\int_{J_n}(1-\theta_n)d\lambda\right|
 \leq  
\left|\int_{\cup_{0 \leq j <q_n}T^jJ_n}
(1-\theta_n)d\lambda\right|
+ 2\pi |sb_n+b'_n|\alpha_nq_n(q_n\alpha_{n-1}) \\
+4q_n\lambda(\{\psi_n\neq 1\}),
\end{multline*}
d'o\`u, comme  $\lambda(\{\psi_n\neq 1\})\leq |b_n|\alpha_n$ et
$\alpha_{n-1}q_n>1/2$,
$$\inv{\lambda(J_n)}\left|\int_{J_n}(1-\theta_n)d\lambda\right| \leq 
2 \|1-\theta_n\|_1+ 2(\pi|sb_n+b'_n|+4|b_n|) \alpha_n q_n .$$ 
Par hypoth\`ese le second terme converge vers $0$ et on en
d\'eduit (ii).
\\
Pour montrer que (ii) entraine (iii), il suffit de 
calculer explicitement l'int\'egrale pour  une suite $(J_n)$ bien
choisie. 
Supposons par exemple que $b_n>0$. On choisit la suite des bases des
tours majeures de la suite $(\C D_n(\beta))$, qu'on note $(B_n)$.  
Comme les discontinuit\'es de $\theta_n$ sont toutes dans $\C I_n$,
celle-ci est continue sur le sous-intervalle de $B_n$,   
$B'_n=T^{-q_n}(]\beta_{n+1},T^{-q_{n-1}}\beta_n[)$ (voir figure
\ref{3tours}).
On peut donc \'ecrire pour tout $x\in B'_n$ 
$$\theta_n(x)=c_n \exp{(2i\pi (sb_n+b'_n)q_nx)}\; ,$$
o\`u $c_n$ est une constante de module 1.
Comme $\lambda(B'_n)/\lambda(B_n) \rightarrow 1$, (ii) entra\^ \i ne
$$ \int_{B_n}{(1-c_n\e{2i\pi(sb_n+b'_n)q_nx})\over \alpha_{n-1}}dx
\tend{n \rightarrow \infty}0.$$
Le calcul  de cette int\'egrale montre alors que 
$|sb_n+b'_n|\tend{n \rightarrow \infty}0$ et que 
$(c_n)$ converge vers 1. On obtient donc (iii).\\
Supposons maintenant (iii)~: dans $\C D_n(\beta)$, $\C I_n$ est
\'egale \`a la r\'eunion des \'etages de la tour majeure priv\'ee de
sa partie principale (cf figure~\ref{3tours})~: par cons\'equent il
suffit de v\'erifier que $\theta_n$ converge uniform\'ement sur la
partie principale et la tour mineure d'ordre $n$, de base $B_n^0$.
Gr\^ ace \`a  la relation (\ref{invariant}), il suffit en fait de
v\'erifier la convergence uniforme sur  l'intervalle $B'_n \cup
B_n^{0}$. Comme par construction $\theta_n$ est continue sur cet
intervalle de longueur inf\'erieure \`a $\alpha_{n-1}$, affine et de
pente constante \'egale \`a $(sb_n+b'_n)q_n$ on  d\'eduit de (iii)
que $\|(1-\theta_n)\unde{B'_n\cup B_n^{0}}\|_\infty$ converge vers
0.\\
Enfin si on suppose (iv), comme la mesure de $\C I_n$ converge vers
$0$, on en d\'eduit imm\'ediatement (i). 
\findem

\rem 
Comme $sb_n+b'_n \rightarrow 0$, si la suite $(b_n)$ est born\'ee
alors $s$
est n\'ecessairement rationnel. Dans ce cas 
pour tout $n$ assez grand, $b_ns \in \M Z$ et $-b'_n=b_ns$. 

\bigskip

Il reste pour terminer \`a \'evaluer la vitesse de convergence de
$(\int_\M T \theta_nd\lambda)$ vers 1. Comme d\'eja dit dans le
paragraphe \ref{cnidee}, l'id\'ee est d'approcher $\int_\M T
f_n\theta_nd\lambda$ par le produit des int\'egrales, puis d'en
d\'eduire $\prod|\int_\M T \theta_nd\lambda|>0$.\\
Notons $B_n$  la  base de la  tour majeure associ\'ee \`a $\C
D_n(\beta)$  (cf paragraphe \ref{df}), et  $\C F_n$ la tribu
engendr\'ee par les \'etages de $\C D_n(\beta)$. Par construction,
$(\C F_n)$ est une filtration qui engendre la tribu  bor\'elienne. On
a d'abord les estimations suivantes~:

\begin{prop}\label{ia} 
On suppose que $(b_nq_n\alpha_n)$ et $(b'_nq_n\alpha_n)$ convergent
vers $0$. On note $(f_n)$ et $(\theta_n)$ les fonctions d\'efinies au
paragraphe \ref{transfert}. Si de plus $(\theta_n)$ converge vers 1
dans $L^2$, alors on peut trouver une suite $(\varepsilon_n)$ qui
converge vers $0$ telle que pour tout $n\in \Lambda$ et pour toute
fonction $g$ $\C F_n$-mesurable de module inf\'erieur \`a 1 on ait~:
$$\int_{\M T}gf_{n+1}d\lambda =
\int_{\M
T}gf_nd\lambda\left(\inv{\lambda(B_n)}\int_{B_n}\theta_nd\lambda\right)+
\Oun(V_{B_n}(f_n)+\varepsilon_n \alpha_nq_n). $$
 En notant $\{n_j, j\in \M N\}$  l'ensemble des \'el\'ements croissants de $\Lambda$ on a aussi pour tout $j$
$$V_{B_{n_j}}(f_{n_j})= o(b_{n_{j-1}}\alpha_{n_{j-1}}q_{n_{j-1}}).$$  
\end{prop} 

\dem
On fixe  $n \in \Lambda $. 
On note $\C T_n$ (respectivement $\C T_n^0$) la  r\'eunion des
\'etages de la  tour majeure (mineure) associ\'ee \`a $\C
D_n(\beta)$. 
Quitte \`a supprimer le premier ou le dernier \'etage de  $\C T_n$, 
on a $\psi_n=1$ sur $\C T_n$  et comme $g$ est constante sur les
\'etages de la tour majeure, on peut  appliquer le lemme~\ref{l1}
(en remarquant que la variation de $f_n$ est la m\^eme sur tous les
\'etages)~:
$$
\int_{\C T_n}gf_{n+1}d\lambda =  \inv{\lambda(B_n)}
\int_{\C T_n}gf_nd\lambda\int_{B_n}\theta_nd\lambda 
+\lambda(\C T_n)\Oun(V_{B_n}(f_n)+{2\pi
|sb_n+b'_n|\alpha_n\over\alpha_{n-1}}). $$
Sur la tour mineure,  on a plus simplement~: 
\begin{flalign*}
\int_{\C T_n^{0}}gf_{n+1}d\lambda
&=  
\int_{\C T_n^{0}}gf_nd\lambda
\inv{\lambda(B_n)}\int_{B_n}\theta_nd\lambda 
+\int_{\C T_n^{0}}gf_n\left(\theta_n-
\inv{\lambda(B_n)}\int_{B_n}\theta_nd\lambda \right)d\lambda\\
&=  
\int_{\C T_n^{0}}gf_nd\lambda
\inv{\lambda(B_n)}\int_{B_n}\theta_nd\lambda \\
& \qquad\quad +\; \lambda(\C T_n^{0})\Oun\left(\|(1-\theta_n)\unde{\C
T_n^{0}}\|_\infty +|1- \inv{\lambda(B_n)}\int_{B_n}\theta_nd\lambda
|\right). 
\end{flalign*}
Posons maintenant $\varepsilon_n=
4\pi|sb_n+b'_n|+\|(1-\theta_n)\unde{\C T_n^{0}}\|_\infty +|1-
\inv{\lambda(B_n)}\int_{B_n}\theta_nd\lambda |$~: cette suite
converge vers $0$ d'apr\`es le lemme~\ref{bs-b'} et 
on obtient  donc la premi\`ere relation de l'\'enonc\'e en remarquant
que $\lambda(\C T_n^0)=\alpha_n q_{n-1}\leq \alpha_nq_n$ et que
$1/\alpha_{n-1}\leq 2q_n$.\\
En ce qui concerne l'estimation de $V_{B_n}(f_n)$, comme 
$f_n$ est   continue sur $B_n$  on peut \'ecrire 
$V_{B_n}(f_n)\leq 2\pi |sk_n+k'_n|\alpha_{n-1}$.
Sachant que lorsque $n$ est assez grand,  $b_n=0$ entra\^ \i ne
$b'_n=0$, on 
obtient pour tout $j$ assez grand
$$sk_{n_j}+k'_{n_j}= \sum_{l=1}^{j-1}(sb_{n_l}+b'_{n_l})q_{n_l}=
sk_{n_{j-1}}+k'_{n_{j-1}}+(sb_{n_{j-1}}+b'_{n_{j-1}})q_{n_{j-1}}.$$
Comme d'apr\`es les hypoth\`eses, $sk_n+k'_n = o(q_n)$ et $b_ns+b'_n
\rightarrow 0$, il vient $sk_{n_j}+k_{n_j}= o(q_{n_{j-1}})$.
En remarquant que  $\alpha_{n_j-1}\leq
|b_{n_{j-1}}|\alpha_{n_{j-1}}$, on trouve 
bien l'estimation annonc\'ee.   
\findem

Explicitons maintenant la proposition~\ref{ia} appliqu\'ee \`a une
fonction $g$ choisie convenablement~:
comme $(f_n)$ converge vers $f$ qui n'est pas nulle, on peut trouver
$N$ et $g$ une fonction $\C F_N$-mesurable de module inf\'erieur \`a
$1$ telle que $\int_\M T f gd\lambda\neq 0$ (et telle que pour tout
$n\geq N$, $\int f_ngd\lambda\neq 0$).
Alors pour tout $j$ tel que  $n_j \geq N$, $g$ est encore $\C
F_{n_j}$-mesurable et la proposition~\ref{ia} donne~:
\begin{multline*}
\int_{\M T}gf_{n_{j+1}}d\lambda=\\
 \qquad\quad \int_{\M T}gf_{n_j}d\lambda
\left(\inv{\lambda(B_{n_j})}\int_{B_{n_j}}\theta_{n_j}d\lambda
+  o(b_{n_{j-1}}\alpha_{n_{j-1}}q_{n_{j-1}})
+o(|b_{n_j}|\alpha_{n_j}q_{n_j})\right).
\end{multline*}
Il en r\'esulte que le produit correspondant
$$\prod_{j>0}\left|\inv{\lambda(B_{n_j})}\int_{B_{n_j}}\theta_{n_j}d\lambda
+
o(b_{n_{j-1}}\alpha_{n_{j-1}}q_{n_{j-1}})+o(|b_{n_j}|\alpha_{n_j}q_{n_j})\right|
$$
est strictement positif ce qui entra\^ \i ne alors la convergence de
la s\'erie  de terme g\'en\'eral
\begin{equation}\label{cni}
1-\left|\inv{\lambda(B_{n_j})}\int_{B_{n_j}}\theta_{n_j}d\lambda+ 
o(b_{n_{j-1}}\alpha_{n_{j-1}}q_{n_{j-1}})+o(|b_{n_j}|\alpha_{n_j}q_{n_j})\right|.
\end{equation}
Il reste \`a calculer
l'int\'egrale de $\theta_n$ sur $B_n$, pour $n \in \Lambda$. 
On a le r\'esultat suivant~:
\begin{lemme}\label{calcul} 
On reprend les notations pr\'ec\'edentes. On suppose que
$(b_nq_n\alpha_n)$,  $(b'_nq_n\alpha_n)$ et $(b'_n+b_ns)$ convergent
vers $0$~; alors on a 
$$\frac{1}{\lambda(B_n)}\left|\int_{B_n}\theta_nd\lambda\right|=
1-(b_n \alpha_nq_n+\frac{\pi^2}{6}(sb_n+b'_n)^2)(1+o(1)).$$
\end{lemme}
\dem
Il suffit de d\'emontrer l'\'egalit\'e pour $n \in \Lambda$. 
On suppose, pour all\'eger les notations que $n$ est pair , $k_n\geq
0 $ et 
$b_n>0$. $B_n$ \'etant un intervalle de longueur $\alpha_{n-1}$ il
s'\'ecrit
$B_n=[x_n,x_n+\alpha_{n-1}[$. Sur $B_n$, on sait par construction que
$\theta_n$ est affine par morceaux  de pente constante \'egale \`a
$q_n(sb_n+b'_n)$ et que, d'apr\`es le choix de $\C D_n(\beta)$, ses
discontinuit\'es se trouvent aux points de la forme $x_n+j\alpha_n$
pour $j \in \{1,..b_n\}$ (cf figure~\ref{3tours}). Comme les sauts
aux discontinuit\'es sont constants et \'egaux \`a $-s$, on peut donc
\'ecrire sur $[x_n,x_n+b_n\alpha_n[$ ($=I_n$)
\begin{eqnarray*}
\int_{I_n}\theta_nd\lambda 
& =& 
\sum_{j=0}^{b_n-1}\int_{x_n+j\alpha_n}^{x_n+(j+1)\alpha_n}\theta_nd\lambda\\
&=&
\sum_{j=0}^{b_n-1}\e{i\pi j(-s+q_n(b_ns+b'_n)\alpha_n)}
\int_{x_n}^{x_n+\alpha_n}\theta_nd\lambda\\
&=& 
\e{i\pi b_n(-s+q_n\alpha_n(b_ns+b'_n))}
\frac{\sin{(\pi b_n(-s+q_n\alpha_n(b_ns+b'_n)}}
{\sin{\pi(-s+q_n\alpha_n(b_ns+b'_n))}}
\int_{x_n}^{x_n+\alpha_n}\theta_nd\lambda
\end{eqnarray*}  
Or, on sait que, $(b_nq_n\alpha_n(b_ns+b'_n))_n$ converge vers 0, que
$s \notin \M Z$ et que $(\|b_ns\|)_n$ converge vers 0. 
Par cons\'equent on obtient 
$$\left| \int_{I_n} \theta_nd\lambda\right|=
o(\int_{x_n}^{x_n+\alpha_n} \theta_nd\lambda)= o(\alpha_n).$$
Sur $B'_n$, $\theta_n$ est continue, affine et de pente
$q_n(sb_n+b'_n)$, d'o\`u, en notant $\alpha'_{n-1}=
\alpha_{n-1}-b_n\alpha_n$,
\begin{eqnarray*}
\left|\int_{B'_n}\theta_nd\lambda\right|
& = &
\int_{-\alpha'_{n-1}/2}^{\alpha'_{n-1}/2} \e{2i\pi q_n(sb_n+b'_n)x}dx 
\\
& = &
\frac{\sin{(\pi (sb_n+b'_n)q_n\alpha'_{n-1})}}
{\pi(sb_n+b'_n)q_n\alpha'_{n-1}} \\
& = & \alpha'_{n-1}
\left(1-{\pi^2\over
6}(b_ns+b'_n)^2(\alpha'_{n-1}q_n)^2(1+o(1))\right).
\end{eqnarray*}
Comme sur $\Lambda$,
$\int_{I_n}\theta_nd\lambda=o(\int_{B'_n}\theta_nd\lambda)$ et que la
suite  $(q_n\alpha_{n-1})$ converge vers 1, 
on obtient finalement le r\'esultat annonc\'e.
\findem

Il r\'esulte de ce lemme et de la convergence de la s\'erie donn\'ee
par  (\ref{cni}) qu'une  s\'erie de terme g\'en\'eral du type
$${\pi^2\over
6}(b_{n_j}s+b'_{n_j})^2(1+o(1))+{|b_{n_j}|q_{n_j}\alpha_{n_j}}(1+o(1))
+o({|b_{n_{j-1}}|q_{n_{j-1}}\alpha_{n_{j-1}}})$$
est convergente. 
Nous en d\'eduisons finalement les conditions du th\'eor\`eme
\ref{marche}~:
\centerline{
$\begin{array}{lll}
\DSM \sum_n \frac{|b_n|}{a_{n+1}}&<& +\infty\\
\DSM \sum_n (b_ns+b'_n)^2 &< & +\infty.
\end{array}$}\\
 Ceci cl\^ot la recherche des conditions n\'ecessaires \`a
l'existence de solutions \`a l'\'equation de cobord.

\subsection{Des conditions suffisantes.}\label{cs}

On suppose maintenant les conditions du th\'eor\`eme~\ref{marche}~:
il existe une suite d'entiers $(b_n)_{n\geq 0}$ telle que
$$
\beta = \sum_0^\infty b_n q_n \alpha  \qquad \hbox{avec } \qquad
\sum_n\frac{|b_n|}{a_{n+1}}<\infty \quad \hbox{et} \quad \sum_n
\|b_ns\|^2<\infty.$$
Posons $t= \sum_n b'_nq_n\alpha$ avec $b'_n=-[b_ns]$~: il suffit de
montrer qu'on peut construire une fonction $f$ telle que l'\'equation
(\ref{eqbis}) soit v\'erifi\'ee.  

On reprend les notations des paragraphes pr\'ec\'edents.
On choisit la constante des fonctions $\theta_n$ du paragraphe
\ref{transfert} de sorte que $ 
\theta_n(m_n)=1$, o\`u $m_n$ d\'esigne le milieu de $B'_n$.
Pour trouver une fonction $f$ v\'erifiant (\ref{eqbis}) il suffit  de
montrer que la suite $(f_n)$ d\'efinie par $f_n =
\prod_0^{n-1}\theta_k$ admet une valeur d'adh\'erence pour la topologie faible $L^2$ qui soit non identiquement nulle.
Pour cela, on utilise la proposition~\ref{ia}~: par construction de
$(\theta_n)$ les hypoth\`eses de la proposition sont bien
v\'erifi\'ees. 
L'estimation de $V_{B_n}(f_n)$ donn\'ee par cette proposition et la
sommabilit\'e de la s\'erie $(|b_n|q_n\alpha_n)$ permettent
d'affirmer qu'il existe une suite sommable $(u_n)$ v\'erifiant~: pour
tout $n$ et pour toute $g$ $\C F_n$-mesurable de module inf\'erieur
\`a $1$, on a 
$$\int_\M T f_{n+1}g d\lambda = \int_\M T f_n g d\lambda
\left(\inv{\lambda(B_n)}\int_{B_n}\theta_nd\lambda\right)  + u_n.$$ 
D'autre part les conditions arithm\'etiques sur $(b_n)$ et le choix
de $(\theta_n)$ donnent gr\^ace au lemme~\ref{calcul}~:
$$\sum_{n>0}\left|1-\inv{\lambda(B_n)}\int_{B_n}\theta_nd\lambda\right|<+\infty.$$
Par cons\'equent le produit
$\prod_{n>0}(\inv{\lambda(B_n)}\int_{B_n}\theta_nd\lambda)$ converge
vers une 
limite non nulle. 
On d\'eduit de ces  remarques  que toute valeur d'adh\'erence faible
de $(f_n)$, $f$, satisfait~: pour tout $n$ et pour toute fonction $g$
$\C F_n$-mesurable de module inf\'erieur \`a 1, 
$$\int_{\M T}fgd\lambda= \int_\M T f_ng d\lambda
\prod_{k=n}^{\infty}\inv{
\lambda(B_k)}\int_{B_k}\theta_kd\lambda+\Oun(\sum_{k=n}^\infty
u_k).$$
Comme 
 $\prod_n^\infty \inv{\lambda(B_k)}\int_{B_k}\theta_kd\lambda \tend{n
\rightarrow \infty}1$ et $\sum_{k=n}^{+\infty}u_k \tend{n \rightarrow
\infty}0$, pour tout $n$ assez grand on a $\sum_{k \geq n} u_k<1/4$ et
$\prod_{k\geq n}\inv{\lambda(B_n)}|\int_{B_k}\theta_kd\lambda|>1/2$. 
Il reste maintenant \`a trouver $n$  et $g$ de sorte que $|\int_{\M
T}gf_nd\lambda|$ soit strictement sup\'erieur \`a 1/2~: on sait par
construction que 
$f_n$ est affine sur chacun des \'etages de $\C D_n(\beta)$ (ses
discontinuit\'es sont sur le bord des \'etages), et de pente
constante \'egale \`a $(sk_n+k'_n)$. Par cons\'equent sa variation
moyenne sur chacun des \'etages de $\C D_n(\beta)$ est de l'ordre de
$\Oun(\alpha_{n-1}(sk_n+k'_n))$ qui converge vers $0$ lorsque $n $
tend vers l'infini. Choisissons $\overline g=\M E[f_n|\C F_n]$ la
projection orthogonale de $f_n$ sur les \'etages de $\C D_n(\beta)$,
on obtient alors
\begin{eqnarray*}
\int_{\M T}f_ngd\lambda &= & \|f_n\|_2^2-\|f_n-\overline g\|_2^2\\
& =& 1-\sum_{j=0}^{q_n-1} \int_{T^jB_n}|f_n-\overline
g|^2d\lambda+\Oun(\alpha_nq_{n-1})\\
 &=
&1-\Oun(q_n\alpha_{n-1}(sk_n+k'_n)\alpha_{n-1}+q_{n-1}\alpha_n)\\
&\tend{n \rightarrow +\infty}& 1.
\end{eqnarray*}
Pour $n$ assez grand on peut donc trouver $g$ $\C F_n$-mesurable
telle que $|g|\leq 1$ et  
$|\int_\M T f_ngd\lambda|>1/2$. 
Ceci montre donc 
que $|\int_\M T fgd\lambda |>0$. Finalement  
$f$  est non nulle et la r\'eciproque est bien compl\`ete.

\noindent{\bf Remarque~:}
En raison des r\'esultats du paragraphe pr\'ec\'edent, on obtient en
fait que la suite $(f_n)$ converge fortement dans $L^1$ vers la
fonction de transfert $f$.

\section{Cas g\'en\'eral des fonctions en escalier.}

L'objet de cette partie est la preuve du th\'eor\`eme~\ref{cn}.
On se donne, pour un entier $m$ non nul $m+1$ \'el\'ements de $\M T$
distincts  $(\beta_j)_{j\geq 0}$ tels que $\beta_0=0$, et $m+1$
r\'eels non entiers de somme nulle $(s_j)_{j\geq 0}$. On note ensuite 
$$\phi= \sum_{j=0}^{m}s_j\phi_{\beta_j}, \quad \hbox{et } \quad
\varphi= \exp{(2i\pi\phi)}.$$ 
Le plan suivi est analogue \`a celui de la partie pr\'ec\'edente~: Le
premier paragraphe est consacr\'e \`a la n\'ecessit\'e des conditions
du th\'eor\`eme~\ref{cn}, sous l'hypoth\`ese du th\'eor\`eme sur
l'unicit\'e de l'\'ecriture de $\phi$. 
le second paragraphe \'etablit la preuve de la r\'eciproque en toute
g\'en\'eralit\'e. 

\subsection{Conditions n\'ecessaires.}

On travaille ici avec l'hypoth\`ese suivante~:

\bigskip

 {\it  Pour toute partie stricte non vide  $ J$ de $\{0,..,m\}$, on a 
$ \sum_Js_j \notin \M Z $.}\hfill (H)

\bigskip

\noindent
Comme pr\'ec\'edemment,  on suppose 
qu'il existe $t\in \M T$ et $F$ mesurable de $\M T$ dans $\M S^1$
v\'erifiant 
\begin{equation} \label{mcob}
\varphi= \e{2i\pi t} \frac{F}{F\rond T},\end{equation}
et on note $\tilde F$ un relev\'e de $F$ sur $\M R$.

\subsubsection{D\'ecomposition d'Ostrowski.}
Nous avons toujours d'apr\`es (\ref{mcob}) la condition  spectrale 
$|\int_\M T \varphi^{(q_n)}d\lambda|\rightarrow
1$. Avec l'hypoth\`ese (H), on \'etablit la proposition suivante~:

\begin{prop}\label{init}
Sous l'hypoth\`ese (H), si  
$|\int_\M T \varphi^{(q_n)}d\lambda|\tend{n \rightarrow\infty} 1$,
alors on a    
pour tout $j \in \{0,..,m\}$
$$\|\beta_jq_n\|\tend{n \rightarrow \infty} 0.$$
\end{prop}

\dem
Supposons par l'absurde qu'il existe $j \in \{1,..m\}$ tel que 
$\|\beta_jq_n\|$ ne converge pas vers $0$. On pose par exemple $j=1$. 
Alors l'\'equivalence (\ref{vbeta}) (voir l'appendice) montre qu'on
peut trouver  $\varepsilon>0$ et un ensemble infini d'entiers,$\Lambda_1$, tels que $\beta_1 \notin \cup_{n\in\Lambda_{1}}\C
V_n(\varepsilon)$,  o\`u on a pos\'e
$h_{n}=\min(\varepsilon/\alpha_{n}, q_{n}/2)$ et 
$$\C V_n(\varepsilon)=\{x_{n}+k_{n}\alpha, |x_{n}|\leq
\varepsilon/q_{n}, \;  \hbox{et }  |k_{n}|\leq h_{n}\}.$$
Remarquons que lorsque $\varepsilon\leq 1/4$, l'\'ecriture dans $\C
V_{n}(\varepsilon)$ est unique et qu'on peut donc le repr\'esenter
dans un domaine fondamental d'ordre $n$ par un rectangle (voir figure
\ref{vois0}). De plus,  ces rectangles sont croissants en fonction de
$\varepsilon$. On   r\'epartit maintenant les $(\beta_i)$ en deux
cat\'egories, selon que $\beta_i$ est dans $\limi_{n \in \Lambda_{1}}
\C V_n(\varepsilon)$ pour
tout $\varepsilon$, ou non. Par extractions successives, on peut donc
trouver un ensemble infini d'entiers, $\Lambda$, et $\varepsilon_0\in
]0,1/8[$ tels que 
pour tout $i \in \{0,..,m\}$ on ait \\
(i) 
ou bien $\beta_i \in \limi_{n\in \Lambda} \C V_n(\varepsilon)$ 
pour tout $\varepsilon>0$,\\
(ii)
ou bien 
$\beta_i \notin \cup_{n \in \Lambda}\C V_n(2\varepsilon_0)$.\\
On notera $I= \{i \in \{0,..,m\}, \; \beta_i \hbox{ satisfait
(i)}\}$. Remarquons que $I$ est une partie non triviale de
$\{0,..,m\}$ (car $0 \in I$ mais $1 \notin I$).
Enfin, on peut encore supposer, quitte \`a restreindre $\Lambda$, que 
$\lim_{n \in \Lambda}a_{n+1}=+\infty$, ou que $\lims_{n \in
\Lambda}a_{n+1}<\infty$.

{\noindent {\bf Premier cas~: $\lim_{n \in \Lambda} a_{n+1}=
\infty$.}}\\
Choisissons $\varepsilon=\varepsilon_{0}$, de sorte que pour tout $i
\notin I$, 
$\beta_i \notin \cup_{\Lambda}\C V_n(2\varepsilon)$.
Comme dans ce cas la suite $(\alpha_nq_n)_n$ converge vers $0$, on
peut aussi supposer que pour tout $n \in \Lambda$, on a $\alpha_n
q_n<\varepsilon/2$. 
Alors $\C V_n(\varepsilon)$ et $\C V_n(2\varepsilon)$ sont deux tours
de
hauteur $q_n$, qu'on peut repr\'esenter dans un domaine fondamental
d'ordre
$n$ (voir figure~\ref{QPNB}).
\begin{figure}
\begin{center}
\caption{$\lim_{n \in \Lambda}a_{n+1}=+\infty$}\label{QPNB}
\input{QPNB.pstex_t}
\end{center}
\end{figure} 
Dans $\C V_n(2\varepsilon)\backslash \C V_n(\varepsilon)$ il n'y a
aucun des $\beta_i$. De plus les discontinuit\'es de $\phi^{(q_n)}$
issues
d'un $\beta_i$ sont contenues dans une tour de largeur $\alpha_n$
limit\'ee \`a
droite par $\beta_i$ (cas $n$ pair, cf figure~\ref{QPNB}). 
On en d\'eduit qu'il existe deux tours d'intervalles, dont on note
$\C T_n'$ et $\C T_n''$ la r\'eunion des \'etages,  
 de hauteur $q_n$ incluses dans $\C V_n(2\varepsilon)$, situ\'ees de
part et d'autre de $\C V_n(\varepsilon)$, de largeur
$\varepsilon/q_n-\alpha_n$, et 
dans lesquelles $\phi^{(q_n)}$ ne poss\`ede aucune discontinuit\'e. 
Par cons\'equent $\phi^{(q_n)}$ est  constante sur chacune de ces tours  
(cf paragraphe \ref{betaqn}).
Enfin, pour passer d'un \'etage de $\C T_n'$ \`a un \'etage de $\C
T_n''$, 
il faut franchir toutes les discontinuit\'es issues des points dans 
$\C V_n(\varepsilon)$, c'est \`a dire que le saut de $\phi^{(q_n)}$
entre $\C 
T_n'$ et $\C T_n''$ est \'egal \`a la somme des sauts $-s_i$ pour $i
\in I$~: 
il en d\'ecoule que 
$$|\int_{\C T_n'}\varphi^{(q_n)}d\lambda- \int_{\C
T_n''}\varphi^{(q_n)}d\lambda|=
\varepsilon |1-\e{-2i\pi\sum_{k\in I}s_k}|.$$
Comme $\alpha_n<\varepsilon/2q_n$, on a 
$\lambda(\C T_n')= \lambda(\C T_n'')>\varepsilon/2$. 
Or, on doit avoir $|\int_\M T \varphi^{(q_n)}d\lambda|\rightarrow 1$,
donc 
n\'ecessairement $\sum_{k \in I}s_k=0\mod 1$. Avec l'hypoth\`ese (H),
ceci implique que $I$ doit \^etre triviale, ce qui est absurde.

\medskip

\begin{figure}
\caption{Lignes de sauts dans $\C V_n(2\varepsilon)$.}\label{QPB2}
\input{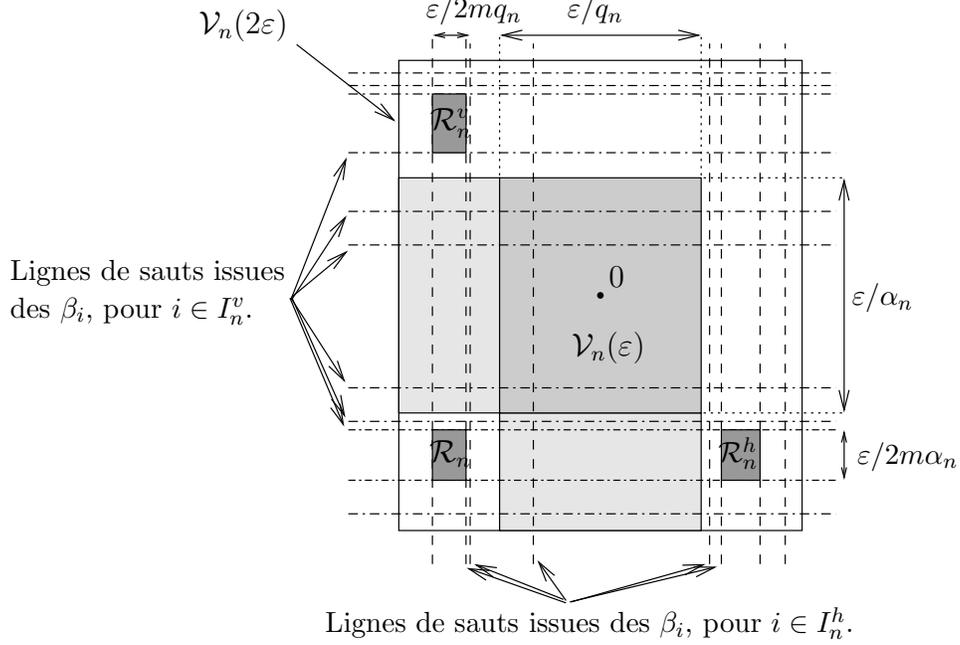}
\end{figure}

{\noindent {\bf Second cas~: $\lims_{n \in
\Lambda}a_{n+1}<\infty$.}}\\
Comme on a $\limi_\Lambda \alpha_n q_n>0$, on choisit
$\varepsilon<\varepsilon_0$ tel que $4\varepsilon/q_n<\alpha_n$ pour
tout $n \in \Lambda$. Alors  pour tout $i \notin I$, $\beta_i \notin
\cup_{\Lambda}\C V_n(2\varepsilon)$. On s'int\'eresse \`a la
r\'epartition des valeurs de $\phi^{(q_n)}$ dans $\C
V_n(2\varepsilon)$. Rappelons que dans un domaine fondamental d'ordre
$n$, $\phi^{(q_n)}$  est constante sur les intervalles des \'etages
qui ne contiennent pas les points $(T^{-j}\beta_i)_{0 \leq j <q_n}$,
et constante sur leur images, s'ils sont disjoints de
$[T^{-q_n}\beta_i, \beta_i[$ (cf paragraphe \ref{betaqn}). Par
cons\'equent, dans n'importe quel domaine fondamental d'ordre $n$,
$\phi^{(q_n)}$ est constante dans toutes les parties connexes
d\'elimit\'ees par  les lignes bris\'ees issues des $\beta_i$ (qui
indiquent les instants de saut $-s_i$, si $n$ pair, pour
$\phi^{(q_n)}$)~: on les appellera lignes de sauts issues de
$\beta_i$ (voir la figure~\ref{QPB1}). 
\begin{figure}
\begin{center}
\caption{ $\lims_{n \in \Lambda}a_{n+1}<\infty$}\label{QPB1}
\input{QPB1.pstex_t}
\end{center}
\end{figure} \\
-- Lorsque $i \in I$, les lignes de sauts issues de $\beta_i$ sont
incluses
dans une bande bris\'ee coupant $\C V_n(2\varepsilon)$ en deux
parties (repr\'esent\'ees par les parties claires de la figure
\ref{QPB2}).\\
--  Si $i \notin I$ et que $T^{q_n}\beta_i \notin \C
V_n(2\varepsilon)$,  
la ligne de sauts issue de $\beta_i$ traverse  \'eventuellement
$\C V_n(2\varepsilon)$  selon une droite horizontale ou verticale 
(le ``ou'' est exclusif ici car sinon, on aurait $\beta_i \in \C
V_n(2\varepsilon)$). \\
--  Si $i \notin I$ et que $T^{q_n} \beta_i \in \C
V_n(2\varepsilon)$, la ligne de saut issue de $\beta_i$ ne coupe que
le bord haut et le bord droit de $\C V_n(2\varepsilon)$ (cas $n$
pair).\\
Comme il y a au plus $m$ lignes de sauts, on peut 
trouver dans $\C V_n(2\varepsilon)\backslash\C V_n(\varepsilon)$ 2
tours l'une \`a gauche et l'autre \`a droite de $\C
V_n(\varepsilon)$, de 
hauteurs $2\varepsilon/\alpha_n$ et de largeurs $\varepsilon/2mq_n$
et 
qui ne contiennent aucune partie verticale des lignes de sauts 
(c'est \`a dire aucun des $T^j\beta_i$ pour $0<j\leq q_n$, $0 \leq i
\leq m$). 
De la m\^eme mani\`ere on trouve 2 tours, l'une en dessous et l'autre
au dessus de $\C V_n(\varepsilon)$, de largeurs $2\varepsilon/q_n$ et
de hauteurs $\varepsilon/2m \alpha_n$   dont les \'etages  sont
disjoints des lignes de sauts (c'est \`a dire disjoints des
intervalles $[T^{q_n}\beta_i, \beta_i[$). Les intersections de ces
tours dans $\C V_n(2\varepsilon)$ forment 4 tours qui sont disjointes
de toutes les lignes de sauts~: 
les 3 qui nous int\'eressent sont repr\'esent\'ees par $\C R_n$, $\C
R_n^h$ et $\C R_n^v$ sur la figure~\ref{QPB2}. 
Sur chacune de ces tours, $\phi^{(q_n)}$ est donc constante, 
et comme la mesure de chacune de ces tours  est sup\'erieure \`a
$\varepsilon^2/4m^2$, la
condition $|\int_\M T\varphi^{(q_n)}d\lambda|\rightarrow 1$ impose
donc que
les valeurs de $\phi^{(q_n)}$ sur chacune d'elles soient \'egales
 modulo 1 pour tout 
$n$ assez grand (car il n'y a qu'un nombre fini de sauts).
Il suffit pour terminer de compter les sauts de $\phi^{(q_n)}$ pour
passer 
d'une tour \`a l'autre.
Pour aller de $\C R_n$ \`a $\C R_n^h$ ou \`a $\C R_n^v$, il faut
traverser toutes les lignes de sauts issues des $\beta_i$ pour $i \in
I$. 
On note ensuite $I_n^v$ (resp. $I_n^h$) l'ensemble des indices $i
\notin I$ des lignes de sauts qu'on doit traverser pour aller de $\C
R_n$ \`a $\C R_n^v$ (resp. $\C R_n^h$). 
On doit donc avoir pour tout $n$ assez grand 
$\sum_{I\cup I_n^h}s_i=\sum_{I \cup I_n^v}s_i=0\mod 1$. Alors comme
$I$ est non vide, on a n\'ecessairement que $I_n^h \cup I$ et $I_n^v
\cup I$ sont \'egaux \`a $\{0,..,m\}$. Mais $\beta_i \notin\C
V_n(2\varepsilon)$ si $i \notin I$, donc 
$I_n^v$ et $I_n^h$ sont disjoints. 
Par cons\'equent ces 2 ensembles sont n\'ecessairement vides, et $I$
contient donc $1$ ce qui est absurde.
\findem

L'existence d'une suite de fonctions de transfert convergeant
fortement dans 
$L^2$ vers $F$ se d\'emontre maintenant comme dans la partie
pr\'ec\'edente.
\\
D'apr\`es la proposition~\ref{init},  on peut appliquer le  lemme
\ref{phiqn}  pour tout $\beta_j$, donc  $(\phi_{\beta_j}^{(q_n)})$
converge vers 0 en mesure. On en d\'eduit que  $(\varphi^{(q_n)})$
tend vers $1$ dans $L^1$. 
De la condition (\ref{rigide}) il vient alors    
$\|tq_n\|\tend{n \rightarrow \infty} 0$. 
De l\`a on obtient les d\'ecompositions~:
$$ \forall j \in \{0,..,m\}, \quad
\beta_j =\sum_0^\infty b^{(j)}_n q_n \alpha
\qquad 
\hbox{et }\quad t =
\sum_1^\infty b'_n q_n\alpha ,$$
o\`u $(b^{(j)}_n)$ et $(b'_n)$ sont des suites d'entiers v\'erifiant
$\lim
b^{(j)}_n/a_{n+1}= \lim b'_n/a_{n+1} =0$. Nous utiliserons pour la
suite les 
notations suivantes~:
\\
\centerline{
$
\begin{array}{rclcl}
\beta_n^{(j)} &=& k_n^{(j)}\alpha   
& \hbox{avec} & \DSM 
k_n^{(j)}= \sum_0^{n-1}b_l^{(j)}q_l,
% \qquad  \<\beta_j-\beta_n^{(j)}\> = \sum_n^\infty b_l^{(j)}\<
%q_l\alpha\> , 
%\hbox{ et} 
\\
 t_n  &=& k'_n\alpha  & \hbox{avec} & 
\DSM  k'_n= \sum_0^{n-1} b'_jq_j.
%,  \qquad \<t-t_n\>= \sum_n^\infty b'_j\< q_j\alpha\> 
\end{array}
$}
\\
On remarquera que $k_n^{(j)}$ et $k'_n$ sont infiniment petits devant
$q_n$, et
que $\|\beta_j-\beta_n^{(j)}\|$ et $\|t-t_n\|$ sont infiniment petits
devant 
$\alpha_{n-1}$.

\subsubsection{Repr\'esentation. Choix des Domaines
Fondamentaux.}\label{+}
Nous avons besoin, comme dans la partie pr\'ec\'edente de
repr\'esenter les suites des fonctions de transfert \`a l'aide des
domaines fondamentaux d\'ecrits au paragraphe \ref{df}. Il suffira,
comme auparavant, d'\'etudier les termes de cette suite pour les
indices $n$ dans l'ensemble 
 $$\Lambda=\{n \in \M N, \; \max_{j}(|b_n^{(j)}|) >0\}.$$ 
Une premi\`ere simplification consiste \`a ne consid\'erer que les
cas o\`u tous les $b_n^{(j)}$ sont positifs~: pour cela il suffit de
translater $\phi$ par $\beta_{\min}=\sum_{n \in
\Lambda}\min(b_n^{(j)}-1)q_n\alpha$. La nouvelle fonction obtenue,
encore not\'ee $\phi$ est toujours un quasi-cobord, associ\'e \`a la
m\^eme valeur propre $t$, et \`a la fonction de transfert
translat\'ee. 
Elle s'\'ecrit encore sous la forme $\sum_0^m s_j\phi_{\beta_j}$,
o\`u les $(s_j)$ sont inchang\'es, et les $\beta_j$ sont les
translat\'es des anciens. On obtient donc bien
$$\hbox{Pour tout } j \in \{0,..,m\}, \quad \beta_j= \sum_n
b_n^{(j)}q_n\alpha  \quad \hbox{avec } \forall n, \; b_n^{(j)}\geq
0.$$
De plus on a la propri\'et\'e suivante~:
$$n \in \Lambda \Longleftrightarrow \left( \min_j(b_n^{(j)})=1 \hbox{
et } \max_j(b_n^{(j)})\geq 2\right).$$ 
Pour tout $n$ on d\'efinit  $b_n=\max_j(b_n^{(j)})$ et $k_n=
\sum_0^{n-1}b_iq_i$, puis on note $\beta=\sum_0^\infty b_nq_n\alpha$.
On a comme dans la partie pr\'ec\'edente  $\|\beta q_n\| \rightarrow
0$. On choisira pour la suite de repr\'esenter les fonctions de
transfert dans les domaines fondamentaux associ\'es \`a $\beta$,
$\C{D}_n(\beta)$ (cf paragraphe \ref{df}). Pour $n \in \Lambda$
fix\'e, assez grand pour que $b_n<a_{n+1}/2$, on localise facilement
les points $(\beta_j)_j$ et $(\beta_n^{(j)})$ selon la figure
\ref{nvtour}.

\begin{figure}
\begin{center}
\caption{$\C D_n(\beta)$}\label{nvtour}
\input{nvtour.pstex_t}
\end{center}
\end{figure}

\subsubsection{Convergence forte des fonctions de transfert.}
 
Posons  avec les notations du paragraphe \ref{cobord} et par analogie
avec 
\ref{transfert},  
$$\tilde F_n= \sum_{j=0}^m s_j \omega_{k_n^{(j)}} + k'_n\omega. $$
D'apr\`es le paragraphe \ref{cobord} on a l'\'egalit\'e
$$\sum_0^ms_j\phi_{\beta_n^{(j)}}= 
k'_n\alpha+\tilde F_n-\tilde F_n \rond T \mod 1.$$  
Pour tout $n$ on d\'efinit $F_n= \e{2i\pi \tilde F_n}$, on obtient la
proposition suivante~:
\begin{prop}
Sous l'hypoth\`ese (H), si $\varphi$ est un quasi-cobord, quitte \`a
multiplier $F_n$ par une constante de module 1, la suite de fonctions
$(F_n)_n$ converge fortement dans $L^2$ vers $F$.
\end{prop}

\dem
Comme $\varphi= \exp{(2i\pi t)}F/F\rond T$, on obtient avec ce qui
pr\'ec\`ede l'\'egalit\'e 
$$ (F \overline {F_n})\rond T = (F\overline{F_n})
\e{2i\pi\left((t-t_n)+\sum_0^m
s_j(\beta_j-\beta_n^{(j)})\right)}\exp{\left(2i\pi
\sum_0^ms_j\unde{[\beta_n^{(j)},\beta_j[}\right)}.$$
On veut maintenant appliquer le lemme~\ref{convergefort} pour la
suite 
de fonctions $(F \overline{F_n})$~: pour cela on choisit la suite des
tours majeures  de $\C D_n(\beta)$, $(T^jB_n)_{0 \leq j <q_n}$.  (cf
figure~\ref{nvtour}).
On v\'erifie ensuite les hypoth\`eses du lemme~\ref{convergefort}~: \\
La condition $(i)$ est bien v\'erifi\'ee pour $N=2$.\\
Les conditions $(ii)$ et $(iii)$ d\'ecoulent de la convergence des suites  
$(\|t-t_n\|q_n)$ et $(q_n\sum_{j=0}^m\|\beta_j-\beta_n^{(j)}\|)$ 
vers $0$.\\ 
Reste la condition $(iv)$~: $V_{B_n}(F \overline{F_n}) \rightarrow
0$. Pour la variation de $F_n$, lorsque $n \in \Lambda$, les
discontinuit\'es de $F_n$ sont contenues dans l'ensemble des  
points $(T^k 0)_{0 \leq k<k_n}$, qui se trouvent sur le bord de la
tour d'indice $n$. Par cons\'equent sur la base $B_n$, $F_n$ est 
affine et sa variation est major\'ee par  
$2\pi \alpha_{n-1} |\sum_{j<p}k_n^{(j)}+k'_n|$, qui converge bien
vers $0$ quand $n \rightarrow \infty$. 
Enfin, pour montrer la convergence vers $0$ de $V_{B_n}(F)$, il
suffit d'appliquer le lemme~\ref{proj} avec la suite des tours
$(T^jB_n)_{0 \leq j <q_n-k_n}$, en remarquant que comme $\varphi$ est
constante sur chacun des \'etages de ces tours, la variation de $F$
sur $T^jB_n$ 
est ind\'ependante de $j<q_n-k_n$ (cf proposition~\ref{l2}). Ceci
assure la condition $(iv)$ du lemme~\ref{convergefort}, et montre,
quitte \`a multiplier $F_n$ par une constante de module 1,  
la convergence de la suite $(F_n)$ vers $F$ dans $L^2$.
\findem

\subsubsection{R\'ecurrence.}

Les conditions arithm\'etiques du th\'eor\`eme~\ref{cn} s'obtiennent
selon le 
plan suivi dans la partie pr\'ec\'edente~:

Comme pr\'ec\'edemment, on note $\tilde\Theta_n= \tilde
F_{n+1}-\tilde F_n$. 
C'est encore une fonction de transfert, qui satisfait l'\'equation~:
$$\sum_{j=0}^ms_j\phi_{b_n^{(j)}q_n\alpha}\rond T^{-k_n^{(j)}}=
b'_nq_n\alpha+\tilde\Theta_n-\tilde \Theta_n\rond T\mod 1.$$
Par construction, $\tilde \Theta_n$ est affine par morceaux et de
pente constante \'egale \`a $(\sum_{j=0}^ms_jb_n^{(j)}+b'_n)q_n$.
L'ensemble de ses  discontinuit\'es est inclus dans  $\{T^k0, \;0
\leq k <k_{n+1}\}$. 
On note, comme dans  la partie pr\'ec\'edente, 
$$\Psi_n=\exp{-\sum_{j=0}^ms_j\unde{[0,b_n^{(j)}\<q_n\alpha\>[}\circ
T^{-k_n^{(j)}}}.$$ 
Alors la relation entre $\Theta_n$ et $\Psi_n$ s'\'ecrit~:
\begin{equation}\label{mtheta}
\frac{\Theta_n \rond T}{\Theta_n}= \Psi_n
\e{2i\pi(\sum_{j=0}^ms_jb_n^{(j)}+b'_n)\<q_n\alpha\>)}.
 \end{equation}

Comme $F_n$ converge dans $L^2$ vers $F$, $\Theta_n$ converge vers 1 
dans $L^2$, et en notant $\C{I}_n$ l'ensemble associ\'e \`a $\beta$
d\'efini par (\ref{in}), on obtient les propri\'et\'es suivantes,
semblables \`a celles du lemme~\ref{bs-b'}~:

\begin{lemme}\label{bsj}Avec les notations pr\'ec\'edentes, si pour
tout $j\in \{0,..,m\}$, $(b_n^{(j)}q_n\alpha_n)$ est une suite
positive convergeant vers $0$ et si $(b'_nq_n\alpha_n)$ converge vers
$0$, alors toutes les propri\'et\'es sont \'equivalentes~:
\begin{enumerate}
\item
$\|1-\Theta_n\|_2 \tend{n\rightarrow \infty} 0,$ 
\item
Pour toute suite d'intervalles $(J_n)_n$ de longueur
$\alpha_{n-1}$, on a 
$$\inv{\lambda(J_n)}\int_{J_n}(1-\Theta_n)d\lambda \tend{n \rightarrow
\infty}0,$$
\item
$\sum_{0\leq j \leq m}s_jb_n^{(j)}+b'_n  \tend{n \rightarrow \infty}
0$,  et $\exists x_n \in B'_n, \; \Theta_n(x_n)\tend{n \rightarrow
\infty}  1$,
\item
$ \|(1-\Theta_n)\unde{\C I_n^c}\|_\infty \tend{n \rightarrow \infty }
0$.
\end{enumerate}
\end{lemme}

\dem
On a $(iv) \Rightarrow (i)\Rightarrow (ii)$ pour les m\^emes raisons
que dans le lemme~\ref{bs-b'} ($\lambda(\C{I}_n^c)=b_n\alpha_nq_n
\rightarrow 0$ et $\lambda(\Psi_n\neq 1)\leq mb_n\alpha_n$).\\
Pour montrer $(ii) \Rightarrow (iii)$, on calcule encore
l'int\'egrale de $\Theta_n$ sur $B'_n$ (cf figure~\ref{nvtour2}).
Comme les discontinuit\'es de $\Theta_n$ sont dans $\C {I}_n$,
$\tilde \Theta_n$ est affine de pente constante $(\sum_0^m
s_jb_n^{(j)}+b'_n)q_n$ sur $B'_n$~: on obtient $(iii)$ par un calcul
identique au lemme~\ref{bs-b'}. \\
Pour montrer que $(iii)$ entra\^ \i ne $(iv)$, on v\'erifie d'abord
que  la convergence est uniforme sur les bases $B'_n$ et $B_n^{0}$
(base de la tour mineure)~: comme pour le lemme~\ref{bs-b'} c'est une
cons\'equence de la continuit\'e de $\Theta_n$ sur l'intervalle
$B'_n\cup B_n^{0}$ et de $(iii)$. Ensuite,  on obtient par
it\'eration de l'\'egalit\'e (\ref{mtheta}) que pour tout $0\leq
k<q_n$
$$\frac{\Theta_n \rond T^k}{\Theta_n}= 
\exp{\left(2i\pi\sum_{j=0}^ms_jb_n^{(j)}+b'_n)k\<q_n\alpha\>\right)}\Psi_n^{(k)}.$$
Lorsque $0 \leq k <q_n$, on sait par $(iii)$ que le terme constant
converge vers 1 uniform\'ement par rapport \`a $k$. De plus $\Psi_n$
vaut 1 sur 
 $\C I_n^c= \cup_{j<q_n}T^jB'_n \bigcup \cup_{j<q_{n-1}}T^jB_n^{0}$,
donc $\Psi_n^{(k)}=1$ sur $B'_n\cup B_n^{0}$ pour $0 \leq k <q_{n-1}$
et sur $B'_n$ si $0\leq k<q_n$. Par cons\'equent la convergence de
$\Theta_n$ vers 1 sur $T^kB'_n$ est bien uniforme pour $0 \leq k
<q_n$, de m\^eme que sur $T^kB_n^{0}$ si $0 \leq k <q_{n-1}$.  
 \findem

\medskip
On \'etablit comme dans la partie pr\'ec\'edente une estimation de la
corr\'elation entre $F_n$ et $\Theta_n$ sur $\C F_n$, la tribu
engendr\'ee par les \'etages de $\C D_n(\beta)$. On note $(B_n)$ la
suite des bases des tours majeures  associ\'ees \`a $\C D_n(\beta)$
(cf \ref{df}).
On obtient la proposition suivante~:

\begin{figure}
\begin{center}
\caption{Distribution de $\Theta_n$ dans $\C D_n(\beta)$ ($n$ pair).}
\label{nvtour2}
\input{nvtour2.pstex_t}
\end{center}
\end{figure}

\begin{prop}\label{ie2}
On reprend les notations pr\'ec\'edentes.
On suppose que pour tout $j$, $(b^{(j)}_nq_n\alpha_n)$ et
$(b'_nq_n\alpha_n)$ convergent vers $0$ (tous les $b_n^{(j)}$ sont
positifs). Si  $(\Theta_n)$ converge vers 1 dans $L^2$, alors on peut
trouver une suite $(\varepsilon_n)$ qui converge vers $0$ telle que
pour tout $n\in \Lambda$ et pour toute fonction $g$ $\C
F_n$-mesurable de module inf\'erieur \`a 1 on ait~:
$$\int_{\M T}gF_{n+1}d\lambda =
\int_{\M
T}gF_nd\lambda\left(\inv{\lambda(B_n)}\int_{B_n}\Theta_nd\lambda\right)+
\Oun(V_{B_n}(F_n)+\varepsilon_n (b_n\alpha_nq_n+k_n\alpha_{n-1})).
$$
 En notant $\Lambda=\{n_j, j\geq 1\}$ on a aussi pour tout $j\geq 1$
$$V_{B_{n_j}}(F_{n_j})=
o(b_{n_{j-1}}\alpha_{n_{j-1}}q_{n_{j-1}})\quad \hbox{et } \quad
k_{n_j}\alpha_{n_j-1}=\Oun(2b_{n_{j-1}}\alpha_{n_{j-1}}q_{n_{j-1}}).$$

\end{prop}

\dem
Pour $n \in \Lambda$, on d\'ecoupe  le domaine $\C D_n(\beta)$ en
trois tours~: la tour mineure d'ordre $n$, et deux tours de 
m\^eme largeur $\alpha_{n-1}$ et de hauteur $q_n-k_n$ et $k_n$, 
qui forment la tour majeure d'ordre $n$.  
On d\'efinit alors 
$\C T_n^{0}= \cup_{0\leq j <q_{n-1}}T^jB_n^{0}$, 
$\C T_n^1=\cup_{0 \leq j <q_n-k_n}T^jB_n$ et 
$\C T_n^2=\cup_{q_n-k_n\leq j <q_n}T^jB_n$ (voir figure
\ref{nvtour2}). 
Sur  
$\C T_n^1$, on sait que $\Psi_n=1$ et que $g$ est constante sur les
\'etages de la tour correspondante. L'\'egalit\'e (\ref{mcob}) permet
d'utiliser 
(\ref{tour}), ce qui donne~:
\begin{multline*}
\int_{\C T_n^1}F_{n+1}gd\lambda= 
\int_{\C T_n^1}F_ngd\lambda
\inv{\lambda(B_n)}\int_{B_n}\Theta_nd\lambda \\
+\lambda(\C T_n^1) \Oun (V_{B_n}(f_n) +\frac{2\pi
|\sum_0^ms_jb_n^{(j)}+b'_n|\alpha_n}{\alpha_{n-1}})\;.
\end{multline*} 
Sur $\C T_n^0$, on a la m\^eme \'egalit\'e que pour la proposition
\ref{ia}~:
$$
\int_{\C T_n^{0}}gF_{n+1}d\lambda =  
\int_{\C T_n^{0}}gF_nd\lambda
\inv{\lambda(B_n)}\int_{B_n}\Theta_nd\lambda 
+ \lambda(\C T_n^{0})\Oun(2\|(1-\Theta_n)\unde{\C T_n^{0}\cup
B_n}\|_\infty).$$
Si l'on d\'ecompose $\C T_n^2$ selon $\C I_n$, on obtient   2
sous-tours de m\^eme hauteur, $\C T_n^2\cap \C I_n$ et $\C T_n^2 \cap
\C I_n^c$~: 
la premi\`ere est de mesure $k_nb_n\alpha_n$. Sur la seconde tour qui
est incluse dans $\C I_n^c$,  on a la m\^eme \'egalit\'e que sur $\C
T_n^0$~: on obtient donc
\begin{multline*}
\int_{\C T_n^2}F_{n+1}gd\lambda= \int_{\C T_n^2} F_ngd\lambda
\inv{\lambda(B_n)}\int_{B_n}\Theta_nd\lambda \\
 + \Oun\left(2b_nk_n\alpha_n +\lambda(\C
T_n^2)(2\|(1-\Theta_n)\unde{\C I_n^c}\|_\infty)\right).
\end{multline*}
Notons de fa\c con analogue \`a la proposition~\ref{ia}, 
$$\varepsilon_n= 2 \|(1-\Theta_n)\unde{\C I_n^c}\|_\infty+2\pi
|\sum_0^m s_jb_n^{(j)}+b'_n|+2k_n/q_n.$$
Cette suite converge bien vers $0$, en raison des hypoth\`eses et du
lemme~\ref{bsj}. En remarquant que $\C T_n^0 \subset \C I_n^c$ et que
lorsque $n \in \Lambda$ on a $\lambda(\C T_n^0)\leq \alpha_nq_n\leq
b_n\alpha_n q_n$, la somme des \'egalit\'es pr\'ec\'edentes donne bien
$$\int_{\M T}F_{n+1}gd\lambda = \int_{\M T}F_ngd\lambda
\inv{\lambda(B_n)}\int_{B_n}\Theta_nd\lambda
+\Oun\left(V_{B_n}(f_n) + \varepsilon_n
(b_n\alpha_nq_n+k_n\alpha_{n-1})\right).$$  
En ce qui concerne la seconde partie de la proposition, on sait que
$\tilde F_n$ est continue et de pente constante sur $B_n$, et sa
variation 
sur $B_n$  est donc major\'ee par 
$\alpha_{n-1}|\sum_{0 \leq j\leq m}s_jk_n^{(j)}+k'_n|$.
Notons $n=n_i$ le $i$-i\`eme indice de $\Lambda$,  par le m\^eme
calcul que  dans la proposition~\ref{ia}, on trouve que
$V_{B_{n_i}}(f_{n_i})=o(b_{n_{i-1}}q_{n_{i-1}}\alpha_{n_{i-1}})$.
Enfin on a aussi $k_{n_i}= b_{n_{i-1}}q_{n_{j-1}}+k_{n_{j-1}}$~:
comme $k_n<q_n$ et que $\alpha_{n_j-1}\leq \alpha_{n_{j-1}}$ on
obtient bien l'\'egalit\'e propos\'ee.
\findem

On applique maintenant cette proposition \`a une fonction $g$, $\C
F_N$-mesurable, v\'erifiant $\int_{\M T}F_ngd\lambda \neq 0$ pour
tout $n\geq N$ (ce qui est toujours possible en raison de la
convergence forte de $(F_n)$ vers $F$ de module $1$, et de la
convergence des $(\C F_n)_n$ vers la tribu bor\'elienne).
On obtient donc que pour tout $i$ assez grand~:
$$
\int_\M T gF_{n_{i+1}}d\lambda =  
\int_\M T gF_{n_i}d\lambda
\bigg(\inv{\lambda(B_{n_i})}\int_{B_{n_i}}\Theta_{n_i}d\lambda
+
o(b_{n_{i-1}}\alpha_{n_{i-1}}q_{n_{i-1}}+b_{n_i}\alpha_{n_i}q_{n_i})\bigg).$$

Par cons\'equent, on doit avoir qu'un produit du type 
$$\prod_{i>i_0}
\left|\inv{\lambda(B_{n_i})}\int_{B_{n_i}}\Theta_{n_i}d\lambda
+
o(b_{n_{i-1}}\alpha_{n_{i-1}}q_{n_{i-1}}+b_{n_i}\alpha_{n_i}q_{n_i})\right|>0.$$
Pour terminer la preuve du th\'eor\`eme~\ref{cn}, il reste  \`a 
estimer $1-|\inv{\lambda(B_n)}\int_{B_n}\Theta_nd\lambda
+ o(b_{n-1}\alpha_{n-1}q_{n-1}+b_n\alpha_nq_n)|$ lorsque $n \in
\Lambda$. Les conditions du th\'eor\`eme~\ref{cn} d\'ecoulent
directement du r\'esultat suivant~:
\begin{lemme}\label{calcul2}
On reprend les notations pr\'ec\'edentes, et on suppose que tous les
$b_n^{(j)}$ sont positifs, que $(b_n\alpha_nq_n)_n$ et $(\sum_0^m
s_jb_n^{(j)}+b'_n)_n$ convergent vers $0$. On a alors pour tout $n
\in \Lambda$ l'\'egalit\'e~:
\begin{multline*}
1-\inv{\lambda(B_n)}\left|\int_{B_n}\Theta_nd\lambda\right|= 
(\sum_0^m s_jb_n^{(j)}+b'_n)^2 (\frac{\pi^2}{6}+o(1))\\ +
b_n \alpha_nq_n
(1+\Oun(\frac{1}{\lambda(I_n)}\int_{I_n}\Theta_nd\lambda)+o(1)). 
\end{multline*} 
 Sous l'hypoth\`ese (H), on a en plus l'in\'egalit\'e suivante~:
Il existe une constante $C<1$ strictement positive telle que pour
tout $n$ assez grand dans $\Lambda$, on ait 
$$\left|\int_{I_n}\Theta_nd\lambda\right|\leq C\lambda(I_n).$$
\end{lemme}

\dem 
Soit $n \in \Lambda$.
On d\'ecompose $B_n$ en deux sous-intervalles $B_n=I_n\cup B'_n$
comme 
pr\'ec\'edemment (cf figure~\ref{nvtour2}). Pour montrer la
premi\`ere partie du lemme,  on \'ecrit
$$
\left|\int_{B_n}\theta_nd\lambda\right|
= 
\left| \int_{B'_n}\theta_nd\lambda\right| + 
\Oun(\int_{I_n}\theta_nd\lambda).$$
Il suffit donc d'\'evaluer $|\int_{B'_n}\Theta_nd\lambda|$~:
Sur $B'_n$, $\tilde \Theta_n$ est continue, affine et de pente
constante \'egale \`a $(\sum_0^ms_j b_n^{(j)}+b'_n)q_n$. En notant
$\alpha'_n= \lambda(B'_n)$, on trouve~:    
\begin{eqnarray*}
\left|\int_{B'_n}\Theta_nd\lambda\right| 
& = &
\int_{-\alpha'_n/2}^{\alpha'_n/2}
\e{2i\pi(\sum_0^m s_jb_n^{(j)}+b'_n)q_nx}dx\\
& = & 
(\lambda(B_n)-\lambda(I_n))\left(1-\frac{\pi^2}{6}(\sum_0^m
s_jb_n^{(j)}+b'_n)^2
(q_n\alpha'_n)^2(1+o(1))\right).
\end{eqnarray*}
Par cons\'equent, comme $(\alpha'_n q_n)$ converge vers 1 et
$\lambda(I_n)=b_n\alpha_n$ on trouve bien la premi\`ere relation
annonc\'ee. 
\smallskip
Il reste  pour terminer \`a d\'emontrer la seconde partie du lemme~:
\\ 
Notons $2\varepsilon= \min\{\|\sum_Js_j\|, \; J\neq \emptyset, J \neq
\{0,..,m\}\}$. D'apr\`es l'hypoth\`ese (H), on a $\varepsilon>0$.
Pour $n \in \Lambda$ (par exemple pair), on note   $B_n=[x_n,
x_n+\alpha_{n-1}[$  et  $I_n=[x_n,x_n+b_n\alpha_n[$. Les
discontinuit\'es de $\Theta_n$ dans  $B_n$ sont de la forme
$x_n+k\alpha_n$ avec $0 <k\leq b_n$ (car lorsque $n \in \Lambda$ on
a $\min_j(b_n^{(j)})= 1$). 
Dans ce cas, on sait aussi que $b_n\geq 2$, et qu'en chaque point de
la forme $x_n+k\alpha_n$ de $I_n$, $\theta_n$ admet un saut de taille
$\sum_{j=0}^m s_j \unde{b_n^{(j)}\geq k}=\gamma_{n,k}$ qui v\'erifie
$\|\gamma_{n,k}\|\geq 2\varepsilon $  pour tout $1<k\leq b_n$ (pour $k=1$ il n'y a pas de saut car tous les $b_n^{(j)}$ sont sup\'erieurs \`a 1).\\
D'autre part, on peut supposer que pour tout $n$ assez grand  on a
$b_n\geq 3$~: en effet, sinon on aurait $b_n=2$ infiniment souvent.
Mais  $(\sum_0^m b_n^{(j)}s_j-b'_n)_n$  converge vers $0$ lorsque $n$
tend vers l'infini. Or si $b_n=2$, on a pour tout $j \in\{0,..,m\}$
$b_n^{(j)} \in \{1,2\}$, donc pour un  $n$ assez grand on doit avoir
$\sum_j b_n^{(j)}s_j\in \M Z$. Comme  $\sum_js_j=0$ on a encore $\sum_j
(b_n^{(j)}-1)s_j\in \M Z$. Notons $J=\{ j \in \{0,..,m\}, \;
b_n^{(j)}=2\}$, on obtient donc $\sum_{j \in J}s_j\in \M Z$ 
ce qui est impossible car $J$ n'est ni vide ($b_n=
\max_j(b_n^{(j)})=2$), ni $\{0,..,m\}$ tout entier (car
$\min_j(b_n^{(j)})=1$).\\
On peut d\'ecomposer maintenant l'int\'egrale sur $I_n$ de la fa\c con suivante~:
$$\left|\int_{I_n}\Theta_nd\lambda\right|\leq
 \sum_{k=1}^{\lfloor\frac{b_n-1}{2}\rfloor}\left|\int_{x_n+(2k-1)\alpha_n}^{x_n+(2k+1)\alpha_n}\Theta_nd\lambda\right|+\alpha_n(b_n-2\left\lfloor\frac{b_n-1}{2}\right\rfloor).$$
Sur les intervalles $[x_n+(2k-1)\alpha_n,x_n+(2k+1)\alpha_n]$, $\tilde \Theta_n$ est affine par morceaux de pente constante \'egale \`a $(\sum_js_jb_n^{(j)}+b'_n)q_n$, avec un saut de taille $\gamma_{n,2k}$ au milieu de l'intervalle. Par cons\'equent on peut \'ecrire 
\begin{eqnarray*}
\left|\int_{x_n+(2k-1)\alpha_n}^{x_n+(2k+1)\alpha_n}\Theta_nd\lambda\right|
&=&
\left|\int_{x_n+(2k-1)\alpha_n}^{x_n+2k\alpha_n}\Theta_nd\lambda(1+\e{2i\pi(\gamma_{n,2k}+(\sum_js_jb_n^{(j)}+b'_n)q_n\alpha_n)})\right|\\
&\leq & 
2 \alpha_n \left|\cos\left(\pi(\gamma_{n,2k}+(\sum_js_jb_n^{(j)}+b'_n)q_n\alpha_n\right)\right|.
\end{eqnarray*}
 Or on
sait que $(q_n\alpha_n(\sum_j b_n^{(j)}s_j+b'_n))_n$ converge vers 0
et que $\|\gamma_{n,2k}\|\geq 2 \varepsilon$.  Donc pour tout $n$ assez
grand on a $\|(q_n\alpha_n(\sum_j b_n^{(j)}s_j+b'_n)+\gamma_{n,2k} \|\geq \varepsilon$, d'o\`u il vient
\begin{eqnarray*}
\left|\int_{I_n}\Theta_nd\lambda\right| & \leq &
\alpha_n\left(b_n-2\left\lfloor\frac{b_n-1}{2}\right\rfloor(1-\cos \pi \varepsilon)\right)\\
&\leq &
\lambda(I_n)\left(1-(1-\cos\pi \varepsilon)(1-\frac 1 {b_n})\right)\\
&\leq&
\lambda(I_n)\cos ^2(\pi \varepsilon/2).
\end{eqnarray*}
\findem

\subsection{R\'eciproque}

On \'etablit dans ce paragraphe la preuve de la premi\`ere partie du
th\'eor\`eme~\ref{cn}, toujours selon le sch\'ema de la partie
pr\'ec\'edente.\\
Ecrivons $\phi= \sum_0^m s_j\phi_{\beta_j}$, et supposons que $\phi$
v\'erifie les hypoth\`eses suivantes~: pour une partition $\C P$, on
a pour tout $J$ de $\C P$ et pour tout $j_0\in J$ les propri\'et\'es
du th\'eor\`eme~\ref{cn}. Alors en notant $\phi_J=\sum_{j\in J}
s_j\phi_{\beta_j}$, on a $\phi=\sum_{J\in \C P}\phi_J$.  Pour  
montrer que $\phi$ est cohomologue \`a $t$, il suffit donc de
d\'emontrer que pour tout $J$, $\phi_J$ est cohomologue \`a $t_J$
(car $t=\sum_J t_J+k'\alpha$). Par cons\'equent, on peut  se
ramener au cas o\`u $\C P$ est triviale, c'est \`a dire qu'on a les
hypoth\`eses suivantes~:
on peut \'ecrire $\phi= \sum_0^m s_j\phi_{\beta_j}$ o\`u $(s_j)_j$
sont des r\'eels de somme nulle et $(\beta_j)_j$ sont des r\'eels
distincts se d\'ecomposant sous la forme $\beta_j= \sum_0^\infty
b_n^{(j)}q_n\alpha \mod 1$, avec $(b_n^{(j)})$ des suites d'entiers
v\'erifiants~:
$$\forall j \in \{0,..,m\}\quad
\sum_n\frac{|b_n^{(j)}|}{a_{n+1}}<\infty \qquad \hbox{et } \quad
\sum_n \|\sum_{j=0}^mb_n^{(j)}s_j\|^2<\infty.$$
Quitte \`a translater $\phi$ par une constante, on peut encore
supposer que tous les $(b_n^{(j)})$ sont positifs et que si
$\Lambda= \{n , \max_j(|b_n^{(j)}|)>0\}$, on a, selon le paragraphe
\ref{+},      
$$n \in \Lambda \Longleftrightarrow \left( \min_j(b_n^{(j)})=1 \hbox{
et } \max_j(b_n^{(j)})\geq 2\right).$$ 
On pose alors, toujours d'apr\`es \ref{+}, $b_n= \max_j(b_n^{(j)})$
et $\beta= \sum b_nq_n\alpha$, et on consid\`erera les domaines
fondamentaux associ\'es \`a $\beta$, $(\C D_n(\beta))$, la filtration
$(\C F_n)$ associ\'ee, ainsi que toutes les quantit\'es li\'ees.\\
Notons pour tout $n$, $b'_n=-[\sum_0^ms_jb_n^{(j)}]$, on a alors
$t=\sum_nb'_nq_n\alpha$. On construit  comme dans le chapitre
pr\'ec\'edent  une fonction $F$ qui v\'erifie (\ref{mcob}). Pour cela
on appelle $ \Theta_n$ l'unique solution  de l'\'equation
(\ref{mtheta}) qui vaut 1 sur le milieu de $B'_n$, puis on d\'emontre
que  la suite de fonctions $F_n= \prod_{k<n}\Theta_k$ admet une
valeur d'adh\'erence non nulle pour la topologie faible $L^2$.
La preuve, compl\`etement analogue \`a celle du chapitre
pr\'ec\'edent, est laiss\'ee au lecteur.

\break
\part{Constructions de facteurs discrets pour des flots sp\'eciaux et
des \'echanges de 3 intervalles.}

Dans cette partie nous d\'emontrons les th\'eor\`emes \ref{flot} et
\ref{3i} \`a l'aide des conditions  du th\'eor\`eme
\ref{marche}.  Le dernier chapitre est consacr\'e au r\'esultat de
r\'egularisation (proposition~\ref{regul}) qui permet de prouver les
parties  des th\'eor\`emes~\ref{flot} et \ref{reci}
concernant les exemples de flot sp\'eciaux avec fonctions plafond
r\'eguli\`eres.

\section{G\'en\'eralit\'es, Tours de Kakutani.}\label{kaku}

Nous rappelons ici les liens existants  entre les \'echanges de 3
intervalles et  les flots sp\'eciaux, notamment par l'interm\'edaire
des tours de Kakutani au dessus d'une rotation irrationnelle. Nous
pr\'ecisons ensuite comment s'appliquera le th\'eor\`eme~\ref{marche}
dans les chapitres suivants.
Dans ce qui suit $\alpha$ est un irrationnel dont la suite des
quotients partiels $(a_n)$ est non born\'ee.  On note toujours $T$ la
translation de $\alpha$ sur $\M T$ et $(q_n)$ la suite des
d\'enominateurs des r\'eduites de $\alpha$.

\subsection{Valeurs propres et fonctions propres des diff\'erentes
transformations.}

\subsubsection{Flots sp\'eciaux.}
On rappelle que le flot sp\'ecial, $\tau_{\alpha,\phi}$ au dessus de
$T$ et de fonction plafond $\phi$ est d\'efini sur le domaine
quotient $D_\phi=\M T \times \M R/\sim$,  pour la relation
d'\'equivalence $(x,y+\phi(x))\sim (Tx, y)$, par
$\tau_{\alpha,\phi}^t(x,y)=(x,y+t)$. Les valeurs propres de
$\tau_{\alpha,\phi}$ constituent le sous-groupe de $\M R$ d\'efini
par~:
$$e(\tau_{\alpha,\phi})=\{t \in \M R, \; \exists f \in L^2(\M
T)\backslash \{0\}\quad \e{2i\pi t\phi} f=f\circ T\}.$$ 
Si $t $ est une valeur propre du flot,  les fonctions propres
associ\'ees \`a $t$ sont de la forme $F_t(x,y)=f(x) \e{2i\pi ty}$,
o\`u $f$ est une fonction de transfert dans l'\'equation de
cohomologie associ\'ee. \\
Pour $\beta$, $\gamma$ r\'eels${}>0$ tels que $\gamma\{\beta\}<1$, on
note $\tau_{\alpha,\beta,\gamma}$ le flot sp\'ecial de fonction
plafond $\phi=1+\gamma \phi_{\beta}$.  Dans ce cas $t$ est une valeur
propre lorsqu'il existe une solution $f$ \`a l'\'equation
(\ref{eqbis}) avec $s=-t\gamma$.
Le th\'eor\`eme~\ref{marche} s'applique donc lorsque $s \notin \M Z$
avec cette condition suppl\'ementaire. 

\subsubsection{\'Echanges de 3 intervalles.}
On rappelle qu'un \'echange de 3 intervalles est identifi\'e \`a
l'induit $T_{\alpha,\beta}$ de $T$ sur $[0,\beta[$ avec $\beta \in
]0,1[$.
Son groupe des valeurs  propres est
$$e(T_{\alpha,\beta})=\{s \in \M T, \; \exists f \in L^2(\M
T)\backslash \{0\}\quad \e{2i\pi s\unde{[0,\beta[} }f=f\circ T\}.$$
$s$ est donc une valeur propre de $T_{\alpha,\beta}$  lorsque
l'\'equation (\ref{eqbis}) admet une solution avec $t=-s\beta$, et le
th\'eor\`eme~\ref{marche} s'applique donc dans ce cas. Les fonctions
propres associ\'ees \`a $s$, $f_\beta$,  sont les restrictions \`a
$[0,\beta[$ des fonctions de transfert $f$ de l'\'equation
(\ref{eqbis}).

\subsubsection{Tours de Kakutani}
On appelle (voir par exemple \cite{nad1}) tour de Kakutani au dessus
de $T$ et de fonction hauteur $h$, o\`u $h$ est une fonction
int\'egrable sur $\M T$ \`a valeurs dans $\M N^*$, le syst\`eme
dynamique $(X,T_h)$ o\`u $X=\{(x,k)\in \M T\times \M Z, \; 0\leq k
<h(x)\}$ et $T_h$ est d\'efinie pour tout $ (x,k)\in X$ par~:
$$T_h(x,k)=\left\lbrace
\begin{array}{lcl} 
(x,k+1) &\hbox{si} & k+1<h(x)\\
(Tx,0) &\hbox{sinon.}&
\end{array}\right.
$$
 On munit $X$ de la probabilit\'e invariante $\lambda_h$ \'egale \`a
la mesure de Lebesgue remormalis\'ee.
Son groupe des valeurs propres est~:
$$ e(T_h)=\{s \in \M T, \exists f \in L^2(\T)\backslash \{0\}, \quad
\e{i2\pi sh}f=f\circ T\}.$$ 
Si $s$ est une valeur propre de $T_h$, les fonctions propres $F$
associ\'ees s'\'ecrivent pour tout $(x,k)\in X$, $F(x,k)=f(x)\e{2i\pi
s k}$, o\`u $f$ est une fonction de transfert de l'\'equation
associ\'ee.\\
Pour tout $\beta>1$ donn\'e, on notera  $(X_\beta,T_{\alpha,\beta})$,
la tour de Kakutani au dessus de $T$ et de fonction hauteur
$h=\lfloor \beta\rfloor +\unde{[0,\{\beta\}[}$. 
Lorsque $\beta \in ]0,1[$, l'induit de la rotation sur $[0,\beta[$
peut s'interpr\'eter comme la tour de ``hauteur'' $\unde{[0,\beta[}$
correspondant \`a la g\'en\'eralisation naturelle au cas o\`u $h$
prend des valeurs nulles.  On notera $\lambda_\beta$ la probabilit\'e
invariante associ\'ee.

Lorsque $\beta>1$ et $\gamma ^{-1}= \beta$, le flot
$\tau_{\alpha,\beta,\gamma}$ est, \`a  homoth\'etie de temps pr\`es,
la suspension de la  tour  $T_{\alpha,\beta}$ (la suspension est le
flot sp\'ecial de fonction plafond \'egale \`a $h/\|h\|_1$).  Si $\tilde
e(T_{\alpha,\beta})$ d\'esigne le relev\'e dans $\M R$ de
$e(T_{\alpha,\beta})$, on obtient simplement
$$e( \tau_{\alpha,\beta,\beta^{-1}})= \beta \tilde
e(T_{\alpha,\beta}).$$ 
Pour montrer les parties (ii) et (iii) du th\'eor\`eme~\ref{flot}
ainsi que le th\'eor\`eme~\ref{reci}, il suffira d'\'etudier la tour
de Kakutani $T_{\alpha,\beta}$ associ\'ee au flot~:   les r\'esultats
montrant la conjugaison de $T_{\alpha,\beta}$ \`a $R_s$ ou \`a un
odom\`etre, entra\^ \i nent la conjugaison de
$\tau_{\alpha,\beta,\beta^{-1}}$ au flot de translations obtenu par
suspension, et renormalis\'e par $\beta^{-1}$, $R_{s\beta,\beta}$, ou
\`a un flot de translations sur un sol\'eno\"\i de (voir paragraphe
\ref{solenoide}).

\subsection{Les conditions $(*)$.}

D'apr\`es ce qui vient d'\^etre dit, nous obtiendrons des valeurs
propres $s$ ou $t$ pour les transformations pr\'esent\'ees, en
construisant des triplets $(\beta,s,t)$ pour lesquels l'\'equation
(\ref{eqbis}) admet des solutions, avec de plus $t=-s\beta$ ou
$s=-t\gamma$ selon les cas.

\subsubsection{D\'efinition.}

On dira que le triplet $(\beta, s,t)$ de $\M R^3$ v\'erifie $(*)$
s'il existe des suites d'entiers $(b_n)$ et $(b'_n)$ et des entiers
$k_1,k'_1,l_1$ et $l'_1$
 tels qu'on ait~:
$$(*)\qquad\left\lbrace\qquad
\begin{array}{l}
\beta  =  k_1\alpha -l_1 +\sum_1^\infty b_n\<q_n\alpha\>\\
t  =  k'_1\alpha -l'_1 +\sum_1^\infty b'_n\<q_n\alpha\>\\
\sum_1^\infty |b_n|/a_{n+1} <\infty \quad \hbox{et}\\
 \sum_1^\infty (b_ns-b'_n)^2<\infty\, .
\end{array}
\right.$$

L'avant-derni\`ere condition peut \'egalement s'\'ecrire~:
$\sum_1^\infty b_nq_n\|q_n\alpha\|<\infty$. Nous utiliserons selon
les besoins, l'une ou l'autre des  deux formes. Autrement dit, $\beta
\in H_{1}(\alpha)$, avec la d\'efinition donn\'ee dans l'introduction.

\smallskip

D'apr\`es le th\'eor\`eme~\ref{marche}, lorsque $s \notin \M Z$,  le
triplet $(\beta, s,t)$ v\'erifie $(*)$ si et seulement si
l'\'equation (\ref{eqbis}) associ\'ee au triplet $(\beta,-s,t)$ admet
une solution. 

\smallskip
Lorsque $s \in \M Z$ (pour le flot), comme $ t (1+\gamma \phi_\beta)=
t-s\beta= s(\gamma^{-1}-\beta) \mod 1$, l'\'equation se r\'eduit \`a
l'\'equation aux valeurs propres de $T$. On trouve donc des solutions
lorsque $(\gamma^{-1}-\beta) \in \M Q+\alpha\M Q$, et ces valeurs
propres  sont les multiples de
$d \gamma^{-1}$ o\`u $d = \min\{k \in \M N^* , \;
k(\gamma^{-1}-\beta)\in \M Z+\alpha\M Z\}$. 

Lorsque $\beta \notin H_1(\alpha)$, les seules valeurs propres sont
ces valeurs propres ``triviales''. 

Lorsque $\beta\in  H_1(\alpha)$, comme les solutions triviales
v\'erifient $t-s\beta\in \M Z+\alpha\M Z$ elles sont incluses dans
les solutions donn\'ees par $(*)$ et on obtient donc~:
\begin{eqnarray*}
e(\tau_{\alpha,\beta,\gamma})&=&\{t \in \M R, \; (\beta, \gamma t, t)
\hbox{ v\'erifie } (*)\},\\
\tilde e(T_{\alpha,\beta})&=& \{ s \in \M R , \; (\beta, s, s\beta)
\hbox{ v\'erifie } (*)\}.
\end{eqnarray*}
En d\'efinitive, les constructions de valeurs propres non triviales se
r\'eduisent \`a la construction de triplets $(\beta, s,t)$ v\'erifiant
$(*)$, o\`u $t= s\beta$ pour l'induit ($T_{\alpha,\beta}$ avec
$\beta\in ]0,1[$) et $s=t\gamma$ pour le flot
($\tau_{\alpha,\beta,\gamma}$).

\subsubsection{Notations et ordres de grandeurs.}\label{notation}
Pour construire des solutions \`a $(*)$, on cherchera des
 $\beta \in H_1(\alpha)$ d\'etermin\'es par des suites $(b_n)$ tr\`es
lacunaires, c'est \`a dire que les $b_n$ seront nuls sauf le long
d'une sous-suite qu'on notera $(n_j)$. Par commodit\'e, on notera
encore $(b_j)$ la sous-suite des termes non nuls . On \'ecrit alors
$$\beta=k_1\alpha-l_1+\sum_{j\geq 1}b_j\<q_{n_j}\alpha\>,$$
et on d\'efinit les suites $(\beta_j)$,
$(k_j)$ et $(l_j)$ par 
$$k_{j+1}=k_j+b_jq_{n_j},\quad l_{j+1}=l_j+b_jp_{n_j}
\quad \hbox{et } \beta_j=
k_j\alpha-l_j=\beta_1+\sum_1^{j-1}b_i\<q_{n_i}\alpha\>.$$
Enfin on imposera $|b_j|/a_{n_j+1}\le\varepsilon_j$, o\`u
$(\varepsilon_j)$ est une suite sommable strictement positive
donn\'ee.
Lorsque $(\beta,s,t)$ v\'erifie $(*)$, on peut toujours, quitte \`a
modifier le premier terme, \'ecrire $t$ sous la forme analogue
$$t=k'_1\alpha-l'_1+\sum_{j\geq 1}b'_j\<q_{n_j}\alpha\>,$$
o\`u $k'_1, l'_1$ et les $b'_j$ sont des entiers, et on notera
\'egalement $t_j=k'_1\alpha-l'_1+\sum_{1\leq
i<j}b'_iq_{n_i}\alpha=k'_j\alpha-l'_j$ pour $j\geq 1$.

Le lemme suivant pr\'ecise alors la vitesse de convergence
des sommes partielles $\beta_j$ ainsi que l'ordre de grandeur de $k_j$
et $l_j$.
\begin{lemme}\label{cv*} 
Soit  $\beta$ dans  $H_1(\alpha)$. Avec les notations
pr\'ec\'edentes,  si  $M=1+\sum_{j\geq 2} \varepsilon_j$ et
$C=(1+\max(|k_1|,|l_1|))\prod_{i\geq 1} (1+\varepsilon_i)$, on a les
in\'egalit\'es
\begin{eqnarray*}
|\beta-\beta_j| & \leq & 
M\alpha_{n_j-1}\varepsilon_j\quad\hbox{pour tout }j\geq 1,\\
\max(|k_j|,|l_j|)&\leq  &C b_{j-1}q_{n_{j-1}}
\quad\hbox{pour tout }j\geq 2.
\end{eqnarray*}
Si de plus $\beta\notin \M Z$ et que  le triplet $(\beta,s,t)$
v\'erifie $(*)$ avec  $t=s\beta$, en posant  $s_j=t_j/\beta_j$  on a
aussi  
$(s-s_j)= o(\alpha_{n_j})$.
\end{lemme}
\dem
Pour tout $j \geq 1$  on a 
$$
|\beta-\beta_j|
\leq  \sum_j^\infty|b_i|\alpha_{n_i}
\leq  \sum_j^\infty a_{n_i+1}\alpha_{n_i}\varepsilon_i
\leq \alpha_{n_j-1}\varepsilon_j +\alpha_{n_j}\sum_{j+1}^\infty
\varepsilon_i.
$$
Comme $b_j \neq 0$, on a $1\leq a_{n_j+1}\varepsilon_j$, d"o\`u
$\alpha_{n_j}<\alpha_{n_{j-1}}\varepsilon_j$ et  on
obtient bien 
$|\beta-\beta_j| \leq  M\alpha_{n_j-1}\varepsilon_j$.

En ce qui concerne la majoration de $|k_j|$, il suffit de montrer par
r\'ecurrence que pour tout $j\geq 2$ on a $|k_j|\leq
C_j|b_{j-1}|q_{n_{j-1}}$ avec 
$$C_j=(|k_1|+1)\prod_{i=1}^{j-2}(1+\varepsilon_{i}).$$
On a bien $|k_2|=|k_1+b_1q_{n_1}|\leq C_1 |b_1|q_{n_1}$.
Supposons que pour $j\geq 2$, on ait $|k_j|\leq
C_j|b_{j-1}|q_{n_{j-1}}$. Alors
$$|k_j|\leq C_j\varepsilon_{j-1}a_{n_{j-1}+1}q_{n_{j-1}}
< C_j\varepsilon_{j-1}q_{n_j}$$
et il vient
\begin{eqnarray*}
|k_{j+1}|&\leq |k_j|+|b_j|q_{n_j}\leq
C_j\varepsilon_{j-1}q_{n_j}+|b_j|q_{n_j}\\
&\leq (C_j\varepsilon_{j-1}+1)|b_j|q_{n_j}\leq C_{j+1}|b_j|q_{n_j}.
\end{eqnarray*}
On obtient donc l'in\'egalit\'e
annonc\'ee pour $|k_j|$. La majoration est analogue pour $|l_j|$.

Supposons maintenant que $(\beta,s,t)$ v\'erifie $(*)$. Comme $\beta
\notin \M Z$, alors $s_j=t_j/\beta_j$ est bien d\'efini pour $j $
assez grand et si $\varepsilon'_j=b_js-b'_j$ on peut \'ecrire
$$
\beta_j(s-s_j)  =  \beta_j s- t_j = s(\beta_j-\beta)+(t-t_j) = 
\sum_{i\geq j} {\varepsilon'}_i\<\alpha q_{n_i}\>.
$$
En remarquant que  pour tout $j$, $\sum_{i\geq j}\alpha_{n_i}^2\leq
2\alpha_{n_j}^2$, on obtient  en notant
$M'=\sum_1^\infty{\varepsilon'_j}^2$~:
\begin{eqnarray*}
|s-s_j| & \leq &
\frac 1 {|\beta_j|}
\left( \alpha_{n_j}|\varepsilon'_j|+\alpha_{n_{j+1}}\sqrt {2M'}\right)
\\ 
&\leq &\frac { {\alpha_{n_j}}} {|\beta_j|}
\left(|\varepsilon'_j| + \varepsilon_{j+1}\sqrt{2M'}\right).
\end{eqnarray*}
\findem

\section{Groupe de valeurs propres de rang  infini.}

Nous d\'emontrons dans ce chapitre le th\'eor\`eme suivant~:

\begin{theo}\label{infini}
Pour tout irrationnel $\alpha$ \`a quotients partiels non born\'es et
pour tout $\gamma>0$, on peut construire un ensemble de $\beta$ non
d\'e\-nom\-bra\-ble et dense dans $]0,1/\gamma[$ pour lesquels le flot
sp\'ecial 
$\tau_{\alpha,\beta,\gamma}$ admet comme facteur discret un flot de
translations ergodique sur un tore de dimension infinie.\\
De m\^ eme, pour tout irrationnel $\alpha$ \`a quotients partiels non
born\'es, on peut construire un ensemble de $\beta$ non d\'e\-nom\-bra\-ble
et dense dans $\M R^+$ pour lesquels $T_{\alpha,\beta}$ admet un
groupe de valeurs propres de rang infini.
\end{theo}

Le premier paragraphe fournit une construction de solutions non
triviales \`a $(*)$ dans les cas de $\tau_{\alpha,\beta,\gamma}$  ou
de $T_{\alpha,\beta}$, qui servira de base pour la suite.\\ 
Le second paragraphe montre comment on peut, \`a l'aide de cette
construction, obtenir plusieurs valeurs propres rationnellement
ind\'ependantes~: nous traiterons d'abord le cas de
$T_{\alpha,\beta}$, puis celui de $\tau_{\alpha,\beta,\gamma}$  en
pr\'ecisant \`a chaque \'etape quelles sont les diff\'erences. \\
Enfin nous expliquerons dans le troisi\`eme paragraphe comment
utiliser ce qui pr\'ec\`ede pour construire des syst\`emes admettant
une infinit\'e de valeurs propres rationnellement ind\'ependantes, ce
qui ach\`evera la preuve du th\'eor\`eme~\ref{infini}. 
 
\bigskip

Dans tout ce chapitre, on suppose donn\'e $\alpha$, un  irrationnel \`a quotients partiels non born\'es.  On suppose aussi donn\'ee une
suite $(\varepsilon_j)_{j\geq 1}$ positive strictement d\'ecroissante
et sommable, avec $\sum_{j\geq 1}\varepsilon_j\leq 1/2$.

\subsection{Premi\`ere construction de valeurs propres non
triviales.}\label{const1}

On construit ici des solutions $(\beta,s,t)$ \`a $(*)$, avec
$t=s\beta$ pour l'induit ou la tour $T_{\alpha,\beta}$ ou $s= t\gamma$
pour le flot $\tau_{\alpha,\beta,\gamma}$.  On se donne aussi une
suite $(\varepsilon'_j)$ strictement positive de carr\'e sommable.  On
part de deux entiers non nuls $k_1$ et $k'_1$ et de deux entiers
quelconques $l_1$ et $l'_1$.\\
Nous construisons par r\'ecurrence les suites $(b_j)$ et $(b'_j)$
d'entiers non nuls et la suite $(n_j)$ strictement croissante
d\'efinissant $\beta$ et $t$.  Les suites $(k_j)$, $(k'_j)$,
$(\beta_j)$, $(t_j)$ introduites au paragraphe \ref{notation} sont
alors d\'efinies simultan\'ement pour $j \geq 1$.  On notera aussi
$s_j=t_j/\beta_j$ dans le cas de $T_{\alpha,\beta}$ ($\beta_j$ sera
non nul) ou $s_j=t_j\gamma$ dans le cas de
$\tau_{\alpha,\beta,\gamma}$, pour tout $j \geq 1$.

\medskip

Soit $j\geq 1$.  Supposons avoir construit $(b_i)_{i<j}$,
$(b'_i)_{i<j}$ et $(n_i)_{i<j}$, et supposons que
$\beta_j=k_1\alpha-l_1+\sum_{1\leq i<j}b_i\<q_{n_i}\alpha\>\neq 0$ donc
que $s_j$ est bien d\'efini.
On choisit alors $b_j$ puis $b'_j$ de sorte que 
\begin{equation}\label{sol} 
\|b_js_j\|<\varepsilon'_j\quad \hbox{ et } \quad b'_j=[b_js_j].
\end{equation} 
et on choisit ensuite $n_j>n_{j-1}$ tel que
\begin{equation}\label{majbq}
q_{n_j}\geq\max{(|b_j|,|b'_j|)},
\end{equation}
et
\begin{equation}\label{majba}
\frac{\max{(|b_j|,|b'_j|)}}{a_{n_j+1}}<\varepsilon_j.
\end{equation}
Pour $j=1$, on demande de plus $n_1>1$ et
\begin{equation}\label{ci}
q_{n_1}>\max{(|k_1|,|k'_1|)}.
\end{equation}

Pour $j=1$, on a bien $\beta_1=k_1\alpha-l_1\neq 0$ puisque $k_1\neq 0$
et d'apr\`es la derni\`ere condition $\|\beta_1\|=\|k_1\alpha\|\geq
\alpha_{n_1-1}$.
Pour $j>1$, si (\ref{majba}) est satisfaite pour $i<j$ on a de m\^eme
que dans la preuve du lemme~\ref{cv*}, 
\begin{equation}\label{cvbis}
|\beta_j-\beta_1|
<\sum_{1\leq i<j}a_{n_{i+1}}\alpha_{n_i}\varepsilon_i
<\alpha_{n_1-1}\sum_{1\leq i<j}\varepsilon_i\leq\frac{\alpha_{n_1-1}}{2},
\end{equation}
d'o\`u $\beta_j\neq 0$ et on peut continuer la construction.

Dans la suite, on sera amen\'e \`a faire cette construction avec
(\ref{sol}) et (\ref{majbq}) satisfaites \`a partir d'un certain rang
seulement.  On obtient les propri\'et\'es suivantes~:

\begin{prop}\label{sol*}
Si la condition (\ref{ci}) est satisfaite et que (\ref{majba}) est
satisfaite pour tout $j\geq 1$, les suites $(\beta_j)$, $(t_j)$ et
$(s_j)$ convergent respectivement vers trois r\'eels $\beta$, $t$ non
entiers et $s\neq 0$, avec $s=t\gamma$ pour le cas de
$\tau_{\alpha,\beta,\gamma}$ ou $t=s\beta$ pour le cas de
$T_{\alpha,\beta}$.\\
D'autre part il existe un r\'eel $M(t_1,n_1)$ qui ne d\'epend que de
$t_1$ et $n_1$ tel que pour tout $j\geq 1$, $|s_j|\leq M(t_1,n_1)$. 
Pour tout $j \geq 2$,  $k_j$ et $k'_j$ sont non nuls et v\'erifient 
\begin{equation}\label{majk}
\max(|k_j|,|k'_j|)<2\varepsilon_{j-1}q_{n_j}.
\end{equation}
Si de plus les conditions (\ref{sol}) et (\ref{majbq}) sont
v\'erifi\'ees pour tout $j $ assez grand, alors le triplet $(\beta,
s,t)$ obtenu satisfait  $(*)$. L'ensemble des limites $\beta$
obtenues dans ce cas est non d\'e\-nom\-bra\-ble et dense dans $\M R^+$.
\end{prop} 

\dem 
Il est clair d'apr\`es $(\ref{majba})$ que la suite $(\beta_j)$
converge et que sa limite $\beta$ appartient \`a $H_1(\alpha)$.  Pour
tout $j\geq 1$, d'apr\`es (\ref{cvbis}), on a
$\|\beta_j\|\geq\|\beta_1\|-|\beta_j-\beta_1|\geq \alpha_{n_1-1}/2$
donc $k_j$ est non nul.  On obtient encore $\|\beta\|\geq
\alpha_{n_1-1}/2$ donc $\beta\notin \M Z$.
 
La suite $(t_j)$ converge vers $t$ pour les m\^emes raisons que
$(\beta_j)$ et on a  les m\^emes in\'egalit\'es~:
pour tout $j\geq 1$,  $|t_j-t_1|\leq \alpha_{n_1-1}/2$ et
$\|t_j\|\geq \alpha_{n_1-1}/2$, d'o\`u  $k'_j\neq 0$ et $t\notin \M
Z$. 

Donc $(s_j)$ converge aussi vers une limite $s=t/\beta$ (pour
$T_{\alpha,\beta}$) ou $s=t\gamma$ (pour $\tau_{\alpha,\beta,\gamma}$)
et dans les deux cas $s\neq 0$.  Comme $|t_j)|\leq
|t_1|+\alpha_{n_1-1}/2$ et que $|\beta_j|\geq \alpha_{n_1-1}/2$ pour
tout $j\geq 1$, la suite $(s_j)$ est bien major\'ee par une constante
$M(t_1,n_1)$ qui ne d\'epend que de $t_1$
et $n_1$.

L'in\'egalit\'e (\ref{majk}) se d\'emontre par r\'ecurrence, 
de fa\c con identique pour les suites $(k_j)$ et $(k'_j)$. 
Pour $j\ge 1$, si l'on suppose seulement que
$|k_j|<q_{n_j}$, ce qui est vrai par hypoth\`ese pour $j=1$,
il suffit d'\'ecrire ~:
$$|k_{j+1}|  =  |k_j+b_jq_{n_j}|
\leq  q_{n_j}(1+a_{n_j+1}\varepsilon_j)
< 2q_{n_j}a_{n_j+1}\varepsilon_j<2\varepsilon_jq_{n_{j+1}}
$$
car  $b_j\neq 0$ et $|b_j|< a_{n_j+1}\varepsilon_j$ donne
$a_{n_j+1}\varepsilon_j>1$.

\smallskip
Supposons maintenant que (\ref{sol}) et (\ref{majbq}) soient
satisfaites pour tout $j$ assez grand.  Pour montrer que $(\beta, s,t)$
satisfait $(*)$ il suffit de v\'erifier que la suite $((b'_j-b_js))$ est
de carr\'e sommable et, comme $(b'_j-b_js_j)$ est de carr\'e sommable
d'apr\`es (\ref{sol}), il suffit de montrer que $(b_j(s-s_j))$
l'est \'egalement.\\
D'apr\`es (\ref{majba}) et le lemme~\ref{cv*}, on a $|\beta-\beta_j|\leq 
M\alpha_{n_j-1}\varepsilon_j$ et de m\^eme $|t-t_j|\leq 
M\alpha_{n_j-1}\varepsilon_j$ pour tout $j\geq 1$, avec $M= 3/2$.
Pour le flot $\tau_{\alpha,\beta,\gamma}$ on a donc pour $j$ assez
grand
$$|b_j(s-s_j)|\leq \gamma
q_{n_j}|t-t_j|\leq M\gamma q_{n_j}\alpha_{n_j-1}\varepsilon_j< M\gamma
\varepsilon_j,
$$
donc la suite est sommable.\\
Pour $T_{\alpha,\beta}$ la majoration est semblable~:  
\begin{eqnarray*}
 |b_j(s-s_j)|&=&
\frac{|b_j|}{|\beta_j|}\left|s(\beta_j-\beta)+(t-t_j)\right|\\
&\leq&
\frac{2M}{\alpha_{n_1-1}}\left(|s|+1\right)q_{n_j}\alpha_{n_j-1}\varepsilon_j
<\frac{2M}{\alpha_{n_1-1}}\left(|s|+1\right)\varepsilon_j.
\end{eqnarray*}  
Enfin l'ensemble des $\beta$ obtenu ainsi est dense~: en effet
l'ensemble des $\beta_1=k_1\alpha-l_1$ possibles est dense, et
$|\beta-\beta_1|\leq \alpha_{n_1-1}/2$ peut \^etre rendu
arbitrairement petit par le choix de $n_1$.
Cet ensemble aussi non d\'e\-nom\-bra\-ble~: d'apr\`es le lemme
\ref{unicite}, l'unicit\'e de la d\'ecomposition en s\'erie de $\beta$
\`a partir de l'indice $j$ est assur\'ee d\`es que $|k_j|<q_{n_j}/2$
et $|b_i|<a_{n_i+1}/4$ pour tout $i \geq j$.  Ces 2 conditions sont
v\'erifi\'ees gr\^ace \`a (\ref{majba}) et (\ref{majk}) d\`es que
$\varepsilon_{j-1}<1/4$.  Comme \`a chaque \'etape il y a une infinit\'e
de choix pour $b_j$ et $n_j$, l'ensemble des sommes obtenues est bien
non d\'e\-nom\-bra\-ble.
\findem

\subsection{Construction avec $r$ valeurs propres rationnellement 
ind\'ependantes.}
 
Soit $r$ un entier${}\geq 2$.  On va construire selon le sch\'ema
pr\'ec\'edent, $\beta$ et $r$ r\'eels $t^1$,.., $t^r$ de sorte que
$T_{\alpha,\beta}$ ou $\tau_{\alpha,\beta,\gamma}$ admette $r$ valeurs
propres ind\'ependantes.
Ces valeurs propres seront, pour $T_{\alpha,\beta}$, les
$s^i=t^i/\beta$ pour $i \in \{1,..,r\}$~: on veut qu'elles soient
rationnellement ind\'ependantes dans $\M R/\M Z$, ce qui revient \`a
construire $\beta$, $t^1$,.., $t^r$ rationnellement ind\'ependants.
Pour $\tau_{\alpha,\beta,\gamma}$, les valeurs propres
seront $t^1$,.., $t^r$~: il suffit donc a priori de  construire
$t^1$,.., $t^r$  rationnellement ind\'ependants.
On posera pour la suite, par commodit\'e, $t^0=\beta$.

On suppose toujours donn\'ee la suite $(\varepsilon_j)_{j\geq 1}$
strictement positive d\'ecroissante et sommable, avec $\sum_1^\infty
\varepsilon_j\leq
1/2$. Nous prenons les notations suivantes ~:
\begin{itemize}
\item $\vect{k_1}\in{\M Z^*}^{r+1}$ et
$\vect{l_1}\in \M Z^{r+1}$.
\item $(\vect{b_j})_{j\geq 1}$ est une  suite dans ${\M Z^*}^{r+1}$ 
\item $(n_j)_{j\geq 1}$ est une suite d'entiers positifs strictement
croissante.
\end{itemize} 
Alors les suites $(\vect{k_j})_{j\geq 1}$ et $(\vect{t_j})_{j\geq 1}$
sont d\'efinies pour tout $j \geq 1$ par  
$$ \vect{k_{j+1}} = \vect{k_j}+q_{n_j}\vect{b_j},
\quad \hbox{ et }\quad \vect{t_j} =
\alpha\vect{k_1}-\vect{l_1}+\sum_1^{j-1}\<q_{n_i}\alpha\>\vect{b_i}.
$$
On note $\vect{k_j}=(k_j^i)_{0\leq i\leq r}$,
$\vect{b_j}=(b_j^i)_{0\leq i\leq r}$ et de m\^eme pour les autres
vecteurs.
On sera aussi amen\'e \`a travailler dans $\M R^r$, avec les
\'el\'ements de $\M R^{r+1}$ priv\'es de leur premi\`ere coordonn\'ee.
Nous noterons, sauf ambiguit\'e, pour tout $\vect x =(x^0,..,x^r)$ de
$\M R^{r+1}$, $x$ le vecteur $(x^1,..,x^r)$ de $\M R^r$.\\
Enfin, $\scal(.,.)$ et $\|.\|_2$ d\'esignent le produit scalaire et la norme
euclidienne usuels, indiff\'eremment dans $\M R^{r+1}$ ou $\M R^r$.

Les vecteurs $\vect{k_1}$ et $\vect{l_1}$ \'etant donn\'es, on
choisira les $\vect{b_j}$ et $n_j$ de fa\c{c}on que (\ref{ci}) et
(\ref{majba}) soient satisfaites pour tout couple $(b_j^0, b_j^i)$
avec $1\leq i\leq r$, c'est-\`a-dire
\begin{equation*}
     \max_i (|k_1^i|)<q_{n_1} \quad \hbox{et}\quad 
     \max_i(|b_j^i|)\leq \varepsilon_ja_{n_j+1}\quad
     \hbox{pour tout }j\geq 1.
\end{equation*}
Alors la premi\`ere partie de la proposition~\ref{sol*} s'applique et
en particulier la suite $(\vect{t_j})$ converge vers un vecteur $\vect
t$.\\
Pour tout $j \geq 1$, on aura comme pr\'ec\'edemment
$\beta_j=t^0_j\neq 0$ et on pose $\vect{s_j}=(t^0_j)^{-1}\vect{t_j}$
dans le cas de $T_{\alpha,\beta}$, ou $\vect{s_j}=\gamma\vect{t_j}$
dans le cas de $\tau_{\alpha,\beta,\gamma}$.
On demandera ensuite que (\ref{majbq}) et (\ref{majbq}) soient
satisfaites pour les couples $(b_j^0, b_j^i)$ pour $1\leq i\leq r$,
\`a partir d'un certain rang et avec une suite $(\varepsilon'_j)$ \`a
pr\'eciser. 

\subsubsection{Crit\`eres d'ind\'ependance.}
Nous pr\'ecisons d'abord comment obtenir l'ind\'ependance
alg\'ebrique des valeurs propres ainsi construites. Nous donnons les
\'enonc\'es pour le cas de $T_{\alpha,\beta}$.
\begin{lemme}\label{libre} 
Avec les notations et les hypoth\`eses ci-dessus, si pour tout $j$
assez grand, la famille $(\vect{t_{j-r}},..,\vect{t_j})$ est libre
dans $\M R^{r+1}$, alors les coordonn\'ees du vecteur $\vect
t=\lim\vect{t_j}$ sont rationnellement
ind\'ependantes.
\end{lemme}
\dem
Supposons qu'il existe $\vect v=(v_i)_{1\leq i\leq r}\in \M
Z^{r+1}\setminus 0$ tel que $\sum_0^r v_it^i=0$.
On obtient alors que pour tout $j$, 
$$\scal(\vect v ,\vect{t_j}) +
\sum_j^\infty \scal (\vect v,\vect{b_m})\<q_{n_m}\alpha\>=0.
$$  
On a
$\scal(\vect v ,\vect{t_j})=(\sum_{i}v_ik^i_{j})\alpha \mod 1$.
L'in\'egalit\'e (\ref{majk}) s'applique \`a
tous les $k_j^i$, ce qui donne $|\sum_{i}v_ik^i_{j}|\leq
2\varepsilon_{j-1}q_{n_j}\sum_i|v_i|$ pour tout $j \geq 2$.
Donc pour tout $j$ assez grand pour que
$2\varepsilon_{j-1}\sum_i|v_i|<1$, on a $|\sum_{i}v_ik^i_{j}|<q_{n_j}$
et si cette somme n'est pas nulle
$\|(\sum_{i}v_ik^i_{j})\alpha\|\geq\alpha_{n_j-1}$.\\
D'autre part, pour tout $j\geq 1$, comme
$\max_i(|b_j^i|)\leq \varepsilon_ja_{n_j+1}$,
$$|\sum_{m\geq j}\scal(\vect v,\vect{b_m})<q_{n_m}\alpha\>|
\leq\sum_i|v_i|\sum_{j\geq m}\varepsilon_ma_{n_m+1}\alpha_{n_m}
\leq M\sum_i|v_i| \alpha_{n_j-1}\varepsilon_j,$$ 
selon le lemme~\ref{cv*}, avec $M=3/2$. Il en
r\'esulte que $\scal(\vect v,\vect{t_j})=0$ pour tout $j$ assez grand. 
Alors la famille 
$(\vect {t_{j}},..,\vect {t_{j+r}})$ admet un orthogonal non trivial
ce qui est contraire \`a l'hypoth\`ese.
\findem

Nous utiliserons ce r\'esultat pour la construction de $(\vect{t_j})$
par l'interm\'ediaire du lemme suivant, qui s'applique lorsque la
condition (\ref{sol}) est v\'erifi\'ee.

\begin{lemme}\label{reclibre}
Soit $j>r$. On suppose que la  famille
$(\vect{t_{j-r}},..,\vect{t_j})$ 
est libre dans $\M R^{r+1}$ et que, pour un vecteur $\vect u$ orthogonal
\`a $(\vect{t_{j-r+1}},..,\vect{t_j})$, on a 
$$\sum_{i=1}^r u_i\<b_j^0s_j^i\>\neq 0, \quad \hbox{ et }\quad
b_j^i=[b_j^0s_j^i] \quad\hbox{  pour tout } i \in \{1,..,r\}.$$ 
alors la famille $ (\vect{t_{j-r+1}},..,\vect{t_{j+1}})$ est
libre.\end{lemme} 
\dem
D'apr\`es l'hypoth\`ese l'orthogonal de la famille
$(\vect{t_{j-r+1}},..,\vect{t_j})$ est de dimension 1, et est donc
engendr\'e par le vecteur $\vect u$.  Pour que la famille
$(\vect{t_{j+1-r}},..,\vect{t_{j+1}})$ soit libre, il suffit que
$\scal(\vect u,\vect{t_{j+1}}) \neq 0$.
Comme $\vect{t_{j+1}}=\vect{t_j}+\<q_{n_j}\alpha\>\vect{b_j}$ et
$\scal(\vect u,\vect{t_j})= 0$ par hypoth\`ese,
Il suffit que $\scal(\vect u,\vect {b_j})\neq 0$.\\   
Or $b_j^i=[b_j^0s_j^i]=b_j^0s_j^i-\<b_j^0s_j^i\>$ pour tout $i$.
Comme $\vect{s_j}=(t_j^0)^{-1}\vect{t_j}$, $b_j^0s_j^0=b_j^0\in\M Z$
et on a encore $\scal({b_j^0\vect{s_j}},\vect u)=0$.  Il en r\'esulte
que $\scal(\vect u,\vect{b_j})=-\sum_0^r u_i\<b_j^0s_j^i>=-\sum_1^r
u_i\<b_j^0s_j^i>\neq 0$.
\findem

Dans le cas des flots $\tau_{\alpha,\beta,\gamma}$, ces lemmes sont encore
valides \`a condition de les adapter 
dans $\M R^r$ en supprimant la premi\`ere coordonn\'ee, et en
rempla\c{c}ant $r$ par $r-1$ pour les familles de vecteurs.

\subsubsection{Construction pour $T_{\alpha,\beta}$.}

Nous proc\'edons maintenant \`a la construction pour le cas des
induits ou des tours $T_{\alpha,\beta}$.
\begin{prop}\label{existe}
Soient $\vect{k_1} \in {\M Z^*}^{r+1}$, $\vect{l_1}\in \M Z$ et
$n_1\in \M N^*$ v\'erifiant $\max_{0 \leq i\leq r}{|k_1^i|}<q_{n_1}$
et $J_{r}$ un entier sup\'erieur \`a $r$.  Alors il existe une
constante $C_r>0$ ne d\'ependant que de $\vect{k_1}$ et $n_{1}$ telle
que pour toute suite $(\varepsilon'_j)$ v\'erifiant
$\varepsilon'_j\geq C_r \varepsilon_{j-r-1}$ pour tout $j>J_r$, on
peut constuire une suite $(\vect{b_j})$ dans ${\M Z^*}^{r+1}$ et une
suite strictement croissante $(n_j)$ dans $\M N$ satisfaisant les
propri\'et\'es suivantes.
\begin{enumerate}
\item\label{cii}
La famille $ (\vect{t_{J_r-r+1}},..,\vect{t_{J_r+1}})$ est libre.
\item\label{majbi} 
 Pour tout $j\geq 1$,  $\max_{0 \leq i\leq
r}{|b_j^i|}<\min{(\varepsilon_ja_{n_j+1},q_{n_j})}$.
\item\label{choix}
 Pour tout $j>J_r$,    $\DSM\max_{0\leq i\leq
r}\|b_j^0s_j^i\|\leq\varepsilon'_j$, \\
 il existe       $\vect{u}$  orthogonal \`a
$\{\vect{t_{j-r+1}},..,\vect{t_j}\}$ tel que $\DSM\;\sum_{i=1}^r
u_i\<b_j^0s_j^i\>\neq 0$  \\
et $ b_j^i=[b_j^0s_j^i]$ pour tout  $i \in \{1, .. ,r\}$. 
\end{enumerate}
Alors la suite $(\vect{t_j})$ ainsi construite converge vers un
vecteur $(t^0, t^1, ..,t^r)$ de coordonn\'ees rationnellement
ind\'ependantes.  Lorsque $\sum_{j\geq 1}{\varepsilon'_j}^2<+\infty$
tous les triplets $(t^0,t^i/t^0, t^i)$ sont solutions de $(*)$, et
l'ensemble des limites $t^0$ ainsi obtenues est non d\'e\-nom\-bra\-ble et
dense.
\end{prop}

%\rem On peut choisir $C_r=2M_rc_r\sqrt r$ o\`u $c_r$ est la
% constante d\'efinie au  lemme~\ref{G} et
% $M_r=\max_i(M(t_1^i,n_1))+1$. 

\dem 
La seconde partie de la proposition r\'esulte clairement des lemmes
\ref{libre} et \ref{reclibre} et de la proposition~\ref{sol*}.
Soient $J_{r}\geq r$ et $(\varepsilon'_{j})$  v\'erifiant l'hypoth\`ese
de l'\'enonc\'e, la constante $C_{r}$ restant \`a choisir. Il
suffit  de d\'emontrer l'existence des suites $(\vect{b_{j}})$,
$(n_{j})$ ou, de fa\c con \'equivalente, d'une suite $(\vect{t_{j}})$
v\'erifiant (i), (ii) et (iii).\\
On commence  par se donner $(\vect{t_j})_{j\leq J_r-r+1}$ de sorte
que 
(\ref{majbi}) soit v\'erifi\'ee pour $j\leq J_r-r$~;
ensuite on choisit  $(\vect{b_{J_r-r+1}},..,\vect{b_{J_r}})$ tels que 
 $(\vect{t_{J_r-r+1}},\vect{b_{J_r-r+1}},..,\vect{b_{J_r}})$ soit
libre dans $\M R^{r+1}$. Comme $\vect
{t_{j+1}}=\vect{t_j}+\<q_{n_j}\alpha\>\vect{b_j}$, on obtient ainsi
la condition (\ref{cii}). Enfin on choisit 
$(n_j)_{j\leq J_r}$ de sorte qu'on ait  (\ref{majbi}) pour $j\leq
J_r$.\\
Soit maintenant $J>J_r$. Supposons  avoir construit $\vect{b_j}$ et
$n_{j}$ v\'erifiant les conditions (ii) et (iii) jusqu'au rang $J-1$,
c'est \`a dire qu'on a d\'efini $(t_{j})_{j\leq J}$. Pour construire
$\vect{t_{J+1}}$,
il suffit de trouver $b^0_J$  satisfaisant  les premi\`eres
conditions de (\ref{choix}), de d\'efinir $\vect{b_J}$ par la
derni\`ere condition de (\ref{choix}), puis  de choisir $n_J$ de
sorte que (\ref{majbi}) soit vraie au rang $J$.  
 Il reste donc \`a montrer l'existence de $b_J^0$.\\
 Soit $\vect u$ un vecteur non nul  orthogonal \`a
$\{\vect{t_{J-r+1}},..,\vect{t_J}\}$. Comme $k_j^0\neq 0$ d'apr\`es
la proposition~\ref{sol*}, $u=(u_i)_{i\geq 1}$ est un vecteur non nul
de $\M R^r$ qu'on peut choisir de sorte que $\|u\|_2=1$. 
Notons $x_{i}=\<b_{J}^0s_{J}^i\>$, alors les premi\`eres conditions
de (\ref{choix}) s'\'ecrivent~:
 $$\max_{1\leq i\leq r}{|x_i|}<\varepsilon'_J \quad\hbox{ et }\quad
\sum_1^r u_ix_i\neq 0.$$
Par densit\'e pour montrer l'existence de $b_{J}^0$,  il suffit de
montrer que le relev\'e dans $\M R^r$ du groupe ferm\'e engendr\'e
dans $\M T^r$ par $s_J=(s_J^1, ..,s_J^r)$, not\'e  $G$,  est
$(\varepsilon'_J/3)$-dense dans $\M R^r$ pour la norme euclidienne. 
En effet,
dans ce cas, on peut choisir $x\in G$ tel que 
$\|x-\frac{\varepsilon'_J}{2}u\|_2\leq \varepsilon'_J/3$. Alors 
$\|x\|_2<\varepsilon'_J$ et 
$$|\scal(x,u)|\geq  \varepsilon'_J/2 -\|x-\frac{\varepsilon'_J}{2} u \|_2 
\geq \varepsilon'_J/6.$$
Soit  $R=\{v\in\Z^r:\;\scal(v,s_J)\in\Z\}$, alors $G=\{x\in\R^r:\;
\scal(v,x)\in\Z, \; \forall v\in R\}$  par dualit\'e des groupes
localement compacts, autrement dit $G$ est  l'annihilateur de $R$
dans $\R^r$.
La densit\'e relative de $G$, avec une constante de la forme
$C_r\varepsilon_{J-r-1}/3$, r\'esulte  des 2 lemmes suivants. 
\findem

\begin{lemme}\label{R}
On reprend les notations et  les hypoth\`eses de la proposition
\ref{existe} et  on suppose que (i), (ii) et (iii) sont v\'erifi\'ees
pour $j< J$ . Notons   $M_r=\max_i(M(t_1^i,n_1))+1$ o\`u
$M(t_1^i,n_1)$ est d\'efini 
selon la proposition~\ref{sol*}. 
Alors  pour tout $v\in R\setminus\{0\}$ on a $\|v\|_2\ge 1/(2M_r\sqrt
r\varepsilon_{J-r-1})$. 
\end{lemme}

\dem
Soit $v\in R\setminus\{0\}$.
Alors il existe un entier $v_0$ v\'erifiant $\sum_{i=0}^rv_is_J^i=0$
($s_J^0=1$). 
Comme pour tout $i \in \{0,..,r\}$ on a $s_J^i=t^i_J/t^0_j$ il vient 
$\sum_0^rv_it_J^i=0$. 
Or, avec les  hypoth\`eses, le lemme~\ref{reclibre} montre que  la
famille 
$(\vect{t_{J-r}},..,\vect{t_J})$   est libre, donc il existe un
entier $m$ entre
$J-r$ et $J-1$ tel que $\sum_{i=0}^rv_it_{m}^i\ne 0$ mais
$\sum_{i=0}^rv_it_{m+1}^i= 0$, c'est-\`a-dire
$$
\sum_{i=0}^rv_it_m^i=-\left(\sum_{i=0}^r
v_ib_m^i\right)\<q_{n_m}\alpha\>
\neq 0.
$$
Comme $ \sum_{i=0}^r v_it_m^i= (\sum_{i=0}^rv_ik_m^i) \alpha \mod 1$,
on doit donc avoir pour un entier $m \in \{J-r,..,J\}$ l'une des deux
conditions suivantes ( cf lemme~\ref{unicite}).
 \begin{eqnarray*}
q_{n_m}\leq
|\sum_{i=0}^rv_ik_m^i|<2\sum_{i=0}^r|v_i|\varepsilon_{m-1}q_{n_m} \quad
\hbox{(d'apr\`es (\ref{majk}))},& &\\ 
\hbox{ou bien }\quad a_{n_m+1}<2|\sum_{i=0}^rv_ib_m^i|\leq 2
\sum_{i=0}^r|v_i|\varepsilon_m a_{n_m+1}\quad \hbox{(d'apr\`es
(\ref{majbi}))}.& &
\end{eqnarray*}
On obtient donc n\'ecessairement $\sum_0^r|v_i|\geq
\inv{2\varepsilon_{J-r-1}}$.  Enfin, comme $v_0+\sum_1^rv_is^i_J=0$ on a

$|v_0|\leq (\max_i |s_J^i|) \sum_1^r|v_i|$ et par cons\'equent  
$$\|v\|_2\geq\inv{2M_r\sqrt r \varepsilon_{J-r-1}}.$$ 
 \findem

Le  lemme  qui suit est  un r\'esultat g\'en\'eral sur la densit\'e
relative de l'annihilateur d'un sous-groupe discret de $\M R^r$.

\begin{lemme}\label{G} 
Pour tout $r\ge 1$, il existe une constante $c_r$
telle que si $R$ est un sous-groupe discret de $\R^r$ et
$d=\min(\|v\|_2:\;v\in R\setminus\{0\})$, son annihilateur $G$ est
$(c_r/d)$-dense dans $\R^r$.
\end{lemme}

\dem
Le r\'esultat est \'evident pour $r=1$ avec $c_r=1/2$. Soit $r>1$, on
suppose la propri\'et\'e montr\'ee pour $r-1$, avec la constante
$c_{r-1}$. Soit $a\in R$ avec $|a|=d$ et soit $H$ l'hyperplan
orthogonal \`a $a$. Alors $G'=G\cap H$ est l'annihilateur dans $H$
de la projection orthogonale $R'$ de $R$ sur $H$ 
Comme
$R$ contient $\Z a$, il ne contient pas d'autres points \`a
distance inf\'erieure \`a $d\sqrt3/2$ de la droite $\R a$, donc les
points de
$R'$ sont espac\'es d'au moins $d\sqrt3/2$ et d'apr\`es
l'hypoth\`ese de r\'ecurrence $G'$ est $(2c_{r-1}/d\sqrt3)$-dense
dans $H$.

D'autre part, on obtient une base de $R$ sur $\M Z$ en ajoutant $a$
au relev\'e d'une base de $R'$. Il en r\'esulte qu'il existe $x\in
G$ avec $\scal(a,x)=1$ 
%({\sl on peut aussi dire que, si on
%d\'ecompose $a$ sur une base de $R$, les coefficients sont
%premiers entre eux, \`a cause de la minimalit\'e de $|a|$, et en
%d\'eduire la m\^eme chose}). 
Alors les hyperplans $\scal(a,x)=b$
pour $b\in\Z$, qui sont distants de $1/d$,  et chacun rencontre  $G$
selon un translat\'e de $G'$ qui y est $(2c_{r-1}/d\sqrt3)$-dense.
Cela prouve le r\'esultat avec $c_r^2=\inv4+\frac43 c_{r-1}^2$.
\findem

\rem Dans la preuve de la proposition~\ref{libre}, le groupe $G$ admet
donc la constante de densit\'e relative $2M_{r}\sqrt r
c_{r}\varepsilon_{J-r-1}$ et il suffit de choisir $C_r=6M_r\sqrt{r} c_r$.

\subsubsection{Construction pour $\tau_{\alpha,\beta,\gamma}$.}
Nous reprenons dans ce paragraphe le sch\'ema  pr\'ec\'edent, ainsi
que ses notations. Ici $\gamma$ est un r\'eel fix\'e strictement
positif.\\
Nous devons distinguer les deux cas suivants.

\medskip
\noindent
{\bf Cas $\gamma^{-1}\notin {\M Q} \alpha+{\M Q}$.}
\\
On construit la suite $(\vect{t_j})$ de sorte que $(t^1,..,t^r)$
soient rationnellement ind\'ependants. 
Nous utilisons \`a nouveau les lemmes \ref{libre} et \ref{reclibre},
en supprimant la premi\`ere coordonn\'ee. Les conditions de la
proposition~\ref{existe} deviennent~:
{\it\begin{enumerate}
\item  
 $(t_{J_r-r+1},..,t_{J_r})$ libre dans $\M R^r$,
\item 
Pour tout $j\geq 1$,  $\max_{0 \leq i\leq
r}{|b_j^i|}<\min{(\varepsilon_ja_{n_j+1},q_{n_j})}$.
\item
 Pour tout $j\geq J_r$,  $\DSM\max_{0\leq i\leq
r}\|b_j^0s_j^i\|\leq\varepsilon'_j$,   il existe  $u$  orthogonal \`a
$\{t_{j-r+2},..,t_j\}$ tel que  $\DSM\;\sum_{i=1}^r
u_i\<b_j^0s_j^i\>\neq 0$  et  $ b_j^i=[b_j^0s_j^i]$ pour tout  $i \in
\{1, .. ,r\}$.
 \end{enumerate}}
La seule modification notable concerne la preuve du lemme~\ref{R}. En
effet, 
si $v \in R$, il existe $v_0 \in \M Z$ tel que $\sum_{i=1}^r
v_is_J^i+v_0=0$. Mais ici on a $s_J=\gamma t_J$, d'o\`u la relation 
\begin{equation}\label{Rflot}
\sum_{i=1}^r v_it_J^i+v_0\gamma^{-1}=0.
\end{equation} 
Comme $\gamma^{-1}$ n'est pas dans $\alpha {\M Z}+{\M Z}$, on a
n\'ecessairement 
$v_0=0$, et le reste de la preuve du lemme est inchang\'e. On
obtient dans ce cas une constante $C_r$ qui ne d\'epend que de $r$.

\medskip
\noindent
{\bf Cas $\gamma^{-1} \in \alpha {\M Q}+{\M Q}$.}
\\
L'hypoth\`ese de r\'ecurrence utilis\'ee pour le lemme~\ref{R} sur la
libert\'e de la famille $(t_{J-r+2},..,t_{J})$  ne permet pas
d'exploiter la relation (\ref{Rflot}) comme ci-dessus. Posons
$\gamma^{-1}=(k\alpha-l)/d$ o\`u $k,l$ et $d$ sont premiers entre
eux. Nous construisons dans ce cas la suite $(\vect{t_j})$ de sorte
qu'en notant $\vect{t'_j}=(k\alpha-l,t_j)$, les familles 
$(\vect{t'_{j-r}},..,\vect{t'_j})$ soient libres.
On peut encore appliquer le lemme~\ref{libre} pour $(\vect{ t'_j})$,
\`a condition d'imposer $q_{n_1}>|k|$. On utilise alors le lemme
\ref{reclibre} pour 
$(\vect{ t'_j})$ en  modifiant la condition sur le premier terme.
L'\'enonc\'e  devient alors~:\\
{\bf Lemme~\ref{reclibre} modifi\'e : }{\it
Soit $j>r$. On suppose que la  famille
$(\vect{t'_{j-r}},..,\vect{t'_j})$ 
est libre dans $\M R^{r+1}$. Si pour un vecteur $\vect u$ orthogonal
\`a $(\vect{t'_{j-r+1}},..,\vect{t'_j})$ on a 
$$\hbox{ pour tout } i \in \{1,..,r\}, b_j^i=[b_j^0s_j^i], \hbox{ et
}\quad
\sum_{i=1}^ru_i\<b_j^0s_j^i\>+b_j^0d u_0\neq 0,$$  
alors la famille $ (\vect{t'_{j-r+1}},..,\vect{t'_{j+1}})$ est libre.
}\\
Et la proposition~\ref{existe} devient~:\\
{\bf  Proposition~\ref{existe} modifi\'ee :} {\it
Soit  $\vect{k_1} \in (\M Z\backslash \{0\})^{r+1}$ et  $n_1 \in \M
N^*$ v\'erifiant    $q_{n_1}>|k|$ et $\max_{0 \leq i\leq
r}{|k_1^i|}<q_{n_1}$. Alors  il existe une constante $C_r>0$ telle
que pour toute suite $(\varepsilon'_j)$ v\'erifiant  $\varepsilon'_j\geq
C_r \varepsilon_{j-r-1}$  pour tout $j>J_r$, on peut constuire  des
suites   $(\vect{b_j})$ et $(n_j)$  satisfaisant les propri\'et\'es
suivantes.
\begin{enumerate}
\item
La famille $ (\vect{t'_{J_r-r+1}},..,\vect{t'_{J_r+1}})$ est libre.
\item
 Pour tout $j\geq 1$,  $\max_{0 \leq i\leq
r}{|b_j^i|}<\min{(\varepsilon_ja_{n_j+1},q_{n_j})}$.
\item
 Pour tout $j>J_r$,   $\DSM\max_{0\leq i\leq
r}\|b_j^0s_j^i\|\leq\varepsilon'_j$, \\
 il existe  $\vect{u}$  orthogonal \`a
$\{\vect{t'_{j-r+1}},..,\vect{t'_j}\}$  tel que   $\DSM\;\sum_{i=1}^r
u_i\<b_j^0s_j^i\>+u_0db_j^0\neq 0$ \\
 et  $ b_j^i=[b_j^0s_j^i]$ pour tout  $i \in \{1, .. ,r\}$. 
\end{enumerate}}
\smallskip
\dem
On reprend les notations et le sch\'ema de la preuve de la
proposition~\ref{existe}~: on distingue encore 2 cas.\\
Soit $u_0 =0$ et la condition (\ref{choix}) est alors inchang\'ee. 
Dans ce cas la preuve est identique \`a celle de la proposition
\ref{existe}. Pour le lemme~\ref{R}, de la relation
$\sum_{0}^rv_{i}s_{J}^i=0$  on obtient l'\'equation (\ref{Rflot})
d'o\`u il r\'esulte que $d$ divise $v_0$. En rempla\c cant $v_{0}$
par $v_{0}/d$, on retrouve bien $\sum_{0}^r v_{i}t_{J}^i=0$ et le
reste de la preuve est exactement le m\^eme. \\
Soit $u_0 \neq 0$ et la preuve de la proposition est directe~:
il suffit de voir qu'il existe des entiers $b$ 
arbitrairement grands  v\'erifiant pour tout $1\leq i\leq r$, 
$\|bs_J^i\|<\varepsilon'_J$,  car  $|\sum_1^ru_i\<bs_J^i\>|<
\sum_1^r|u_i|$ reste born\'ee contrairement \`a $du_{0}b$. 
\findem
 
\subsection{Une infinit\'e de valeurs propres ind\'ependantes.}

Nous ne consid\'ererons que le cas de $T_{\alpha,\beta}$. La
construction pour 
$\tau_{\alpha,\beta,\gamma}$ s'adapte sans difficult\'e.\\
Nous reprenons les notations du paragraphe pr\'ec\'edent.
Il s'agit de modifier la construction pr\'ec\'edente de fa\c con \`a
obtenir une infinit\'e de valeurs propres rationnellement
ind\'ependantes. C'est \`a dire qu'on va construire $(t_j^i)_{i\geq
0\atop j\geq 1}$ telle que pour tout $r \geq 1$, la suite
$(t^0_j,..,t^r_j)_j$ v\'erifie  les conditions  de la proposition
\ref{existe} pour un entier $J_r>J_{r-1}$ choisi ci-apr\`es. \\
Supposons avoir constuit une telle suite, d'apr\`es la proposition
\ref{existe}, la suite $(\varepsilon'_j)$ v\'erifie alors pour tout
$r\geq 1$ et pour tout 
$j \in \{J_r+1,..,J_{r+1}\}$
$$\varepsilon'_j= C_r \varepsilon_{j-r-1}.$$
Par cons\'equent, pour que $(\varepsilon'_j)$ soit de carr\'e sommable, 
il est d'abord n\'ecessaire que $(C_r\varepsilon_{J_r-r})$ converge vers
$0$~: posons par exemple $J_r=2r$ pour tout $r\geq 1$, dans ce cas la
suite $(\varepsilon'_j)$ est de carr\'e sommable d\`es que
$(C_r\varepsilon_r)$ l'est.
\medskip
 
On suppose maintenant donn\'ee une suite $(\varepsilon_j)$ positive
strictement d\'ecroissante telle que $(c_r\sqrt r \varepsilon_r)$ soit
de carr\'e sommable ($c_r$ est la constante du lemme~\ref{G}). Il est
imm\'ediat de v\'erifier que dans ce cas $(\varepsilon_j)$ est sommable.
On choisit ensuite pour tout $r\geq 1 $ de poser $J_r=2r$. \\
Soit $r\geq 1$ fix\'e, supposons construite la suite
$(t^0_j,..,t_j^r)_{1\leq j\leq J_{r+1}}$  selon le paragraphe
pr\'ec\'edent. Les conditions 
(\ref{cii}), (\ref{majbi}) et (\ref{choix}) de la proposition
\ref{existe} sont alors v\'erifi\'ees pour tout $j\leq J_{r+1}$ et
$i\leq r$.\\
Pour passer \`a $r+1$, il suffit de construire $(t^{r+1}_j)_{1\leq
j\leq J_{r+1}+1}$  de sorte qu'on ait encore (\ref{majbi}) et
(\ref{cii})~: on pourra ensuite continuer la construction selon le
paragraphe pr\'ec\'edent. 
On pose pour tout $j\leq J_{r+1}$, $t_j^{r+1}=t_j^r$~: ceci assure la
condition (\ref{majbi}) pour $j<J_{r+1}$. On obtient \'egalement
$M_{r+1}=M_r=M$ ind\'ependant de $r$,  d'o\`u $C_r=2M\sqrt r c_r$ (cf
proposition~\ref{existe}). Ceci  assure la convergence de la s\'erie
$((C_r\varepsilon_r)^2)_r$, donc  que  $\sum_j{\varepsilon'_j}^2<\infty$.
\\  
Pour v\'erifier (\ref{cii}), il reste \`a choisir $b^{r+1}_{J_{r+1}}$
de sorte que la matrice
$$
\begin{pmatrix}
t_{J_{r+1}-r}^0&\ldots&t_{J_{r+1}-r}^r&t_{J_{r+1}-r}^{r+1}\\
\vdots&\ddots&\vdots&\vdots\\
t_{J_{r+1}}^0&\ldots&t_{J_{r+1}}^r&t_{J_{r+1}}^{r+1}\\
t_{J_{r+1}+1}^0&\ldots&t^r_{J_{r+1}+1}&t_{J_{r+1}+1}^{r+1}
\end{pmatrix}
$$
soit r\'eguli\`ere. Or, par construction, le premier 
mineur d'ordre r est non nul, donc  en d\'eveloppant le d\'eterminant
de la matrice par rapport \`a la derni\`ere ligne on trouve que
celui-ci s'annule pour au plus une seule valeur de
$t^{r+1}_{J_{r+1}+1}$. Comme
$t^{r+1}_{J_{r+1}+1}=t^{r+1}_{J_{r+1}}+b^{r+1}_{J_{r+1}}\<q_{n_{J_{r+1}}}\alpha\>$,
ceci montre que l'une des deux valeurs $b^{r+1}_{J_{r+1}}=\pm 1$
convient. 
La condition (\ref{majbi}) est alors automatiquement v\'erifi\'ee
pour $J_{r+1}$. \\
On obtient donc par r\'ecurrence la suite $(t_j^i)_{i,j}$ voulue. 
Ceci cl\^ot la preuve du th\'eor\`eme~\ref{infini}.

\section{Constructions de flots speciaux et d'\'echanges de 3
intervalles \`a spectre purement discret.}

Le but de ce chapitre est de donner les constructions des parties
(ii) et (iii) des th\'eor\`emes \ref{flot} et \ref{3i}, ainsi que
leurs r\'eciproques (th\'eor\`eme~\ref{reci} et (iv) et (v) du
th\'eor\`eme~\ref{3i}). Pour les flots sp\'eciaux, nous restreindrons
ici notre \'etude aux  cas des flots $\tau_{\alpha, \beta,
\beta^{-1}}$. 
D'apr\`es ce que nous avons fait remarquer au chapitre~\ref{kaku},
il nous suffira donc de construire des tours de Kakutani
$T_{\alpha,\beta}$ \`a spectre discret. 
Nous donnons dans le premier paragraphe des constructions dans
lesquelles on peut identifier le facteur Kronecker des tours de
Kakutani obtenues.  Dans les paragraphes suivants, nous montrons que
pour certains cas de ces constructions, $T_{\alpha,\beta}$ est
isomorphe \`a son facteur de Kronecker.

\subsection{Constructions avec facteur de Kronecker
explicite.}\label{kro}

Soit $\alpha$ un irrationnel et $\beta \in H_{1}(\alpha)$. D'apr\`es
le chapitre~\ref{kaku}, le relev\'e dans $\M R$ du groupe des valeurs
propres de $T_{\alpha,\beta}$, $\tilde e(T_{\alpha,\beta})$, est
l'ensemble de tous les r\'eels $s$ tels que $(\beta,s,t)$ v\'erifie
$(*)$ avec $t=s\beta$. 
On obtient le r\'esultat suivant~:

\begin{theo}\label{identifier} 
Pour tout $\alpha$ \`a quotients partiels non born\'es, il existe un
ensemble de $\beta$ dense dans $\M R^+$ pour lesquels $\tilde
e(T_{\alpha,\beta})$ est isomorphe \`a la limite inductive d'une suite
de groupes $\M Z^2$ par une suite de matrices $(A_n)$ \`a coefficients
entiers.
\end{theo}
 
%\rem{
%Dans ce cadre, le groupe des valeurs propres de $T_{\alpha,\beta}$
%est, contrairement au chapitre pr\'ec\'edent, de rang  1 au plus.}

\subsubsection{Analyse de la construction.}
On  reprend  les notations du paragraphe \ref{notation}~:
\begin{itemize}
\item $\alpha \notin \M Q$ est  donn\'e, \`a quotients partiels
non born\'es.
\item  $(n_j)_{j\geq 1}$ est une suite d'entiers positifs strictement
croissante,
\item  $(b_j)_{j\geq 1}$ et $(b'_j)_{j \geq 1}$ sont deux suites
d'entiers non nuls, et on suppose toujours que
$|b_j|/a_{n_j+1}=\O(\varepsilon_j)$, o\`u $(\varepsilon_j)$ est
une suite sommable.
\item $k_1$, $k'_1$ sont deux entiers non nuls, et $l_1$ et $l'_1$ deux
entiers quelconques.
\item Les suites $(k_j)_{j\geq 1}$, $(k'_j)_{j\geq 1}$, $(l_j)_{j\geq
1}$, $(l'_j)_{j\geq 1}$, $(\beta_j)_{j\geq 1}$, $(t_j)_{j\geq 1}$
sont d\'efinies  pour tout $j \geq 1$ par
\end{itemize}
$$
\begin{array}{lll}
k_{j+1}= k_j+b_jq_{n_j},\quad & l_{j+1}= l_j+b_jp_{n_j},\quad &
\beta_j= k_j\alpha-l_j,\\
k'_{j+1}= k'_j+b'_jq_{n_j}, \quad& l'_{j+1}= l'_j+b'_jp_{n_j},
\quad& 
t_j= k'_j\alpha-l'_j.
\end{array}$$
On suppose aussi que $b'_j=[b_js]=[b_jt\beta^{-1}]$ pour tout $j$ assez
grand.  Alors les suites $(\beta_j)$ et $(t_j)$ convergent et on note
$\beta=\lim \beta_j$, $t=\lim t_j$.  Il s'agit de trouver des
conditions sur les suites $(b_j)$ et $(n_j)$ permettant de
caract\'eriser les suites $(b'_j)$ telles que $(\beta, s,t)$ v\'erifie
$(*)$ avec $t=s\beta$.\\
\`A l'aide des estimations du lemme~\ref{cv*}, on peut \'ecrire ~:
\begin{eqnarray*}
b _js&=&b_jt_j\beta^{-1} +b_j(s-t_j\beta^{-1})\\
&= & (k'_j \alpha-l'_j)b_j\beta^{-1} + \beta^{-1} b_j(t-t_j)\\
&=&
k'_j[b_j\beta^{-1}\alpha]-l'_j[b_j\beta^{-1}]+k'_j\<b_j\beta^{-1}\alpha\>
-l'_ j\<b_j\beta^{-1}\>+ \O(b_j\alpha_{n_j-1}\varepsilon_j)\\
&=&
k'_j[b_j\beta^{-1}\alpha]-l'_j[b_j\beta^{-1}]+
\O(b_{j-1}q_{n_{j-1}}(\|b_j\beta^{-1}\alpha\|+\|b_j\beta^{-1}\|))\\
&&+\O(b_j\alpha_{n_j-1}\varepsilon_j).
\end{eqnarray*}
Supposons alors que $(b_j)$ et $(n_j)$  v\'erifient de plus~:
$$|b_j| =\O(q_{n_j}) \quad \hbox{ et }\quad 
q_{n_{j- 1}}b_{j-1}\max{(\|b_j\beta^{-1}\|,\|b_j\beta^{-1}\alpha\|)}
=\O(\varepsilon_j).$$ 
Dans ce cas, en posant $d_j=[b_j\beta^{-1}]$  et
$c_j=[d_j\alpha]$, on obtient pour $j$ assez 
grand $[b_j\beta^{-1}\alpha]=c_j$, d'o\`u 
$$b_js=  k'_jc_j-l'_jd_j+\O(\varepsilon_j) \qquad \hbox{ et }
\qquad b'_j= k'_jc_j-l'_jd_j.$$
Ceci donne $k'_{j+1}= k'_j+(k'_jc_j-l'_jd_j)q_{n_j}$ et
$l'_{j+1}= l'_j+(k'_jc_j-l'_jd_j)p_{n_j}$.

Par cons\'equent les valeurs $t$ telles que $(\beta,s,t)$ v\'erifient
$(*)$ sont des limites des suites $(k'_j\alpha-l'_j)$ o\`u les
vecteurs $\begin{pmatrix}k'_j\\ l'_j \end{pmatrix}$ satisfont pour
tout $j$ assez grand la r\'ecurrence:
$$
\begin{pmatrix}
k'_{j+1}\\l'_{j+1}\end{pmatrix}
=A_j \begin{pmatrix}k'_{j}\\l'_{j}\end{pmatrix},
\quad \hbox {avec }\quad 
A_j=\begin{pmatrix}
1+c_jq_{n_j} & -d_jq_{n_j}\\ {c_j} p_{n_j} & 1-d_jp_{n_j}
\end{pmatrix}.
$$
La valeur propre triviale $s=1$ de $T_{\alpha,\beta}$ correspond \`a
$t=\beta$ et $(\beta,1,\beta)$ est une solution. Du coup $\beta$
s'obtient \'egalement par la relation de r\'ecurrence ci-dessus. On
pourra donc partir de la donn\'ee des suites $(d_j)$ et $(n_j)$.

Comme
$\|d_j\alpha\|=\|b_j\beta^{-1}\alpha-\<b_j\beta^{-1}\>\alpha\|\leq
\|b_j\beta^{-1}\alpha\|+\|b_j\beta^{-1}\|$, ces suites doivent
v\'erifier
$$
d_j= \O(\max(q_{n_j},\varepsilon_ja_{n_j+1}))\qquad \hbox{et }\qquad
d_{j-1}q_{n_{j-1}}\|d_j\alpha\|=\O(\varepsilon_j).
$$

\subsubsection{Construction.}
Il s'agit maintenant de montrer que les conditions pr\'ec\'edentes sur
$(d_j)$ et $(n_j)$ permettent de d\'efinir un ensemble de r\'eels
$\beta$ pour lesquels le groupe des r\'eels $t$ tels que
$(\beta,t/\beta,t)$ v\'erifie $(*)$ est isomorphe \`a la limite
inductive de la suite des groupes $\M Z^2$ avec les matrices de
transition $(A_j)$.  On notera $\lim_{\rightarrow}(\M Z^2,A_{j})$
cette limite, qu'on identifiera au groupe des suites $(k_j,l_j)_n$
dans $\M Z^2$, v\'erifiant pour tout $j$ assez grand la relation
$\pmat(k_{j+1},l_{j+1})=A_j\pmat(k_{j},l_{j})$, modulo les suites
nulles \`a partir d'un certain rang.
\medskip

On se donne $(\varepsilon_j)$ un suite positive sommable  strictement
d\'ecroissante, avec $\varepsilon_1<\frac12$.
On construit par r\'ecurrence deux suites d'entiers strictement 
croissantes, $(d_j)$ et  $(n_j)$ de la fa\c con suivante~:
on pose  $d_0=q_0=1$ et $ n_0=0$, puis pour tout $j \geq 1$ on
choisit $d_j$ et $n_j$ v\'erifiant
\begin{equation}\label{majd-1} 
d_{j-1}q_{n_{j-1}}\|d_j\alpha\|\leq \varepsilon_j,
\end{equation} 
\begin{equation}\label{majd}
\min{(q_{n_j},\varepsilon_ja_{n_j+1})}\geq d_j.
\end{equation}
On pose  pour tout $j\geq 1 $, $c_j=[d_j\alpha]$ puis
\begin{equation}\label{matrice}
A_j=\begin{pmatrix}
1+c_jq_{n_j}&-d_jq_{n_j}\\{c_jp_{n_j}}&1-d_jp_{n_j}
\end{pmatrix} \end{equation}
\`A tout couple d'entiers $( k',l')$ et tout $j_{0}\geq 1$ et on associe les
suites $(k'_j)$, $(l'_j)$ et $(t_j)$ d\'efinies \`a partir de indice
$j_0$ par~: pour tout $j \geq j_0$,  
$$\begin{pmatrix}k'_j\\l'_j\end{pmatrix}=
\prod_{i=j_0}^{j-1}A_i\;\begin{pmatrix}k'\\l'\end{pmatrix} \qquad
\hbox{et }\qquad t_j=k'_j\alpha-l'_j.$$ 
On pose pour tout $j\geq j_0$
$$b'_j= k'_jc_j-l'_jd_j$$ 
de sorte que
$k'_{j+1}=k'_j+b'_jq_{n_j}$, $l'_{j+1}=l'_j+b'_jp_{n_j}$ et
$t'_{j+1}=t'_j+b'_j\<q_{n_j}\alpha\>$ pour tout $j\geq j_0$.\\
On remarque qu'on a aussi
$$
b'_j=k'_j[d_j\alpha]-l'_jd_j=d_j(k'_j\alpha-l'_j)-k'_j\<d_j\alpha\>
=d_jt_j-k_j\<d_j\alpha\>.
$$

Le lemme suivant donne les ordres de grandeur ou vitesses de
convergence des suites ainsi d\'efinies~:
\begin{lemme}\label{taille} 
Pour tout couple $(k',l')$ donn\'e  il existe une constante $C'$
telle que, quel que soit $j_0\geq 1$,  pour tout   $j\geq j_0$
$$|t_j|\leq C'\qquad \hbox{et }\qquad \max{(|k'_j|,|l'_j|)}\leq
2C'd_{j-1}q_{n_{j-1}}.$$
La suite $(t_j)_j$ converge vers un r\'eel $t\in H_1(\alpha)$, et on
a les in\'egalit\'es suivantes pour tout $j\geq j_0$~:
$$|b'_j|\leq 2C'd_j\qquad \hbox{ et } \qquad |t-t_j|\leq
2C'\varepsilon_j\alpha_{n_j-1}.$$
\end{lemme}

\dem
Pour tout $j \geq 1$ on pose
$C_j=(|k'|+|l'|)\prod_{i<j}(1+\varepsilon_i)$, qui est major\'ee par une
constante $C'$ par hypoth\`ese.  Soit $j_0\geq 1$, on montre par
r\'ecurrence que pour tout $j\geq j_0$, on a $$|t_j|\leq C_j \qquad
\hbox{ et } \qquad \max{(|k'_j|,|l'_j|)}\leq 2C_jd_{j-1}q_{n_{j-1}}.$$
Les in\'egalit\'es sont v\'erifi\'ees pour $j=j_0$.  Supposons
l'hypoth\`ese vraie au rang $j$, on a alors~:
$$t_{j+1}=t_j+b'_j\<q_{n_j}\alpha\>
=t_j(1+d_j\<q_{n_j}\alpha\>)-k'_j\<d_j\alpha\>\<q_{n_j}\alpha\>.$$
Avec l'hypoth\`ese de r\'ecurrence et les in\'egalit\'es
(\ref{majd-1}) et (\ref{majd}) on obtient alors
$$|t_{j+1}|\leq C_j (1+\varepsilon_j\alpha_{n_j}(a_{n_j+1}+2)).$$
Comme $q_{n_1}>1$ par construction, on voit facilement que
$\alpha_{n_j}(a_{n_j+1}+2)\leq 2\alpha_{n_j}+\alpha_{n_j-1}\leq 1$,
d'o\`u $|t_{j+1}|\leq C_{j+1}$.\\
En ce qui concerne $k'_j$ (et le raisonnement est identique pour
$l'_j$), on a~:
$$k'_{j+1}=k'_j+ b'_jq_{n_j}= k'_j
+q_{n_j}(t_jd_j-k'_j\<d_j\alpha\>).$$
Comme pr\'ec\'edemment on obtient
\begin{eqnarray*}
|k'_{j+1}|& \leq & 2C_j(d_{j-1}q_{n_{j-1}}+\frac 1 2 q_{n_j}d_j +
q_{n_j}\varepsilon_j)\\
&\leq & 2C_jd_jq_{n_j}( \frac{q_{n_j-1}}{q_{n_j}}+\frac 1 2
+\varepsilon_j).
\end{eqnarray*}
Comme $a_{n_j+1}\geq 2$, on obtient $q_{n_j-1}<2q_{n_j}$ d'o\`u 
l'in\'egalit\'e cherch\'ee.\\
Les derni\`eres in\'egalit\'es se d\'eduisent imm\'ediatement de ce
qui pr\'ec\`ede. Pour $b'_j$ on a
$$|b'_j|\leq |t_j|d_j+|k'_j|\|d_j\alpha\|\leq
2C'd_j(1/2+\varepsilon_j)<2C'd_j.$$
Puisque $d_j/a_{n_j+1}<\varepsilon_j$, cela entra\^{\i}ne que $t_j$
converge vers une limite $t\in H_1(\alpha)$, et on a
$$|t-t_j|\leq \sum_j^\infty |b'_j|\alpha_{n_j}\leq
2C'\varepsilon_j\sum_j^\infty a_{n_j+1}\alpha_{n_j}, $$
d'o\`u l'in\'egalit\'e annonc\'ee.
\findem

Le th\'eor\`eme~\ref{identifier} r\'esulte de la proposition
suivante~:
\begin{prop}\label{Gamma}
Avec les notations pr\'ec\'edentes, pour tout irrationnel $\alpha$
\`a quotients partiels non born\'es, pour toutes les suites
$(\varepsilon_j)$, $(n_j)$, $(d_j)$  v\'erifiant les conditions
(\ref{majd-1}) et (\ref{majd}),  la limite inductive des groupes $\M
Z^2$ avec les matrices de transition $(A_j)$ donn\'ees par (\ref{matrice})
s'identifie \`a un sous-groupe $\Gamma_\alpha$ de $H_1(\alpha)$
dense dans $\M R$, via l'homomorphisme $\pi$ d\'efini par~: 
\begin{eqnarray*}
\lim_{ \rightarrow }(\M Z^2,A_{j}) &\longrightarrow& \Gamma_\alpha\\ 
(k'_j,l'_j)_j & \pi \atop \longmapsto & \lim_\infty (k'_j\alpha-l'_j).
\end{eqnarray*}
De plus pour tout $\beta>0$ dans $\Gamma_\alpha$, on a $\beta \notin
\M Z$, et le relev\'e $\tilde e(T_{\alpha,\beta})$ du groupe des
valeurs propres de la tour de Kakutani associ\'ee dans $\M R$,
$T_{\alpha,\beta}$ est \'egal \`a $\beta^{-1}\Gamma_\alpha$.
\end{prop} 

\dem
Nous  justifions d'abord  la d\'efinition de $\Gamma_\alpha$, ainsi
que ses propri\'et\'es. 
La convergence de $(k'_{n}\alpha-l'_{n})$ r\'esulte du lemme
\ref{taille}.
 La densit\'e de $\Gamma_\alpha$ dans $\M R$  r\'esulte de la
derni\`ere in\'egalit\'e du lemme~\ref{taille}~: en effet, comme
$\alpha$ est irrationnel, il suffit de montrer qu'on peut approcher
arbitrairement pr\`es un r\'eel de la forme $k'\alpha-l'$. En
construisant la suite $(t_{j})$ en partant de $k',l'$ et d'un indice
$j_0$ assez grand, comme $C'$ est ind\'ependante de $j_0$, on peut
obtenir $C'\varepsilon_{j_0}\alpha_{n_{j_0-1}}$  et donc
$|t-t_{j_{0}}|=|t-(k'\alpha-l')|$ arbitrairement petit.\\
Montrons ensuite l'injectivit\'e de $\pi$ et
l'absence d'entiers non nuls dans $\Gamma_\alpha$. Supposons avoir
trouv\'e une suite $(k'_j,l'_j)$ v\'erifiant $\lim(k'_j\alpha-l'_j)
\in \M Z$. En utilisant les notations du lemme on a pour tout
$j$ assez grand~:
$$k'_j\alpha+\sum_j^\infty b'_iq_{n_i}\alpha  = 0 \mod 1.$$  
Or, le lemme~\ref{taille} donne d'une part $|b'_j|\leq 2 C' d_j\leq
2C'\varepsilon_ja_{n_j+1}$ et d'autre part~: 
$$ |k'_j|\leq 2C' d_{j-1}q_{n_{j-1}}\leq
2C'\varepsilon_{j-1}a_{n_{j-1}+1}q_{n_{j-1}}<
2C'\varepsilon_{j-1}q_{n_j}.$$ 
Par cons\'equent, d\`es que $2C'\varepsilon_{j-1}<1/2$ on doit avoir
$b'_j=k'_j=0$  (cf lemme~\ref{unicite}).   On trouve donc  que 
$$0=b'_j=k'_jc_j-l'_jd_j=-l'_jd_j,$$ 
c'est \`a dire que $l'_j=0$ pour tout $j$ assez grand. On en conclut
que $\pi$ est bien injective et $\Gamma_\alpha \cap \M Z=\{0\}$.

Il reste maintenant \`a d\'emontrer la derni\`ere partie de la
proposition.
On choisit $\beta>0$ dans $\Gamma_\alpha$. Par d\'efinition c'est
l'image par $\pi$  d'une suite d'entiers $(k_j,l_j)$ telle que
$k_j \neq 0$ \`a partir d'un certain rang $J$, et on peut donc \'ecrire 
$\beta= \beta_J+\sum_J^\infty b_j\<q_{n_j}\alpha\>$ avec pour tout
$j\geq J$,  
$$b_j=k_jc_j-l_jd_j \qquad \hbox{ et }\qquad
\beta_j=k_j\alpha-l_j \neq 0.
$$
Il faut montrer $\tilde e(T_{\alpha,\beta})=\beta^{-1}\Gamma_\alpha$.   
Supposons que $t\in \Gamma_\alpha$ et posons $s=t/\beta$.  Avec les
notations pr\'ec\'edentes, pour que $(\beta,s,t)$ v\'erifie $(*)$ il
suffit de montrer que $\sum_1^\infty (b'_j-sb_j)^2<\infty$. On a
$b'_j=d_jt_j-k'_j\<d_j\alpha\>$, de m\^eme
$b_j=d_j\beta_j-k_j\<d_j\alpha\>$, et d'apr\`es le lemme~\ref{taille},
pour tout $j$ assez grand
\begin{eqnarray*}
|b'_j-b_js|&\leq  & |t_j-s\beta_j|d_j+|k'_j-sk_j|\|d_j\alpha\|\\
& \leq & d_j(|t-t_j|+|s||\beta_j-\beta|)+\O(\varepsilon_j) \\
&\leq & \O((d_j\alpha_{n_j-1}+1)\varepsilon_j).
\end{eqnarray*}
Comme $d_j\leq q_{n_j}$ on obtient $d_j\alpha_{n_j-1}\leq 1$, d'o\`u
le r\'esultat puisque $(\varepsilon_j)$ est sommable.
\\
On suppose maintenant que $s\in \tilde e(T_{\alpha,\beta})$ et
$t=s\beta$, donc que $(\beta, s, t)$ satisfait $(*)$.  On
\'ecrit $t=k'_1\alpha-l'_1+\sum_{j\geq 1}b'_j\<q_{n_j}\alpha\>$, avec
$b'_j=[b_js]$ pour tout $j$ assez grand,  et
$t_j=k'_j\alpha-l'_j=k'_1\alpha-l'_1+\sum_{1\leq i<j}
b'_i\<q_{n_i}\alpha\>$ pour tout $j\geq 1$. Alors
$k'_{j+1}=k'_j+b'_jq_j$, $l'_{j+1}=l'_j+b'_jp_j$ pour $j\geq 1$ et,
pour montrer que $t \in \Gamma_\alpha$, il suffit de v\'erifier
que  pour tout $j$ assez grand  
$$[b_js]=k'_jc_j-l'_jd_j.$$
On a toujours $b_j=d_j\beta_j-k_j\<d_j\alpha\>$ et aussi
$k'_jc_j-l'_jd_j=k'_j[d_j\alpha]-l'_jd_j=t_jd_j-k'_j\<d_j\alpha\>$
donc, avec $s_j=t_j/\beta_j$,
$$
b_js-(k'_jc_j-l'_jd_j)=d_j\beta_j(s-s_j)-(sk_j-k'_j)\<d_j\alpha\>.
$$
Comme pour $j$ assez grand $b'_j=[b_js]=\O(d_j)$ d'apr\`es le lemme
\ref{taille}, il r\'esulte encore de la r\'ecurrence sur les $k'_j$
que $k'_j=\O(d_{j-1}q_{j-1})$. D'autre part $s-s_j=o(\alpha_{n_j})$
d'apr\` es  le lemme  \ref{cv*}, donc
$$
b_js-(k'_jc_j-l'_jd_j)= 
\O(\varepsilon_ja_{n_j+1}\alpha_{n_j})+\O(d_{j-1}q_{n_{j-1}}\|d_j\alpha\|)
= \O(\varepsilon_j).
$$
On en d\'eduit que $b_js-(k'_jc_j-l'_jd_j)$ tend vers
$0$, d'o\`u l'\'egalit\'e cherch\'ee.
\findem

\rem On v\'erifie facilement que la r\'eunion de tous les groupes
$\Gamma_\alpha$ obtenus pour tous les choix de suites $(n_j)$ et
$(d_j)$ v\'erifiant (\ref{majd}) et (\ref{majd-1}) est non
d\'e\-nom\-bra\-ble.

\subsection{Tours de Kakutani et odom\`etres.}\label{solenoide}
 
Soit $(\delta_n)$ une suite d'entiers sup\'erieurs ou \'egaux \`a 2 et
$D_n=\prod_{j<n}\delta_j$. L' odom\`etre associ\'e \`a $(D_n)$ peut se
se d\'efinir comme l'action de la translation de 1 sur la limite
projective de la suite~: 
$$
\M Z/D_{1}\M Z\cdots \leftarrow \M Z/D_{n}\M Z
\leftarrow \M Z/D_{n+1} \M Z \leftarrow\cdots.
$$
Son groupe de valeurs propres est la r\'eunion des $\left(\frac
1{D_{n}}\M Z\right)/ \M Z$ pour $n\geq1$. La suspension de cette
translation est un flot de translations sur le sol\'eno\"\i de
associ\'e \`a la m\^eme suite $(\delta _{n})$, c'est \`a dire la
limite projective des groupes $\M T$ avec les homomorphismes de
transition $x \mapsto \delta _{n} x$. En effet, la suspension de
l'odom\`etre est un flot de translations sur un groupe compact dont le
groupe des valeurs propres est le relev\'e dans $\M R$ du groupe des
valeurs propres de l'odom\`etre, c'est-\`a-dire $\bigcup_{n\geq
1}\frac 1{D_{n}}\M Z $, qui est le groupe dual du sol\'eno\"\i de.
Plus pr\'ecis\'ement le plongement de $\bigcup_{n\geq 1}\frac
1{D_{n}}\M Z $ dans $\M R$ d\'efinit par dualit\'e le flot de
translations sur le sol\'eno\"\i de. Remarquons qu'un tel plongement
est unique \`a homoth\'etie pr\`es donc qu'il existe un unique flot de
translations sur le sol\'eno\"\i de, \`a homoth\'etie de temps pr\`es.

\medskip

Ce paragraphe montre  comment on peut obtenir des conjugaisons entre
tours de Kakutani et odom\`etres et donc, avec les m\^emes conditions, la
conjugaison de flots $\tau_{\alpha,\beta,\beta^{-1}}$ avec des flots
de translations sur des sol\'eno\"\i des.  L'id\'ee consiste \`a
remarquer que lorsque les matrices $(A_n)$ de (\ref{matrice}) sont de
rang 1, le groupe $\Gamma_{\alpha}$ est le dual d'un sol\'eno\" \i
de. Le r\'esultat suivant donne les constructions de  (iii) et (v) du
th\'eor\`eme~\ref{3i}, et par suspension les constructions de (iii)
du th\'eor\`eme~\ref{flot} pour $\tau_{\alpha,\beta,\beta^{-1}}$.  La
partie (ii) du th\'eor\`eme~\ref{reci} repose \'egalement 
sur cette construction.
\begin{theo}\label{odo}
Pour tout $\alpha$ tel que  $\inf_{q \neq 0} q^2\|q\alpha\|=0$,  il
existe un ensemble non d\'e\-nom\-bra\-ble dense de $\beta$ tels que la
tour de Kakutani $T_{\alpha,\beta}$ soit conjugu\'ee \`a un
odom\`etre. R\'eciproquement, tout odom\`etre est conjugu\'e \`a une
tour de Kakutani $T_{\alpha,\beta}$, pour des irrationnels $\alpha$
avec des propri\'et\'es d'approximation diophantiennes arbitrairement
bonnes, et pour un choix dense de $\beta$.
\end{theo}
 
\subsubsection{Construction d'odom\`etres.}  

On se donne $\alpha$ un irrationnel tel que $\inf_{q \neq 0}
q^2\|q\alpha\|=0$, et une suite $(\varepsilon_n)$ positive strictement
d\'ecroissante sommable.
On choisit une sous-suite strictement croissante $(q_{n_j})$ des
d\'enominateurs de la fraction continue de $\alpha$, v\'erifiant pour
tout $n\geq 1$
\begin{equation}\label{type1}
    q_{n_j}^2\alpha_{n_j}\leq \varepsilon_j/2.
\end{equation}
Les matrices $(A_j)$ d\'efinies par (\ref{matrice}) sont de rang
1 si
$$
\det(A_j)=\det\pmat(1+c_jq_{n_j}&-d_jq_{n_j},c_jp_{n_j}&1-d_jp_{n_j})
=1+c_jq_{n_j}-d_jp_{n_j}=0.
$$
Cette \'egalit\'e est satisfaite pour un couple d'entiers positifs
$(c_j,d_j)$ minimal, celui des coefficients de Bezout de $(q_{n_j},p_{n_j})$
c'est \`a dire, selon la parit\'e de $n_j$,
\\
\centerline{$ 
\left\lbrace  
\begin{array}{lc l}
c_j=p_{n_j-1} & \hbox{si $n_j$ est impair et } &p_{n_j}-
p_{n_j-1}\quad \hbox{sinon}\\
d_j=q_{n_j-1} & \hbox{si $n_j$ est impair et } &q_{n_j}-q_{n_j-1}
\quad \hbox{sinon}.
\end{array} 
\right.
$ } 
En particulier la condition (\ref{type1}) permet d'assurer la
validit\'e de la proposition~\ref{Gamma}. 
En ommettant les sous-suites par souci de lisibilit\'e,  on  \'ecrit maintenant
$A_n=\pmat(d_n,c_n)(p_n,-q_n)$. Alors  toute suite $\pmat(k_n,l_n)_n$ de la
limite inductive de $\M Z^2$ par $(A_n)$  v\'erifie pour $n$ assez
grand 
\begin{equation}\label{recrg1}
\Pmat(k_{n+1},l_{n+1})=A_n\Pmat(k_n,l_n)=(k_np_n-l_nq_n)\Pmat(d_n,c_n).
\end{equation}
Si on note pour tout $n \geq 2$, $u_{n}=\pmat(d_{n-1},c_{n-1})$,
alors $u_{n}$ est invariant par $A_{n-1}$ et on a  
\begin{equation}\label{u}
A_{n}u_{n}=\delta_{n}u_{n+1}\qquad \hbox{o\`u} \qquad
\delta_{n}=d_{n-1}p_{n}-c_{n-1}q_{n}.
\end{equation}
Dans ces conditions, le groupe $\Gamma_\alpha$ de la proposition
\ref{Gamma} est de rang  1, et  engendr\'e par des r\'eels
$(t^{(j)})_{j\geq 2}$ associ\'es aux  suites $(v_n^{(j)})_n$
d\'efinies pour tout $j\geq 2$ par  
$$v_{j-1}^{(j)}=v_{j}^{(j)}=u_j \quad \hbox{ et}\quad
v_{n+1}^{(j)}=A_nv_n^{(j)} \quad \hbox{ si }n\geq j.$$ 
Alors,  pour tout $j\geq 2$ donn\'e, on a
$v^{(j)}_{j+1}=A_ju_j=\delta_jv_{j+1}^{(j+1)}$, d'o\`u pour tout
$n\geq j+1$~:
\begin{equation}\label{v}
v^{(j)}_n=\delta_jv^{(j+1)}_n.
\end{equation}
Par cons\'equent on peut \'ecrire 
$\Gamma_\alpha= \cup_{j\geq 2}t^{(j)}\M Z$ avec $t^{(j)}=\delta
_{j}t^{(j+1)}$.

\medskip
Pour d\'efinir $\beta$,  on part d'un  vecteur  non nul  de $\M Z^2$,
$u_1=\pmat(k_1,l_1)$, tel que $k_{1}\alpha-l_{1}>0$ et on note $C'$
la constante du lemme~\ref{taille}~: comme $\Gamma_\alpha$  ne
d\'epend pas des premiers indices, on peut toujours quitte \`a
d\'ecaler d'un indice $n_0$ supposer que
\begin{equation}
4C'\varepsilon_1/q_1 <\|k_1\alpha\|.\label{beta>0}
\end{equation}
On d\'efinit alors pour tout $n\geq 1$,  $v_n^{(1)}=\pmat(k_n,l_n)$
par la r\'ecurrence $v_{n+1}^{(1)}=A_nv_n^{(1)}$ et on note
$\beta=\lim(k_n\alpha-l_n)$. La condition (\ref{beta>0}) et le lemme
\ref{taille} permettent d'assurer que $\beta>0$. En posant
$\delta_1=k_1p_1-l_1q_1$,  on obtient \`a l'aide de (\ref{u}) que
pour $ n\geq 1$  
$$v_n^{(1)}=D_nu_n \quad \hbox{avec } \quad  D_n=\prod_{1\leq j< n}
\delta_j. $$ 
En identifiant les coefficients avec (\ref{recrg1}), on trouve aussi
pour tout $n\geq 1$~:
 \begin{eqnarray*}
D_{n+1}&=&k_np_n-l_nq_n\\
&=&q_n\beta_n-k_n\<\alpha q_n\>\\
&=& q_n\beta +\O((\beta-\beta_n)q_n+k_n\alpha_n)\\
&=& q_{n}\beta +\O(4C'\varepsilon_n)\qquad  \hbox{( cf lemme
\ref{taille})}
\end{eqnarray*}  
Par cons\'equent, on peut \'ecrire
$$\Gamma_\alpha= \bigcup_{n\geq
1}\left(\frac{\beta}{D_{n+1}}\right)\M Z \quad \hbox{avec }
D_{n+1}=[\beta q_n] \hbox{ pour tout $n$ assez grand}. $$
Le groupe des valeurs propres de la tour de Kakutani
$T_{\alpha,\beta}$ est alors donn\'ee par  la proposition
\ref{Gamma}~:  
$$  e(T_{\alpha,\beta})=\cup_{n\geq 2} \frac{1}{D_n}\M Z \mod \M Z
.$$ 
Cette construction montre donc que le facteur de Kronecker de
$T_{\alpha,\beta}$ est bien l'odom\`etre $(\Omega,S)$, translation de
1 sur la limite projective de  la suite de groupes $(\M Z/D_n \M
Z)_n$~:
$$ \M Z/D_{1}\M Z\cdots \leftarrow \M Z/D_{n}\M Z \leftarrow \M
Z/D_{n+1} \M Z \ldots.$$
 
\noindent{\bf Remarque~:}
 Il est clair par cette construction que l'ensemble des $\beta$
obtenu est non d\'e\-nom\-bra\-ble (en modifiant les choix des suites
$(q_n)$) et dense dans $\M R^+$ (en d\'ecalant, pour une condition
initiale donn\'ee \`a l'avance,  les indices des suites v\'erifiant
(\ref{type1}) ).
  
\subsubsection{Injectivit\'e des fonctions propres sur l'odom\`etre.}

On cherche maintenant \`a montrer qu'une tour  de Kakutani,
$T_{\alpha,\beta}$, obtenue par la construction pr\'ec\'edente  est
conjugu\'ee \`a  l'odom\`etre. Notons $F$ la fonction  qui
semi-conjugue $(X_\beta,  T_{\alpha,\beta})$ \`a son facteur
Kronecker $(\Omega,S)$. 
 Pour montrer la conjugaison, il suffit de prouver l'injectivit\'e en
mesure de $F$~: autrement dit, il suffit de voir que sur un ensemble
de mesure pleine, 2 points distincts  ont des images distinctes par
$F$.  Comme $\Omega$ est la limite projective de $(\M Z/D_n \M Z)_n$
on peut \'ecrire $F=(F_n)_n$, o\`u pour tout $n\geq 1$, $F_n$ est, 
\`a  un facteur de normalisation pr\`es,  la fonction propre de
$T_{\alpha,\beta}$ associ\'ee \`a $s_n=1/D_n$. 
Il s'agit alors de trouver un ensemble de mesure pleine $X'_\beta$
tel que pour tout couple de points distincts dans $X'_\beta$ il
existe $n\geq 1$  pour lequel $F_n(x)\neq F_n(x')$. On doit donc
\'etudier le comportement des fonctions propres de
$T_{\alpha,\beta}$~: comme on sait que celles-ci sont li\'ees aux
solutions $\tilde f$ \`a l'\'equation
$$s\phi_\beta=t+\tilde f-\tilde f\circ T\mod 1, \qquad \hbox{(voir
chapitre~\ref{kaku})}$$ 
on va d'abord  d\'ecrire celles-ci sur $\M T$, lorsque $s$ est une
valeur propre de $T_{\alpha,\beta}$ et $t=s\beta$. Comme $s$ est
rationnel, en reprenant les notations  des paragraphes
\ref{transfert} et \ref{csidee}, on peut \'ecrire $\tilde f$ comme
une limite $L^1$ des fonctions $(\tilde f_n)$, soit en posant
$\tilde \theta_n=\tilde f_{n+1}-\tilde f_n \mod 1$, pour tout $n\geq
1$~: 
$$\tilde f=\tilde f_n +\sum_n^\infty \tilde \theta_m \mod 1.$$
Supposons maintenant $n\geq 1$  fix\'e, on pose $s=1/D_{n+1}$,  et on
reprend les notations habituelles~: pour tout $m \geq n$ \\
$\begin{array}{ll} 
\displaystyle \pmat(k_m,l_m)=v_m^{(1)},&\beta_n=k_n\alpha-l_n,\\
\displaystyle \pmat(k'_m,l'_m)= v_m^{(n+1)} , & t_m=k'_n\alpha-l'_n,
\quad \hbox{et }s_m=t_m/\beta_m.
\end{array}
$\\
 D'apr\`es (\ref{v}), on a en  fait  pour tout $m> n$, 
$$\pmat(k_m,l_m)=D_{n+1} \pmat(k'_m,l '_m) \quad \hbox{d'o\`u }\quad
\frac{b'_m}{ b_m}=\frac{t_m}{\beta_m}= s_m=s=1/D_{n+1}.$$ 
A l'aide du paragraphe \ref{idee}, on sait d\'ecrire l'allure de
$\tilde f$ sur la tour majeure d'ordre $n$, $(T^jB_{n})_{0\leq
j<q_{n}}$, associ\'ee \`a $\alpha$ (voir figure~\ref{2tours})~:
on sait, par \ref{transfert}, que  $\tilde f_n$  et $\tilde \theta_n$
sont  affines sur les \'etages de la partie principale de la tour
majeure d'ordre $n$. On en d\'eduit, comme $\tilde f_{n+1}=\tilde
f_n+\tilde\theta_n$ que $\tilde f_{n+1}$ est encore affine sur ces
\'etages et  de pente $sk_{n+1}-k'_{n+1}=0$. De plus,  si $m>n$,
$\tilde \theta_m$ est aussi affine par morceaux et de pente constante
$sb_m-b'_m=0$.  A l'aide de  (\ref{invariant}) et de la convergence
de la s\'erie $(\tilde \theta_m)$ on en d\'eduit que $ \tilde
\theta_m$ est nulle en dehors de $\I_m$ (voir lemme~\ref{bs-b'}). Par
cons\'equent, en dehors de $\J_n=\cup_{m\geq n}\I_m$ on a $\tilde
f=\tilde f_{n+1}$.  Il en r\'esulte que, en dehors de $\J_n$, $\tilde
f$ envoie chaque \'etage de la partie principale de la tour majeure
d'ordre $n$ sur une constante.\\
\begin{figure}
\begin{center}\caption{Tour d'ordre $n$ de
$T_{\alpha,\beta}$.}\label{tourodo} 
\input{tourodo.pstex_t}
\end{center}\end{figure}
Supposons pour l'instant que  $\beta>1$, et reprenons l'expression de
la fonction propre de $T_{\alpha,\beta}$ associ\'ee \`a $s$~: elle
s'\'ecrit pour tout $(x,y)\in X_\beta$ sous la forme  
$F_{n+1}(x,y)= sy+ \tilde f(x)\mod 1$ (cf  chapitre~\ref{kaku}).  Par
cons\'equent, en dehors du relev\'e $\tilde {\J_n} $ de   $\J_n$ dans
$X_\beta$, $F_{n+1}$ prend des valeurs constantes sur tous les
ensembles de la forme $T^jB_n\times \{k\}$ inclus dans $X_\beta$. \\
On peut pr\'eciser la nature de ces ensembles pour tout $n\geq 1$~:
comme $\beta$ n'appara\^ \i t que sur le dernier \'etage de la tour
majeure d'ordre $n$, ces ensembles constituent une tour de base
$\tilde B_n=B_n\times\{0\}$ associ\'ee \`a $T_{\alpha,\beta}$, pour
tous les indices $(j,k)$ tels que $0\leq k<\unde{[0,\beta[}(T^jB_n)$
si $0\leq j<q_n-1$ et $0\leq k<\lfloor \beta\rfloor$ si $j=q_n-1$
(voir figure~\ref{tourodo}).\\
Posons $\tilde q_{n-1}$ le d\'enominateur pr\'ec\'edent imm\'ediatement $q_n$ dans la suite des r\'eduites de $\alpha$ et $\tilde
\alpha_{n-1}=\|\tilde q_{n-1}\alpha\|$. 
Si on note $h_n$ la hauteur de cette tour,  en consid\'erant sa masse
totale   dans $X_\beta$ on obtient~:
$$\beta= \tilde
\alpha_{n-1}h_n+|\beta-\beta_n|+\O_1((\beta+1)\alpha_n\tilde
q_{n-1}).$$
D'apr\`es le lemme  \ref{taille}, on a
$|\beta-\beta_n|=\O(\varepsilon_n\tilde \alpha_{n-1})$. De plus avec
(\ref{type1}) on a aussi $1/\tilde\alpha_{n-1}=q_n+\O(\varepsilon_n)$,
ce qui donne finalement 
$$h_n=\beta q_n+\O(\varepsilon_n)=D_{n+1}.$$
Pour tout  $n$ assez grand, le relev\'e dans $X_\beta$ de la tour
majeure d'ordre $n$ pour la rotation d'angle $\alpha$ est donc une
tour associ\'ee \`a $T_{\alpha,\beta}$, de base $\tilde B_n$, et de
hauteur $D_{n+1}$. \\
Comme $F_{n+1}$ envoie, en dehors de $\tilde \J_n$, chacun de ces
$D_{n+1}$ \'etages sur une valeur dans $1/D_{n+1}\M Z$, et que
$F_{n+1}$  conjugue $T_{\alpha,\beta}$ \`a la translation de
$1/D_{n+1}$ sur cet ensemble, on en d\'eduit qu'elle est  injective,
en dehors de $\tilde \J_{n}$, sur les \'etages de la tour
$(T_{\alpha,\beta}^j\tilde B_n)_{0\leq j<D_{n+1}}$ .\\
Pour tout $n\geq 1$, la tour recouvre $X_{\beta}$ \`a un ensemble
$\Delta_{n}$ pr\`es, constitu\'e du relev\'e dans $X_\beta$ de la
tour mineure d'ordre $n$  (voir \ref{releve}) et de l'intervalle
$[\beta_{n},\beta]\times \{0\}$. Alors
$\lambda_{\beta}(\Delta_{n})\leq|\beta_{n}-\beta|+(\beta+1)\alpha_{n}q_{n}=\O(\varepsilon_{n}/q_{n})$
qui est  sommable. Par cons\'equent $\Delta=\lims \tilde \I_n\cup
\Delta_{n}$ est de mesure nulle.
Comme la suite des tours $(T_{\alpha,\beta}^j\tilde B_n)_{0\leq
j<D_{n+1}}$ engendre la tribu de $X_\beta$,  pour tout couple de
points distincts de $X_\beta \backslash \Delta$ on peut trouver $n$
tel que $x$ et $x'$ soient sur 2 \'etages disjoints de la tour
d'ordre $n$, en dehors de $\tilde \J_n$~: dans ces conditions, on a
bien $F_{n+1}(x)\neq F_{n+1}(x')$.

\medskip

Lorsque $0<\beta<1$, le raisonnement est analogue~: la fonction
propre associ\'ee \`a $s=1/D_{n+1}, F_{n+1}$, est la restriction \`a
$X_{\beta}=[0,\beta[$ de $\tilde f$. Elle est donc constante sur les
\'etages $T^jB_{n}$ inclus dans $[0,\beta[$, en dehors de $\J_{n}$.
Mais comme les it\'er\'es successifs de $B_{n}$ dans $X_{\beta}$ sont
exactement les \'etages de la tour pour l'induit, il suffit de
calculer la hauteur de celle-ci, $h_{n}$~: en d\'ecomposant
l'intervalle $[0,\beta[$ dans $\C D_{n}(\beta)$ (cf \ref{releve}), on
a encore $\beta=\tilde
\alpha_{n-1}h_{n}+\O_{1}(|\beta-\beta_{n}|+\alpha_{n}q_{n})$ et le
reste de la preuve est inchang\'e.

\subsubsection{Les odom\`etres comme tours de Kakutani.}

On se donne \`a pr\'esent un odom\`etre $(\Omega,S)$ d\'efini par
$$ \M Z/D_{1}\M Z\cdots \leftarrow \M Z/D_{n}\M Z \leftarrow \M
Z/D_{n+1} \M Z \ldots, $$
o\`u $(D_n)$ est une suite d'entiers donn\'ee telle que $D_n$ divise
$D_{n+1}$. Bien entendu, quitte \`a choisir  une sous-suite de
$(D_n)$, on peut choisir  une suite de diviseurs
$\delta_n=D_{n+1}/D_n$ aussi grande qu'on veut sans modifier
l'odom\`etre.\\
On veut construire une suite $(p_n/q_n)$ de r\'eduites et une suite
$(k_n,l_n)$ de sorte qu'on se trouve dans le cadre de la construction
pr\'ec\'edente~: il suffit pour cela de v\'erifier que si $(p_n/q_n)$
sont les r\'eduites d'un nombre $\alpha$, la condition (\ref{type1})
est satisfaite pour une suite $(\varepsilon_{n})$ sommable et que
$(k_n,l_n)$ est bien dans la limite inductive de la suite des groupes
$\M Z^2$ avec les matrices $(A_n)$ donn\'ees par la construction.
Ceci permet d'assurer que la tour de Kakutani ainsi construite est
isomorphe \`a un odom\`etre. Il faut ensuite voir qu'on peut choisir
les suites pr\'ec\'edentes de sorte que l'odom\`etre obtenu soit
exactement $\Omega$~: pour cela, il suffit d'avoir 
$$\delta_1=k_1p_1-l_1q_1 \qquad \hbox{et pour tout }n\geq 1 \quad
\delta_n=p_nd_{n-1}-q_nc_{n-1}.   $$
Notons $(a_n)$ la suite des quotients partiels associ\'ee \`a
$\alpha$~: on a vu que lorsque $n$ est impair, pour que $A_n$ soit de
rang 1 il faut que $d_n=\tilde q_{n-1}$ et $c_n=\tilde p_{n-1}$ o\`u
$\tilde p_{n-1}/\tilde q_{n-1}$ est la r\'eduite pr\'ec\'edent
imm\'ediatement $p_n/q_n$. Dans la construction suivante, on  prend
la suite compl\`ete des r\'eduites $(p_{n}/q_{n})$.  On obtient dans
ce cas que  pour tout  $n$ impair on a
$\delta_{n+1}=p_{n+1}q_{n-1}-p_{n-1}q_{n+1}=a_{n+1}$. Lorsque $n$ est
pair et non nul, on a alors $\delta_{n+1}=a_{n+1}+1$. \\
\smallskip

On construit maintenant $(a_n)$, et $(k_n,l_n)$ par r\'ecurrence~:
pour $n=0$ on d\'efinit comme d'habitude $p_0=a_0=0$ et $q_0=1$. Pour
$n=1$, on choisit $\delta_1$ et $a_1$ de sorte que $a_1$ ne divise pas
$\delta_1$. On pose $p_1=1$ et $q_1=a_1$ puis on choisit $(k_1,l_1)$
tels qu'on ait $\delta_1=k_1-l_1a_1$.\\
On choisit maintenant une suite $(\varepsilon_n)$ strictement
d\'ecroissante positive sommable de somme inf\'erieure \`a 1, de
sorte que $10(|k_1|+|l_1|)\varepsilon_1<1$~: si on note $C'$ la
constante du lemme~\ref{taille} associ\'ee \`a $(k_1,l_1)$,
l'in\'egalit\'e donne simplement $3C'\varepsilon_1<1$.\\
Supposons avoir d\'efini $a_j$, $(k_j,l_j)$ jusqu'\`a l'ordre $n$ de
sorte que (\ref{type1}) soit v\'erifi\'ee~:
pour d\'efinir $a_{n+1}$, il suffit de choisir $\delta_{n+1}$ tel que
$q_n<\varepsilon_n(\delta_{n+1}-1)$ et de prendre $a_{n+1}=\delta_{n+1}$
si $n$ est impair et $\delta_{n+1}-1$ sinon. 
$p_{n+1}$ et $q_{n+1}$ sont alors d\'efinies par la r\'ecurrence
classique des fractions continues, et on pose  bien s\^ur
$\pmat(k_{n+1},l_{n+1})=A_n\pmat(k_n,l_n)$.

\smallskip
On a bien d\'efini dans cette construction un irrationnel $\alpha$ et
une suite de  r\'eduites associ\'ee $(p_n/q_n)$ qui v\'erifient pour
tout $n\geq 1 $, $q_n\leq \varepsilon_n a_{n+1}$, ce qui entra\^\i ne
(\ref{type1}).   Posons $\beta=\lim_n (k_n\alpha-l_n)$,  pour montrer
que l'odom\`etre $\Omega$ est isomorphe \`a la tour de Kakutani
$T_{\alpha,\beta}$, il suffit de v\'erifier que $\beta>0$. 
D'apr\`es le lemme~\ref{taille}, on a
$\beta=\beta_1+\O_1(2C'\varepsilon_1/q_1)$~; d'autre part on avait 
$\delta_1=k_1p_1-l_1q_1$ d'o\`u 
$$\beta_1=k_1(q_1\alpha-p_1)/q_1+\delta_1/q_1=\delta_1/q_1+\O_1(k_1/(a_2q_1^2))
=\delta_1+\O_1(C'\varepsilon_1/q_1^3).$$
On obtient finalement la relation 
$\beta=\delta_1/a_1+\O_1(3C'\varepsilon_1/a_1)$. Comme on a
$3C'\varepsilon_1<1$ et qu'on a choisi $\delta_1/a_1$ positif, non
entier, on obtient bien par cette construction que $\beta>0$. Cette
derni\`ere relation montre \'evidemment la densit\'e des param\`etres
$\beta$ possibles. 

\subsection{Isomorphismes non triviaux avec des rotations
irrationnelles.}

Dans ce paragraphe, nous obtenons pour un choix particulier de la
suite $(d_j)$ d\'efinie au premier paragraphe  des exemples de tours
de Kakutani $T_{\alpha,\beta}$ conjugu\'ees \`a des rotations
irrationnelles. 
 Le r\'esultat suivant montre les parties (ii) et (iv) du
th\'eor\`eme~\ref{3i}, et par suspension le point (ii) du
th\'eor\`eme~\ref{flot} et (i) du th\'eor\`eme~\ref{reci}. 
\begin{theo}\label{rotation}
Pour tout $\alpha$  tel que $\inf_{q \neq 0} q^2\|q\alpha\|=0$, il
existe un ensemble non d\'e\-nom\-bra\-ble dense dans $\M R_+$ de $\beta$
tels que $T_{\alpha,\beta}$ soit conjugu\'e \`a une rotation
irrationnelle sur le cercle.  R\'eciproquement, toute rotation
irrationnelle du cercle  d'angle $s$  v\'erifiant $\inf_{q\neq 0}
q^2\|qs\|=0$ est conjugu\'ee \`a une  tour de Kakutani
$T_{\alpha,\beta}$, pour un choix non d\'e\-nom\-bra\-ble de param\`etres,
et dense en  $\beta$.
\end{theo}

Le paragraphe est compos\'e de trois sous-paragraphes~: le premier
donne, pour $\alpha$ donn\'e v\'erifiant $\inf_{q \neq 0}
q^2\|q\alpha\|=0$, une construction particuli\`ere de $\Gamma_\alpha$
qui est dans ce cas un sous-groupe de rang  2.  
Le second sous-paragraphe prouve la premi\`ere partie du th\'eor\`eme
\ref{rotation}, pour cette construction. Le troisi\`eme
sous-paragraphe \'etablit la partie r\'eciproque.

\subsubsection{Construction, propri\'et\'es diophantiennes de la
valeur propre.}\label{1vap}

Soit  $\alpha$ un irrationnel v\'erifiant $\inf_{q \neq 0}
q^2\|q\alpha\|=0$.
On garde, aux indices pr\`es, les notations pr\'ec\'edentes. Ici
$(\varepsilon_n)$ est une suite positive sommable strictement
d\'ecroissante dont la s\'erie des  restes converge. On choisit une
sous-suite $(p_n/q_n)$ de r\'eduites de la fraction continue de
$\alpha$ v\'erifiant  pour tout $n \geq 1$~:
$$q_{n}^2\alpha_{n}<\varepsilon_n$$ 
On pose $d_n=q_n$ (les relations (\ref{majd-1}) et  (\ref{majd}) sont
donc v\'erifi\'ees) et  $c_n=[d_n\alpha]=p_n$. La matrice de la
relation (\ref{matrice}) s'\'ecrit alors pour tout $n\geq 1$~:
$$A_n= \begin{pmatrix}
1+p_{n}q_{n} & -q_n^2 \\
p_n^2 & 1-p_nq_n 
\end{pmatrix}.$$
Dans ce cas particulier,  on a $\det(A_n)=1$ pour tout $n$~: la
limite inductive de $\M Z^2$ par la suite  $(A_n)$ est donc isomorphe
\`a $\M Z^2$. \\
Le groupe $\Gamma_\alpha$ de la proposition~\ref{Gamma} est alors de
rang 2 sur $\M Z$, et l'isomorphisme  de $\M Z^2$ dans
$\Gamma_\alpha$ est donn\'e par 
$(k'_1,l'_1)\mapsto t=\pi((k'_n,l'_n)_n) $.  Par cons\'equent, tout
choix d'entiers  $(k_1,l_1)$, $(k'_1,l'_1)$ tels que
$k'_1l_1-k_1l'_1=1$ donne  une base  $(\beta,t)$ de $\Gamma_\alpha$.\\
On choisit  $(k_1,l_1)$ un couple d'entiers premiers entre eux, on
associe la suite $(k_n,l_n)$ de la limite inductive de $\M Z^2$ par
$(A_n)$ et le r\'eel $\beta$ de $\Gamma_\alpha$ correspondant. Alors
on peut \'ecrire $\Gamma_\alpha=\beta\M Z+t\M Z$ o\`u $t$ est
l'\'el\'ement associ\'e \`a une suite $(k'_n,l'_n)$ issue d'un couple
$(k'_1,l'_1)$ v\'erifiant $k'_1l_1-k_1l'_1=1$.\\
Lorsque $\beta >0$, la proposition~\ref{Gamma} donne alors~:
$$e(T_{\alpha,\beta})=\beta^{-1}\Gamma_\alpha \mod \M Z=  t
\beta^{-1} \M Z \mod \M Z.$$
Notons $s =t/\beta$, le but des 2 sous-paragraphes suivants est
d'\'etablir une  conjugaison entre $T_{\alpha,\beta}$ et la rotation
sur $\M T$ d'angle $s$, qu'on notera $R_s$.  Nous donnons dans un
premier temps les propri\'et\'es diophantiennes de $s$. \\
On reprend les notations pr\'ec\'edentes~: les suites $(k_n,l_n)$,
$(k'_n,l'_n)$ et $(b_n, b'_n)$ sont donn\'ees pour tout $n \geq 1$
par les relations matricielles suivantes~:
\begin{equation}\label{defmat}
\Pmat(k_{n+1}&k'_{n+1},l_{n+1}&l'_{n+1})= A_n\Pmat(k_n&k'_n,l_n&l'_n),
\qquad  \Pmat( b_n,b'_n)=\Pmat(-l_n&k_n, -l'_n&k'_n)\Pmat(q_n,p_n) .
\end{equation}
 Enfin, $(\beta_n)$, $ (t_n)$ et $(s_n)$ sont d\'efinies par 
$$\beta_n=k_n\alpha-l_n,\qquad 
 t_n=k'_n\alpha-l'_n,\qquad s_n=t_n/\beta_n.$$
Notons  que $ b_n=\beta_nq_n-k_n<q_n\alpha>$ et que  pour $n$ assez
grand on a aussi~:
$$b_n=[\beta q_n] \qquad \hbox{  et} \qquad b'_n=[b_ns].$$
On rappelle les r\'esultats du lemme~\ref{cv*}~:
\begin{eqnarray}
 (t-t_n)=\O(\frac{\varepsilon_n}{q_n}),& &
 (\beta-\beta_n)=\O(\frac{\varepsilon_n}{q_n})\label{tb}\\
\max(k_n,k'_n) =\O(q_{n-1}^2), &\hbox{et} & (s-s_n) =
o(\alpha_n).\label{ks} 
\end{eqnarray}

\begin{lemme}\label{diophantien}
Avec les notations pr\'ec\'edentes, si $\alpha$ v\'erifie $\inf_{q
\neq 0} q^2\|q\alpha\|=0$,  alors le
groupe  $\Gamma_\alpha$  est engendr\'e par deux r\'eels $\beta$ et
$t$ rationnellement ind\'ependants. Pour  $s=t/\beta$, on a alors
$b_n=[q_n\beta]$ et $|b_ns-b'_n|\sim\alpha_n/\beta$. En particulier,
$(b'_n/b_n)$ est une sous-suite des r\'eduites de  $s$, qui satisfait
donc aussi $\inf_{q \neq 0} q^2\|qs\|=0$.
\end{lemme}
\dem
D'apr\`es la proposition~\ref{Gamma}, on sait que $(\beta,s,t)$
satisfait $(*)$~: par cons\'equent, $(b_ns-b'_n)$ converge vers $0$.
Comme $b_n=[\beta q_n]$ est arbitrairement grand, si $s$ \'etait
rationnel $b_n$ et $b'_n$ auraient n\'ecessairement un facteur
commun. Il suffit donc de voir que dans cette construction $b_n$ et
$b'_n$ sont premiers entre eux. Comme $A_n$ est de d\'eterminant $1$,
on obtient $\Det{ \mat(k'_n&k_n,l'_n&l_n)}= \Det{
\mat(k'_1&k_1,l'_1&l_1)}=1$.
D'autre part on a, d'apr\`es les relations pr\'ec\'edentes, 
$$\Pmat(b_n,b'_n)=\Pmat(k_n&-l_n,k'_n&-l'_n)\Pmat(p_n,q_n).$$ 
$(b_n,b'_n)$ est donc l'image par une isom\'etrie de $\M Z^2$ de
$(p_n,q_n)$ qui est irr\'eductible dans $\M Z^2$, donc $b_n$ et
$b'_n$ sont bien premiers entre eux.\\ 
Pour montrer la suite, on commence par une premi\`ere estimation de
$b_ns-b'_n$~: comme $s-s_n=o(\alpha_n)$, on a
\begin{eqnarray*}
b_ns-b'_n &= & b_n(s-s_n)+ \frac{1}{\beta_n}(b_nt_n-b'_n\beta_n)\\
& = & o(q_n\alpha_n)+ \frac{<\alpha
q_n>}{\beta_n}(k'_n\beta_n-k_nt_n)\\
& = & o(q_n\alpha_n) + \frac{<\alpha
q_n>}{\beta_n}(k_nl'_n-k'_nl_n))\\
& = & o(q_n\alpha_n)+ \frac{<\alpha q_n>}{\beta_n}.
\end{eqnarray*}
On obtient donc que $(b_ns-b'_n)=o(\varepsilon_n/q_n)$. Comme d'apr\`es
le lemme~\ref{taille}, $b_n=\O(q_n)$, on en d\'eduit que
$b_n|b_ns-b'_n|=o(\varepsilon_n)$, ce qui permet d'affirmer que
$(b'_n/b_n)$ sont bien des r\'eduites de $s$. Pour affiner notre
estimation, on r\'einjecte notre premi\`ere  relation dans
l'expression de $(s-s_n)$~:
$$
\beta_n(s-s_n)  =  \sum_{j\geq n} (b_js-b'_j)<\alpha q_j>
 =  o(\sum_{j\geq n} \varepsilon_j\alpha_j/q_j)
 =  o(\varepsilon_n \alpha_n/q_n).$$
Par cons\'equent, $b_n(s-s_n)=o(\varepsilon_n \alpha_n)$, et la
premi\`ere estimation devient donc bien $|b_ns-b'_n|\sim
\alpha_n/\beta$.\\ 
Comme $b_n= [\beta q_n]$,  il en r\'esulte imm\'ediatement que $b_n^2
\|b_ns\|=\O(q_n^2\alpha_n)$ tend vers $0$ lorsque $n \rightarrow
\infty$. 
\findem

\subsubsection{Tours de Kakutani conjugu\'ees \`a une rotation.}

Soit $\alpha$ un irrationnel donn\'e, v\'erifiant $\inf_{q \neq 0}
q^2\|q\alpha\|=0$. On note comme pr\'ec\'edemment  $\beta>0$ et $t$
des g\'en\'erateurs de $\Gamma_\alpha$ et on pose $s=t/\beta$. \\
Notons $R_s$ la rotation de $s$ sur $\M T$, et $f_\beta$ la fonction
propre de $T_{\alpha,\beta}$ associ\'ee \`a $s$.  Il s'agit de
montrer que $\tilde f_\beta$ est un isomorphisme de
$(X_\beta,\lambda_\beta, T_{\alpha,\beta})$ sur $(\M T ,\lambda,
R_s)$ (en notation additive).
Comme $s$ est irrationnel,  l'unique ergodicit\'e de $R_s$ montre que
la mesure image de $\lambda_\beta$ par $\tilde f_\beta $ est \'egale
\`a $\lambda$. $\tilde f_\beta$ est donc surjective sur un ensemble
de mesure de Lebesgue pleine, et  il reste \`a montrer que $\tilde
f_\beta$ est injective en dehors d'un ensemble de mesure 0. 
\\ 
On reprend les notations du paragraphe \ref{transfert}. Pour \'eviter
les confusions, on note $\tilde q_{n-1}$ le d\'enominateur
pr\'ec\'edent $q_n$ dans la suite des r\'eduites de $\alpha$ et
$\tilde \alpha_{n-1}=\|\tilde q_{n-1}\alpha\|$. Les conditions
diophantiennes sur $\alpha$ donnent l'approximation~: 
$$\tilde \alpha_{n-1}=1/q_n(1+\O(\varepsilon_n/q_n)).$$ 

\begin{figure}
\caption{Tour induite d'ordre $n$.}\label{tourbeta} 
\input{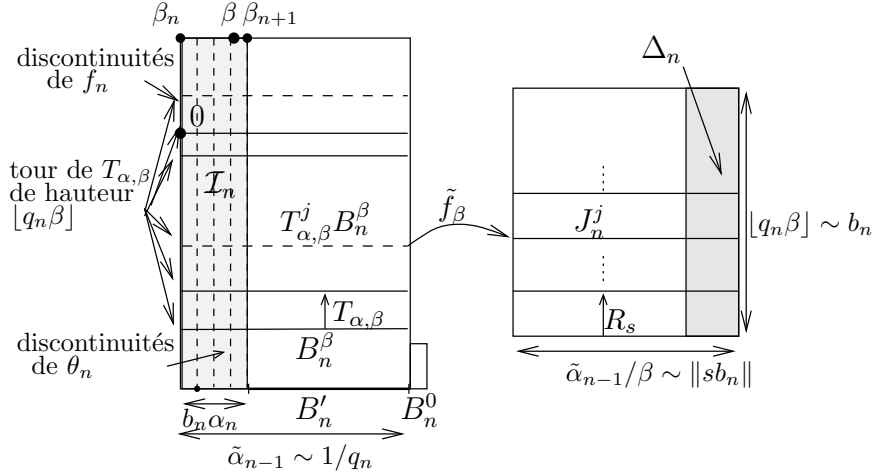}
\end{figure}
On rappelle que $f_\beta$ est donn\'ee \`a l'aide  de la fonction de
transfert $f$ o\`u $\tilde f$ est solution de $s\phi_\beta= -t+\tilde
f-\tilde f\rond T \mod 1$. On sait d'apr\`es le paragraphe  
\ref{cs} qu' en choisissant  $\tilde\theta_n=0$ au milieu de $B'_n$, 
la suite des fonctions $(f_n)$ converge vers $f$. Dans cette
construction, la convergence est tr\`es forte~: 
regardons  l'allure de $\tilde f_n$ sur la tour majeure associ\'ee
\`a 
$\C D_n(\beta)$ (cf figure~\ref{tourbeta}). 
 $\tilde f_n$ est affine sur chaque \'etage de la tour et de pente
$sk_n-k'_n$. \`A l'aide des estimations (\ref{tb}) et (\ref{ks}) on
obtient ici~:
$$
sk_n-k'_n = \frac{1}{\beta_n}(t_nk_n-\beta_nk'_n)+k_n(s-s_n) 
= \frac{1}{\beta_n}+o(k_n\alpha_n)
= \frac{1}{\beta}+\O(\varepsilon_n/q_n).
$$
Par cons\'equent, la pente de $\tilde f_n$ est de l'ordre de
$\frac{1}{\beta}$.\
D'autre part, en dehors de $\I_n$ (d\'efini selon (\ref{in})), on
sait que
$\tilde \theta_n$ est uniform\'ement petit (cf lemme~\ref{bs-b'}). 
Plus pr\'ecis\'ement ici on a, d'apr\`es (\ref{invariant}),  
$$\|\tilde \theta_n \unde{\I_n^c}\|_\infty\leq
q_n\alpha_n|b'_n-b_ns|+\|\tilde \theta_n\unde{B'_n\cup
B^0_n}\|_\infty.$$
Comme on a choisi $\tilde \theta_n=0$ au milieu de $B'_n$, et que
$\tilde
\theta_n$ est affine de pente $q_n(b_ns-b'_n)$ sur $B'_n \cup B_n^0$,
on
obtient que 
$$\|\tilde \theta_n \unde{\I_n^c}\| \leq  |b_ns-b'_n|
\leq   \O(\alpha_n).$$
De plus on a $\lambda(\I_n)=\alpha_n|b_n|q_n=\O(\varepsilon_n)$ dont la
s\'erie converge, d'o\`u $\lambda(\lims \I_n)=0$. On en d\'eduit que
$f_n$ converge 
vers $f$ en dehors de cet ensemble, et plus pr\'ecis\'ement, si
$\J_n= \cup_{m\geq n}\I_n$, on a 
$$\|(\tilde f-\tilde f_n)\unde{\J_n^c}\|_\infty\leq \sum_{m\geq
n}\alpha_m\leq \O(\alpha_n).$$ 
Il en r\'esulte, $\tilde f_n$ \'etant affine sur chaque \'etage de la
tour  majeure $T^jB_n$, que $\tilde f$ envoie $T^jB_n\cap \J_n^c$ sur
un ensemble inclus dans un intervalle de longueur $\tilde
\alpha_{n-1}(\frac{1}{\beta}+\O(\varepsilon_n))$ ($\O$ est bien s\^ur
uniforme pour $0\leq j<q_n$).\\
Revenons maintenant \`a $\tilde f_\beta$~:

Lorsque $\beta \in ]0,1[$ (cas de l'induit), $f_\beta$ est la
restriction de $\tilde f$ \`a l'intervalle $[0,\beta[$. Comme la tour
$\C D_n(\beta)$ est constitu\'ee d'intervalles ne contenant ni $0$ ni
$\beta$ (sauf le dernier \'etage), $[0,\beta[$ est \'egal \`a une
r\'eunion d'\'etages de $\C D_n$, \`a un ensemble de mesure $\tilde
\alpha_{n-1}$ pr\`es (le dernier \'etage). De plus la tour mineure
est de mesure $\alpha_n\tilde q_{n-1}\leq \varepsilon_n/q_n$, donc
$[0,\beta[$ est \'egal, \`a un ensemble de mesure $ \O(1/q_n)$
pr\`es, \`a une r\'eunion de $[\beta q_n]$ \'etages de la tour
majeure. Notons $B_n^\beta$ correspondant \`a l'\'etage le plus bas
dans la tour majeure, alors $(T_{\alpha,\beta}^jB_n^\beta)_{0\leq j<[
q_n\beta]}$ est une tour d'intervalles pour la transformation induite
$T_{\alpha,\beta}$. (cf figure~\ref{tourbeta}) \\
D'apr\`es ce qui pr\'ec\`ede,  on peut associer \`a chaque
$T_{\alpha,\beta}^jB_n^\beta$,  un intervalle $J_n^j$ de longueur
$\tilde\alpha_{n-1}(\frac{1}{\beta}+\O(\varepsilon_n))$ (uniforme en
$j$) contenant $\tilde f_\beta(T_{\alpha,\beta}^jB_n^\beta\cap
\J_n^c)$. Il reste \`a v\'erifier que l'ensemble de tous ces
intervalles  forme une partition, \`a un ensemble de mesure sommable
pr\`es~: comme la somme des mesures de ces ensembles est inf\'erieure
\`a $1+\O(\varepsilon_n)$, on minore la mesure de leur
r\'eunion~:
\begin{eqnarray*} 
 \lambda(\cup_{0\leq j<[\beta q_n]} J_n^j)&\geq &
\lambda_\beta( \cup_{0\leq j <[\beta
q_n]}T_{\alpha,\beta}^jB_n^\beta\cap \J_n^c)\\
&\geq & \sum_{j=0}^{[\beta
q_n]}\lambda_\beta(T_{\alpha,\beta}^jB_n^\beta)-\lambda_\beta(\J_n)\\
&\geq & 1-\O(1/q_n)-\O(\sum_n^\infty\varepsilon_m).
\end{eqnarray*}
Par cons\'equent, si on note $\Delta_n=\cup_j(J_n^j\cap \cup_{i\neq
j}J_n^i)$, les ensembles $(J_n^j\cap \Delta_n)_j$ sont  2 \`a 2
disjoints, et on a 
$$0\leq \lambda(\Delta_n)\leq  \sum_0^{[\beta
q_n]}\lambda(J_n^j)-\lambda(\cup_{o\leq j <[\beta q_n]} J_n^j)\leq
\O(\sum_n^\infty \varepsilon_m)+\O(1/q_n).$$
En dehors de $f_\beta^{-1}(\lims \Delta_n)\cup\lims \I_n \cup \lims
(\cup_0^{q_{n-1}}T^jB_n^0)$ qui est de mesure nulle (car
$\sum_n^\infty \varepsilon_m$ est sommable), l'injectivit\'e de
$f_\beta$ est imm\'ediate, ce qui montre la conjugaison avec $R_s$.
\\
Lorsque $\beta>1$, le raisonnement est identique au cas
pr\'ec\'edent, en remarquant comme dans la preuve du th\'eor\`eme
\ref{odo} que les relev\'es des tours majeures dans $X_\beta$ sont de
base $B_{n}\times\{0\}$ et de hauteur $h_{n}=[\beta q_{n}]$ (voir
figure~\ref{tourodo}).
\findem

\subsubsection{Des rotations comme tours de Kakutani.}

Ce paragraphe est consacr\'e \`a la partie r\'eciproque du
th\'eor\`eme~\ref{rotation}.\\
On se donne ici un irrationnel $s$  v\'erifiant $\inf_{q \neq 0}
q^2\|qs\|=0$.  Pour $(\varepsilon_n)$ une suite sommable strictement
d\'ecroissante, de restes sommables, on choisit, parmi les r\'eduites de $s$ 
deux suites strictement croissantes $(b_n)$ et
$(b'_n)$ v\'erifiant pour tout $n$
\begin{equation}\label{epsilon}
b_n^2|b_ns-b'_n|<\varepsilon_n.
\end{equation}
On cherche \`a construire deux r\'eels $\alpha$ et $\beta$ de sorte
que $(b_n)$ et $(b'_n)$ v\'erifient, avec les notations
pr\'ec\'edentes,  la construction du paragraphe \ref{1vap}. Dans  ce
cas, on pourra  appliquer la premi\`ere partie du th\'eor\`eme
\ref{rotation}, et assurer ainsi la conjugaison de $T_{\alpha,\beta}$
avec son facteur discret (celui-ci est bien la rotation d'angle $s$,
car $(b_ns-b'_n)$ tend vers $0$).  
Pour construire $\alpha$ et $\beta$, il nous suffit de construire une
sous-suite  $(p_n/q_n)$ des r\'eduites de $\alpha$ et les suites
$(k_n,l_n)$. Supposons pour l'instant avoir trouv\'e $\alpha$ et ses
r\'eduites~; d'apr\`es \ref{1vap}, sous la condition initiale que
$\pmat(k_1^{},l_1^{}),\pmat(k'_1,l'_1)$ soit une base (directe) de
$\M Z^2$,  les relations (\ref{defmat}) caract\'erisent les suites
$\pmat(k_n&k'_n,l_n&l'_n)$ et la suite $\pmat(b_n,b'_n)$.  Ces
relations sont \'equivalentes aux \'egalit\'es pour tout $n \geq 1$~:
$$ \Pmat(b_n,b'_n)=\Pmat(k_n&l_n,k'_n & l'_n)\Pmat(p_n,-q_n)\quad
\hbox{et
}\quad\Pmat(k_{n+1}&l_{n+1},k'_{n+1}&l'_{n+1})=\Pmat(k_n&l_n,k'_n&l'_n)+\Pmat(b_n,b'_n)(q_n,p_n).
$$
Comme par construction,  $A_n$ est une isom\'etrie directe,
$\pmat(k_n&l_n,k'_n & l'_n)$ est encore une isom\'etrie et on peut
donc inverser la premi\`ere relation. On obtient alors que les
conditions (\ref{defmat}) sont encore \'equivalentes \`a~: pour tout
$n\geq 1$
$$ \Pmat(q_n,p_n)=\Pmat(k'_n&k_n,l'_n& l_n)\Pmat(b_n,-b'_n) \qquad
\hbox{et }\qquad
\Pmat(k'_{n+1}&l'_{n+1},k_{n+1}&l_{n+1})=B_n\Pmat(k'_n&l'_n,k_n&l_n),$$
o\`u $B_n$ est la matrice d\'efinie par~:
$$B_n=\begin{pmatrix}1+b_nb'_n&-{b'_n}^2\\b_n^2& 1-b_nb'_n
\end{pmatrix}.$$
On reconnait \`a pr\'esent la construction du sous-paragraphe
\ref{1vap} pour l'irrationnel $1/s$~: comme $(b_n,b'_n)$ sont premier
entre eux et qu'on a \\
\centerline{${b'_n}^2|b'_n/s-b_n|=s
\left(\frac{b'_n}{sb_n}\right)^2|b_ns-b'_n|b_n^2=\O(\varepsilon_n),$}\\
on peut appliquer la proposition~\ref{Gamma} pour la suite de
matrices $(B_n)$. La limite inductive de $\M Z^2$ par $(B_n)$ est
donc un sous-groupe dense de $\M R$, de rang  2, et dont le quotient
de 2 g\'en\'erateurs quelconques est un nombre irrationnel, qu'on
note $\alpha$. De plus, d'apr\`es le lemme~\ref{diophantien},
$\alpha$ satisfait $\inf_{q \neq 0} q^2\|q\alpha\|=0$ , et
$(p_n/q_n)$ est bien une sous-suite des r\'eduites de $\alpha$ qui
v\'erifie
$q_n^2|q_n\alpha-p_n|=\O(\varepsilon_n)$. \\
Par cons\'equent, la suite $\pmat(k_n&k'_n,l_n&l'_n)$  construite \`a
l'aide de la suite $(B_n)$ v\'erifie aussi les relations
(\ref{defmat}) pour $\alpha$. En appliquant la premi\`ere partie du
th\'eor\`eme~\ref{rotation}, on a donc bien d\'efini un r\'eel
$\beta$ de sorte que $T_{\alpha,\beta}$  soit conjugu\'ee \`a la
rotation irrationnelle d'angle $s$. \\
La densit\'e du param\`etre $\beta$ r\'esulte de la densit\'e de
$\{k\alpha-l, \; (k,l)=1\}$  et du lemme~\ref{taille}~: il suffit de
remarquer que la construction pr\'ec\'edente peut se faire \`a partir
d'un rang $n_1$ choisi comme on veut~; comme d'apr\`es le lemme
\ref{taille} on a 
$|\beta-\beta_{n_1}|\leq 2C'\varepsilon_{n_1}$, et que $C'$ ne d\'epend
pas de $\alpha$, on peut choisir $n_1$ de sorte que le majorant soit
arbitrairement petit.\\ 
Expliquons maintenant pourquoi l'ensemble des param\`etres
$(\alpha,\beta)$ convenant est non d\'e\-nom\-bra\-ble~: pour un
irrationnel $s$ fix\'e, il est clair que le choix des sous-suites des
r\'eduites v\'erifiant l'hypoth\`ese (\ref{epsilon}) est non
d\'e\-nom\-bra\-ble. Revenons maintenant au paragraphe \ref{1vap}~:  avec
la construction de $\alpha$ \`a l'aide de $(B_n)$, on obtient des
d\'enominateurs $(q_n)$ qui s'\'ecrivent pour tout $n$ assez grand
sous la forme $q_n=[\beta_s b'_n]$, o\`u $\beta_s$ est un
g\'en\'erateur quelconque de $\Gamma_{1/s}$.  Comme l'hypoth\`ese
(\ref{epsilon}) impose que $b'_n=o(\tilde b'_{n+1})$ ($\tilde
b'_{n+1}$ est le d\'enominateur suivant $b'_n$ dans la suite des
r\'eduites de $1/s$), on en d\'eduit facilement que l'ensemble des
suites $(q_n)$ obtenues pour un $s$ donn\'e est non d\'e\-nom\-bra\-ble. Il
en r\'esulte que, ou bien l'ensemble des irrationnels $\alpha$ est
non d\'e\-nom\-bra\-ble, ou bien qu'on peut trouver un $\alpha$ pour
lesquel on sait construire un ensemble non d\'e\-nom\-bra\-ble de
r\'eduites associ\'ees~: dans ce dernier cas, l'unicit\'e de la
d\'ecomposition de $\beta$ le long de la suite $(q_n)$ montre que
l'ensemble des $\beta$ obtenus est non d\'e\-nom\-bra\-ble.
\findem

\section{R\'egularisation des flots sp\'eciaux.}

Dans ce chapitre nous montrons la proposition suivante, qui permettra
de terminer les preuves des th\'eor\`emes~\ref{flot} et \ref{reci}. 

\begin{prop}\label{regul}
Soit $\alpha\notin \M Q$, $(q_n)_n$ une sous-suite des
d\'enominateurs des r\'eduites de $\alpha$ et $(b_n)_n$ une suite
d'entiers non nuls.  Pour tout   $\beta \in \M R$ v\'erifiant
$\beta=\sum_1^\infty b_nq_n\alpha \mod 1$, on a les propri\'et\'es
suivantes~:
\begin{itemize}
\item Si     
$\displaystyle\sum_0^\infty
\left(|b_n|q_n^{k+1}\|q_n\alpha\|\right)^{1/(k+2)}<\infty, $  
alors $\phi_\beta$ est additivement cohomologue \`a une fonction
$\psi$ de $C^k(\M T)$ de norme infinie arbitrairement petite.
\item
Si $\displaystyle \sum_0^\infty \frac{q_n\ln (n|b_n|)}{-\ln
(\|q_n\alpha\|)}<\infty$,  alors la fonction $\phi_\beta$ est
additivement cohomologue \`a une fonction $\psi$ de $C^\omega(\M T)$
de norme infinie arbitrairement petite.
\end{itemize}
 \end{prop}

\subsection{Notations}

On rappelle ici les notations utiles dans ce paragraphe~: 
$\alpha$ est un irrationnel fix\'e et $T$ la rotation du cercle
d'angle $\alpha$. On note $(q_{n})$ une sous-suite des
d\'enominateurs des r\'eduites associ\'ees \`a $\alpha$ et
$\alpha_{n}=\|\alpha q_{n}\|$.  
Pour un r\'eel $\beta$  de $H_{1}(\alpha)$, on note $(b_n)$ une
suite d'entiers non nuls telle que 
$\beta=\sum_{n \geq 0}b_{n}\<q_{n}\alpha\>$, avec $\sum_n
|b_{n}|q_n\alpha_n<\infty$.
On d\'efinit $(k_n)_n$  par $k_n=\sum_0^{n-1}b_jq_j$ pour tout $n\geq
1$ et $k_0=0$. \\
D'apr\`es \ref{cobord} , on a en dehors de $\beta$ 
$$\phi_\beta=\sum_0^\infty \phi_n \qquad \hbox{o\`u }\quad
\phi_n=\phi_{ b_nq_n\alpha}\circ T^{k_n}.$$
On a vu (cf \ref{cobord}) que pour tout $n$,  $\phi_n$ est un cobord
additif, de fonction de transfert $\omega_{b_nq_n}\circ T^{k_n}$ qu'on notera
$h_n$. Alors on peut \'ecrire pour tout $n$, $\phi_n=h_n-h_n\circ T$
o\`u $h_n$ est une fonction affine par morceaux, de pente constante
\'egale \`a $b_nq_n$, admettant $|b_nq_n|$ sauts de taille 1, en des
points situ\'es dans $\I_n$ (cf figure~\ref{2tours}).
On note $\C T_n=\cup_{0\leq j<q_n}T^jB'_n$ la r\'eunion des \'etages
de la partie principale de la tour majeure d'ordre $n$ (cf paragraphe
\ref{transfert} et figure~\ref{2tours}). $h_{n}$ est d\'efinie
\`a une constante pr\`es, qu'on choisira  de sorte que $h_{n}$
s'annule sur la base $B'_{n}$.\\
On appelle maintenant $(K_n)$ une suite de noyaux continus, positifs
sym\'etriques, norm\'es dans $L^1(\M T)$. Lorsque la s\'erie
converge, on d\'efinit
$$\psi= \sum_0^\infty \phi_n\star K_n.$$
Nous pr\'ecisons dans ce qui suit les conditions g\'en\'erales sur
$(K_n)$ pour que $\psi$ soit r\'eguli\`ere, puis pour qu'elle soit
cohomologue \`a $\phi_\beta$. La synth\`ese de ces deux crit\`eres
permet ensuite d'obtenir la proposition~\ref{regul}.

\subsection{Crit\`eres de cohomologie et de r\'egularit\'e.} 
\begin{lemme}\label{cohom}
 Avec les notations pr\'ec\'edentes,  les fonctions $\psi$ et
$\phi_\beta$ sont cohomologues s'il existe une suite $(\varepsilon_n)$
strictement positive pour laquelle on ait 
\begin{eqnarray}
\label{cohom1}
\sum_{n}\varepsilon_nq_n&<&+\infty, \hbox{ et}\\
\label{cohom2}
\sum_n |b_n|\int_{|t|\geq \varepsilon_n} K_n(t)dt &<&+\infty.
\end{eqnarray}
\end{lemme}

\dem 
Comme on a $\phi_n=h_n-h_ncirc T$, il suffit de montrer que
$ \sum_{n}(h_n-h_n\star K_n)$ converge presque s\^urement.
Comme $h_n$ est affine sur les \'etages de $\C T_n$ et que $K_{n}$
est sym\'etrique, $h_n\star K_n$ est proche de $h_n$ sur une grande
proportion de $\C T_n$. \\
Si $\varepsilon_n\in ]0,\alpha_{n-1}/2[$ on appelle $B''_n$ le
sous-intervalle centr\'e de $B'_n$ de longueur
$\lambda(B''_n)=\lambda( B'_n)-2\varepsilon_n$, puis on note $\Delta_n$
la  sous-tour  de $\C T_n$ de base $B''_n$. Alors on a 
$$\lambda(\Delta_{n}^c)=
\varepsilon_{n}q_{n} +b_{n}\alpha_{n }q_{n}+\alpha_{n}q_{n-1} ,$$
qui est sommable d'apr\`es (\ref{cohom1}) et le choix de $\beta$
(dans $H_{1}(\alpha)$). Par cons\'equent  $\lims \Delta_{n}^c$ est de
mesure nulle.  
\\
Soit maintenant $x\in \Delta_n$ fix\'e, alors par
construction de $\Delta_n$, $h_n$ est affine sur un intervalle sym\'etrique autour de
$x$ de longueur au moins $\varepsilon_n$, et on peut  \'ecrire~:
\begin{eqnarray*}
h_n(x)-h_n\star K_n(x) &=& \int_{\M T} (h_n(x)-h_n(x-t))K_n(t)dt\\
&=& \int_{|t|\leq \varepsilon_n} b_nq_n t K_n(t)dt 
+\O(\|h_n\|_\infty \int_{|t|\geq \varepsilon_n}K_n(t)dt)\\
&=& \O(\|h_n\|_\infty)\int_{|t|\geq \varepsilon_n}K_n(t)dt. \quad ( K_n \hbox{
\'etant sym\'etrique).}
\end{eqnarray*}
On reste \`a estimer la norme infinie de $h_n$~: comme sur $B'_n\cup B_n^0$, 
$h_n$ est une fonction affine par
morceaux de pente $b_nq_n$ et
qu'elle s'annule sur $B'_{n}$, on obtient
$\|h_{n}\unde{B'_{n}\cup B_n^0}\|_{\infty}=\O(b_{n})$.  
Sur $B_{n}$, $h_{n}$ admet $b_{n}$ sauts de taille 1, ce qui donne encore
$\|h_{n}\unde{B_{n}}\|_{\infty}=\O(b_{n})$.
Enfin, la relation  $\phi_{n}=h_{n}-h_{n}\circ T$ montre que pour
tout $k \in\{1,..,q_{n}-1\}$ et  $x \in B_{n}$ on a l'\'egalit\'e 
$h_{n}(T^kx)=h_{n}(x)+kb_{n}\<\alpha q_{n}\>$, et que celle-ci est 
encore vraie pour $x \in B_n^0$ et $k\in \{1,..q_{n-1}-1\}$. Il en r\'esulte 
que
$\|h_{n}\|_{\infty}=\O(b_{n})$.  \\
On en d\'eduit avec l'hypoth\`ese (\ref{cohom2}) que pour tout
$x\notin \lims \Delta_{n}^c$, la s\'erie
$s(x)=\sum_{n}(h_{n}-h_{n}\star K_{n})(x)$ est convergente. La
fonction $s$ est donc bien d\'efinie presque partout et v\'erifie  la
relation $\phi_{\beta}-\psi=s-s\circ T$.  
\findem

\begin{lemme}\label{lisse}
Avec les notations pr\'ec\'edentes, pour tout $r \in \M N$,  $\psi$
est de classe $C^r$ si $K_n$ est de classe $C^r$ pour tout $n$ et si
on a pour tout $k\leq r$, 
\begin{equation}\label{ck}
\sum_n |b_n|\alpha_n \|K_n^{(k)}\|_\infty <+\infty.
\end{equation}
De plus, si $K_n$ est analytique pour tout $n$, alors $\psi$ est
analytique si 
\begin{equation}\label{ana}
\lim_{k\rightarrow +\infty }
\inv{k!}\sum_n|b_n|\alpha_n\|K_n^{(k)}\|_\infty =0.
\end{equation}
\end{lemme}

\dem 
La preuve est compl\`etement \'el\'ementaire~: il suffit d'estimer la
norme infinie des d\'eriv\'ees $k$-i\`eme de $\phi_n\star K_n$~:
remarquons  que $\phi_n=b_n\<q_n\alpha\>$ en dehors de $[\beta_n,
\beta_{n+1}]$ qui est un intervalle de longueur $|b_n|\alpha_n$. On
obtient donc pour tout $k$ tel que $K_n^{(k)}$ soit bien d\'efini, 
$$
\|(\phi_n\star K_n)^{(k)}\| \leq
\|\phi_{n}\|_{1}\|K_{n}^{(k)}\|_{\infty}\leq
2|b_n|\alpha_n\|K_n^{(k)}\|_\infty.$$
\findem

\subsection{Preuve de la proposition~\ref{regul}.}

\subsubsection{$\phi_\beta$ cohomologue \`a une fonction de  $C^r(\M
T)$.}
On suppose que $\beta=\sum_n b_n q_n\alpha \mod 1$ avec 

$$\sum_n(q_n^{r+1}|b_n|\alpha_n)^{1/(r+2)}<\infty.$$
Nous allons choisir convenablement la suite des noyaux $(K_n)$ pour
que les conditions (\ref{cohom1}), (\ref{cohom2}) et (\ref{ck}) des
lemmes \ref{lisse} et \ref{cohom} soient v\'erifi\'ees. 
Pour cela on se donne un noyau, $K$, positif sym\'etrique de classe
$C^r(\M R)$, nul en dehors de $]-1/2,1/2[$ et tel que $\|K\|_1=1$. Si
$(\delta_n)_n$ est une suite de $[0,1]$  de limite nulle, on
d\'efinit pour tout $x \in \M T $, $K_{n}(x)=K(\<x\>/\delta
_{n})/\delta _{n}$. Alors $K_n$ est un noyau positif sym\'etrique de
$C^r(\M T)$, nul en dehors de $]-\delta_n,\delta_n[$ et tel que
$\|K_n\|_1=1$. \\
Posons maintenant $\varepsilon_n=\delta_n$, les conditions
(\ref{cohom1}) et (\ref{cohom2}) sont v\'erifi\'ees 
 d\`es que 
$$\sum_{n\geq 1} \delta_n q_n<\infty.$$ 
Pour la r\'egularit\'e, pour tout $k\leq r$,  on a 
$\|K_{n}^{(k)}\|_{\infty}=\O(\delta _{n}^{-(k+1)})$.
Par cons\'equent, pour que (\ref{ck}) soit satisfaite il suffit que
$(\delta _{n})$ v\'erifie la condition~:
$$\sum_{n\geq 1 }|b_{n}|\alpha_{n}/\delta _{n}^{r+1}<\infty.$$
Choisissons maintenant $\delta_n=(|b_n|\alpha_n/q_{n})^{1/(r+2)}$ (de
sorte que $q_n\delta_n=|b_n|\alpha_n/\delta_n^{r+1}$)~: on trouve
comme condition suffisante~:
$$ \sum_{n\geq 1}(q_{n}^{r+1}|b_{n}|\alpha_{n})^{1/(r+2)}< \infty,$$
qui est l'hypoth\`ese de la premi\`ere partie de la proposition
\ref{regul}.
Pour montrer que $\psi$ peut \^etre choisie arbitrairement petite, il
suffit de remarquer que pour tout $n_0$, $\phi_{k_{n_0}\alpha}$ est
un cobord additif~: par cons\'equent $\phi_\beta$ est cohomologue \`a
$\sum_{n_0}^\infty \phi_n\star K_n$, pour tout $n_0$. Comme d'apr\`es
le lemme~\ref{lisse}, $(\|\phi_n\star K_n\|_\infty )$ est sommable,
pour tout $\varepsilon>0$, on peut toujours trouver $n_0$ de sorte que
$\|\sum_{n_0}^\infty \phi_n \star K_n\|_\infty <\varepsilon$.

\subsubsection{ $\phi_\beta$ cohomologue \`a une fonction
analytique.}
On suppose maintenant que $\beta=\sum_n b_n q_n\alpha \mod 1$ avec 
$$\sum_n \frac{q_n \ln (n|b_n|)}{-\ln(\alpha_n)}<+\infty.$$
Il s'agit de choisir une  suite $(K_n)$ de noyaux positifs
sym\'etriques analytiques  et norm\'es dans $L^1(\M T)$ qui
v\'erifient les conditions (\ref{cohom1}), (\ref{cohom2}) et
(\ref{ana}).
 On pose pour tout $n \geq 1$ et $t \in \M T$,  
$$K_n(t)={c_n \left(\sin(\pi m_n t) \over m_n\sin(\pi
t)\right)^{r_n}},$$
o\`u $r_n$ est un entier pair et $m_n$ un entier positif . $c_n$ est
donn\'e  de sorte que $\|K_n\|_1=1$.  On a~:
\begin{eqnarray*}
\|K_n\|_1 &\geq & c_n
\int_{1/2\geq |t|} \left(\frac{\sin(\pi m_n t)}{\pi m_n
t}\right)^{r_n}dt\\
&\geq &\frac{c_n}{m_n}
\int_{m_{n}/2\geq |t|} \left(\frac{\sin(\pi  t)}{\pi
t}\right)^{r_n}dt
\sim C \frac{c_n}{m_n\sqrt{r_n}}.
\end{eqnarray*}
Par cons\'equent $c_n =\O(m_n\sqrt{r_n})$. Les noyaux $(K_n)$ ainsi
choisis sont  bien positifs, sym\'etriques et norm\'es dans $L^1(\M
T)$. De plus, comme $K_n$ est proportionnel \`a $(F_{m_n})^{r_n/2}$,
o\`u $F_m$ est le noyau de Fejer d'ordre $m$, on a aussi
$\widehat K_n(j)=0$  pour tout $|j|\geq m_nr_n/2$ . En particulier
pour tout $k \in \M N$ on a 
$$\|K_n^{(k)}\|_\infty\leq (2\pi)^k\sum_{|j|< m_nr_n/2}j^{k}\leq
\frac{1}{k+1}(\pi m_{n}r_{n})^{k+1}.$$
Pour v\'erifier la condition (\ref{ana}), il suffit donc de choisir
$(m_n)$ et $(r_n)$ de sorte que 
$$\sum_n |b_n|\alpha_n \e{\pi m_nr_n}<+\infty.$$
Supposons $r_{n}$ donn\'e, et posons $m_n= \left\lfloor -\ln(\alpha_n)/(2\pi
r_{n})\right\rfloor$. On obtient alors
$$|b_n|\alpha_n \e{\pi m_nr_n}\leq |b_n|\alpha_n^{1/2}.$$
D'apr\`es l'hypoth\`ese, on a $\ln(|b_n|)=o(|\ln(\alpha_n)|)$, d'o\`u 
$ |b_n|\alpha_n^{1/2}=o(\frac{1}{-\ln(\alpha_{n})})$ qui est
sommable, 
et la condition (\ref{ana}) est donc satisfaite.\\
En ce qui concerne les conditions de cohomologie du lemme
\ref{cohom}, on doit choisir $(\varepsilon_n)$ le plus petit possible
(pour minimiser la condition (\ref{cohom1})) de sorte qu'on ait
(\ref{cohom2}). On a l'estimation suivante~:
\begin{eqnarray*}
\int_{1/2\geq |t|>\varepsilon_n}K_n(t)dt 
&\leq & 
c_n\int_{ |t|>\varepsilon_n} \left(\frac{1}{2 m_n t}\right)^{r_n}dt\\
&\leq &\O\left(
\frac{1}{\sqrt{r_n}}(2m_n\varepsilon_n)^{-r_n+1}\right).
\end{eqnarray*}
Posons $\varepsilon_n= \e{}/(2m_n)$, alors la condition (\ref{cohom2})
est satisfaite si $\sum_n |b_n|\e{-r_n}<+\infty$. Il suffit donc de
poser $r_n=2[\ln(n|b_n|)+1]$, ce qui donne 
$$ q_n\varepsilon_n=\frac{q_n\e{}}{2m_n}\sim  \frac{2\pi\e{} \,
q_n\ln(n|b_n|)}{-\ln(\alpha_n)}.$$
La condition (\ref{cohom1}) est donc v\'erifi\'ee par hypoth\`ese, ce
qui termine la preuve de la proposition~\ref{regul}. La norme infinie
des fonctions cohomologues \`a $\phi$ est arbitrairement proche de 0
pour les m\^emes raisons que dans le cas pr\'ec\'edent.

\subsection{Preuve des parties r\'eguli\`eres des th\'eor\`emes
\ref{flot} et \ref{reci}.}
 
D'apr\`es les chapitres pr\'ec\'edents, nous savons construire des
flots sp\'eciaux de fonction plafond $1+\gamma\phi_\beta$  avec les
propri\'et\'es (i), (ii) ou (iii) du th\'eor\`eme~\ref{flot}. Pour
montrer qu'on peut r\'egulariser cette fonction plafond  sans
modifier le type de la transformation, il suffit de montrer que
$\phi_{\beta}$ est cohomologue (additivement) \`a une fonction
r\'eguli\`ere, $\psi$ telle que $\|\psi\||\gamma|<1$. D'apr\`es la
proposition~\ref{regul}, cette derni\`ere contrainte n'en est pas
une. Il reste donc \`a v\'erifier dans chacune des constructions des
propri\'et\'es (i), (ii) et (iii), qu'on peut choisir les suites
$(b_n)$ et les sous-suites $(q_n)$ des d\'enominateurs des r\'eduites
de $\alpha$ de sorte que les conditions de la proposition~\ref{regul}
soient v\'erifi\'ees. Comme on a dans toutes les constructions
$b_n=\O(q_n)$, ces conditions conduisent bien \`a celles \'enonc\'ees
dans le th\'eor\`eme~\ref{flot}, et se v\'erifient facilement en
modifiant les conditions sur les $q_n$. 
\\
En ce qui concerne la preuve du th\'eor\`eme~\ref{reci}, (ii) se
v\'erifie facilement \`a l'aide du lemme~\ref{diophantien}, et (i)
s'obtient sans difficult\'e \`a partir de la preuve du th\'eor\`eme
\ref{odo}.

\section{Appendice. Les sous-groupes $H_\gamma(\alpha)$ des
approximations par les fractions continues.}

On rappelle d'abord un lemme classique, qui montre l'unicit\'e de la
d\'ecomposition d'Ostrowski~:
\begin{lemme}
Soit $k$ et $(b_n)$ des entiers v\'erifiant
$$k\alpha= \sum_{n=n_1}^\infty b_n\<q_n\alpha\> \mod 1.$$
Si $|k|<q_{n_1}$ et $|b_n|\leq a_{n+1}/2$ pour tout $n \geq n_1$,
alors on doit avoir $k=b_n=0$ pour tout $n\geq n_1$.
\end{lemme}
\dem
On montre d'abord que $k=0$~: supposons $k\neq 0$, comme
$|k|<q_{n_1}$ on aurait donc $\|k\alpha\|\geq \alpha_{n_1-1}$. Mais
on a aussi 
$$|\sum_{n=n_1}^\infty b_n\<q_n\alpha\>|\leq
\frac{1}{2}\sum_{n_1}^\infty a_{n+1}\alpha_n \leq \frac 1 2
(\alpha_{n_1-1}+\alpha_{n_1})<\alpha_{n_1-1},$$ 
ce qui est impossible. \\
Maintenant on peut \'ecrire $b_{n_1}q_{n_1}\alpha=
\sum_{n_1+1}^\infty b_n\<q_n\alpha\> \mod 1$, o\`u
$|b_{n_1}q_{n_1}|<a_{n_1+1}q_{n_1}<q_{n_1+1}$~: on obtient en
appliquant ce qui pr\'ec\`ede que $b_{n_1}=0$ puis que $b_n=0$ pour
tout $n \geq n_1$.
\findem 
 
 Il s'agit maintenant de d\'emontrer la proposition~\ref{hsat}, dont
on rappelle l'enonc\'e~:\\
{\bf Proposition~\ref{hsat}} 
{\it
Soit $\alpha \notin \M Q$, alors on a pour tout $\gamma>0$, \\
\centerline{$\begin{array}{lcl}
H_\gamma(\alpha)&=&\{x \in \M T, \; \sum_n \|q_nx\|^\gamma <+\infty\}
\; \hbox{et}\\
H_\infty(\alpha)&=&\{x \in \M T, \; \|q_nx\|\tend{n \rightarrow
\infty}0\}.
\end{array}
$}
}\\
\dem
Posons $\theta=(\sqrt 5-1)/2$, gr\^ace aux relations de r\'ecurrence
de $(q_{n})$ et $(\alpha_{n})$, on a imm\'ediatement les
in\'egalit\'es, pour tout $n\leq m$~:  
$$\alpha_{m}\leq \theta^{m-n-1}\alpha_{n}, \qquad q_{n}\leq
\theta^{m-n-1}q_{m}.$$
Supposons que $\beta=\sum_{n\geq 0}b_{n}\<q_{n}\alpha\>$. 
 On a dans ce cas $\|q_{n}\alpha_{m}\|=\|q_{m}\alpha_{n}\|\leq
\theta^{|m-n|-1}q_{m}\alpha_{m}$, d'o\`u on obtient 
 $$\|\beta q_{n}\|\leq \sum_{m=0}^\infty
|b_{m}|\|\alpha_{m}q_{n}\|\leq
\sum_{m}(|b_{m}|\alpha_{m}q_{m}\theta^{|m-n|-1}).$$
 Si $\sum_{m } (|b_{m}|\alpha_{m}q_{m})^\gamma<\infty$ avec
$\gamma\geq 1$, le terme de droite est la convolution d'une suite
sommable avec une suite de $l^\gamma(\M Z)$, ce qui montre bien  que
$\sum_{n}\|q_{n }\beta\|^\gamma <\infty$.
Si $\gamma <1$, on a $\|q_{n}\beta\|^\gamma\leq
\sum_{m}(|b_{m}|\|q_{n}\alpha_{m}\|)^\gamma$, d'o\`u on tire 
$$\sum_{n}\|q_{n}\beta\|^\gamma\leq
\sum_{m}(|b_{m}|q_{m}\alpha_{m})^\gamma
\sum_{n}\theta^{\gamma(||m-n|-1)}.$$
 Pour la r\'eciproque, on suppose maintenant que $\sum\|\beta
q_{n}\|^\gamma <\infty$. On note $k_{n}\in \{-q_{n}/2,..q_{n}/2\}$
tel que $\|\beta -k_{n}p_n/q_n\|=\min(|\beta-l/q_n|, l \in \M N)$.
Alors $\|\beta-k_{n}p_{n}/q_n\|=\|\beta q_{n}\|/q_n$ d'o\`u 
$$\|\beta-k_n\alpha\|\leq \|\beta- k_np_n/q_n\|+|k_n|\alpha_n/q_n
\leq  \|q_n\beta\|/q_n+\alpha_n/2.$$
En particulier $k_n\alpha \rightarrow \beta$. D'autre part on a aussi 
\begin{eqnarray*}
k_{n+1}p_{n+1}/q_{n+1}-k_np_n/q_n &=&
\frac{1}{q_{n+1}q_n}(k_{n+1}((-1)^n+p_nq_{n+1})-k_np_nq_{n+1})\\
& =&
(-1)^n\frac{k_{n+1}}{q_{n+1}q_n}+(k_{n+1}-k_n)\frac{p_n}{q_n}.
\end{eqnarray*}
On obtient que 
$$\|(k_{n+1}-k_n)p_n/q_n\|\leq
1/(2q_n)+\|q_{n+1}\beta\|/q_{n+1}+\|q_n\beta\|/q_n<1/q_n$$ pour tout
$n$ assez grand (car $\|q_n\beta\|\rightarrow 0$ par hypoth\`ese).
Par cons\'equent, on a bien $k_{n+1}=k_n+b_nq_n$ pour $n$ assez
grand, et $\beta =\sum_n b_n\<q_n\alpha\>$.
En r\'einjectant dans l'\'egalit\'e pr\'ec\'edente, on obtient
lorsque $b_n \neq 0$
\begin{eqnarray*}
\frac{|b_n|}{2a_{n+1}} &\leq&
\frac{|b_n+k_n/q_n|}{a_{n+1}}\\
&\leq
&\frac{q_{n+1}}{a_{n+1}}\left|\frac{k_n}{q_nq_{n+1}}+\frac{b_n}{q_{n+1}}\right|\\
&\leq &
2q_n(\|q_n\beta\|/q_n+\|q_{n+1}\beta\|/q_{n+1})\leq
2(\|q_n\beta\|+\|q_{n+1}\beta\|).
\end{eqnarray*}
 Ceci cl\^ot  la preuve de la proposition.
\findem
 
\rem
 Nous obtenons gr\^ace \`a ce r\'esultat une repr\'esentation
g\'eom\'etrique simple 
 de $H_\infty(\alpha)$ dans les tours associ\'ees aux fractions
continues.
 On reprend les notations du paragraphe \ref{betaqn}. Pour tout
$\varepsilon>0$, on appelle $\C V_n(\varepsilon)$ le rectangle de $\M
T\times \M R$ centr\'e en 0, 
$$\C V_n(\varepsilon)=
]-\varepsilon/q_n,\varepsilon/q_n[\times\{-h_n,..,h_n\}, \quad
\hbox{o\`u } h_n=\min(q_n/2,[\varepsilon/\alpha_n]).$$
Si l'on note indiff\'eremment un \'el\'ement $\beta \in \M T$ et ses
repr\'esentants dans un domaine fondamental d'ordre $n$, $\C
D_{n}(x)$, 
 la proposition~\ref{hsat} montre  qu'on a 
\begin{equation}\label{vbeta}
\|\beta q_n\|\tend{n \rightarrow \infty}0 \quad
\Longleftrightarrow\quad 
\hbox{pour tout }\varepsilon>0, \; \beta \in \limi \C
V_n(\varepsilon).
\end{equation}

\begin{figure}
\begin{center}
\caption{Repr\'esentation de $\C V_n(\varepsilon)$ dans un domaine
d'ordre $n$.}\label{vois0} 
\input{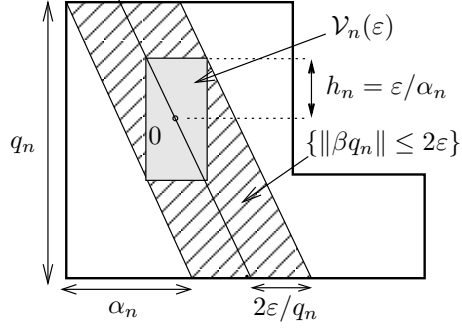}
\end{center}\end{figure}

En effet
 $\|\beta q_n\|\rightarrow 0$  signifie que $\beta \in H_\infty
(\alpha)$ d'o\`u $\beta= \sum b_{n}q_{n}\alpha$ avec
$b_{n}=o(a_{n+1})$. Par cons\'equent en notant
$k_n=\sum_0^{n-1}b_kq_k$, on a aussi $\beta=k_{n}\alpha+x_{n}$ o\`u 
$k_{n}=o(q_{n})$ et $x_{n}=o(1/q_{n})$. Il en r\'esulte que pour tout
$\varepsilon>0$ ( plus petit que 1/2), on a   pour tout $n$ assez
grand, $|k_{n}|\leq \varepsilon q_{n}\leq h_{n}$ et $|x_{n}|\leq
\varepsilon /q_{n}$ c'est \`a dire que $\beta \in 
\C V_{n}(\varepsilon)$.\\
Inversement, soit $\varepsilon>0$ et $\beta \in \C V_n(\varepsilon)$
pour tout $n$ assez grand. C'est \`a dire que $\beta=k_n\alpha+x_n$
o\`u $|x_n|\leq \varepsilon/q_n$ et $|k_n|\leq
\min(\varepsilon/\alpha_{n}, q_{n}/2)$. 
Alors $\|\beta q_n\|=\|k_n q_n\alpha+x_nq_n\|\leq
|k_n|\alpha_n+|x_n|q_n$. On a $|k_{n}|\alpha_{n}\leq
h_{n}\alpha_{n}\leq \min(\varepsilon, q_{n}\alpha_{n}/2)\leq
\varepsilon$ et donc $\|\beta q_{n}\|\leq 2\varepsilon$, ce qui
montre bien l'\'equivalence. Ces propri\'et\'es  peuvent se
repr\'esenter selon la figure~\ref{vois0}.\\ 

\bibliographystyle{plain}

%\bibliography{bibli}
 
\end{document}